\documentclass[reqno]{amsart}
\usepackage{amsmath,amssymb,amsthm,mathrsfs}

\newcommand{\BA}{{\mathbb A}}
\newcommand{\BC}{{\mathbb C}}\newcommand{\BD}{{\mathbb D}}

\newcommand{\BK}{{\mathbb K}}

\newcommand{\BP}{{\mathbb P}}

\newcommand{\BZ}{{\mathbb Z}}
\newcommand{\cD}{{\mathcal D}}
\newcommand{\cE}{{\mathcal E}}\newcommand{\cF}{{\mathcal F}}

\newcommand{\cK}{{\mathcal K}}\newcommand{\cL}{{\mathcal L}}

\newcommand{\cO}{{\mathcal O}}
\newcommand{\cR}{{\mathcal R}}

\newcommand{\cU}{{\mathcal U}}
\newcommand{\cW}{{\mathcal W}}\newcommand{\cX}{{\mathcal X}}
\newcommand{\cY}{{\mathcal Y}}\newcommand{\cZ}{{\mathcal Z}}
\newcommand{\bA}{{\mathbf A}}\newcommand{\bB}{{\mathbf B}}
\newcommand{\bC}{{\mathbf C}}\newcommand{\bD}{{\mathbf D}}
\newcommand{\bF}{{\mathbf F}}
\newcommand{\bG}{{\mathbf G}}

\newcommand{\bK}{{\mathbf K}}\newcommand{\bL}{{\mathbf L}}
\newcommand{\bM}{{\mathbf M}}

\newcommand{\bQ}{{\mathbf Q}}
\newcommand{\bS}{{\mathbf S}}

\newcommand{\bX}{{\mathbf X}}
\newcommand{\bZ}{{\mathbf Z}}
\newcommand{\wtilD}{\widetilde{D}}

\newcommand{\wtilG}{\widetilde{G}}

\newcommand{\wtilL}{\widetilde{L}}
\newcommand{\wtilN}{\widetilde{N}}

\newcommand{\wtilU}{\widetilde{U}}\newcommand{\wtilV}{\widetilde{V}}
\newcommand{\wtilW}{\widetilde{W}}\newcommand{\wtilX}{\widetilde{X}}
\newcommand{\wtilY}{\widetilde{Y}}
\newcommand{\al}{\alpha}

\newcommand{\ga}{\gamma}
\newcommand{\de}{\delta}\newcommand{\De}{\Delta}

\newcommand{\Th}{\Theta}

\newcommand{\la}{\lambda}\newcommand{\La}{\Lambda}

\newcommand{\si}{\sigma}\newcommand{\Si}{\Sigma}

\newcommand{\Om}{\Omega}
\newcommand{\mat}[2]{\ensuremath{\left[\begin{array}{#1}#2\end{array} \right]}}

\newcommand{\sbm}[1]{\left[\begin{smallmatrix} #1\end{smallmatrix}\right]}
\newcommand{\ov}[1]{{\overline{#1}}}

\newcommand{\tu}[1]{\textup{#1}}
\newcommand{\ands}{\quad\mbox{and}\quad}

\newcommand{\DS}{{\bf \Delta}}

\newcommand{\wtil}{\widetilde}

\newcommand{\bcD}{{\boldsymbol \cD}}
\newcommand{\boldcX}{{\boldsymbol \cX}}
\newcommand{\bdelta}{{\boldsymbol \delta}}

\newcommand{\btau}{{\boldsymbol \tau}}
\newcommand{\bSigma}{{\boldsymbol \Sigma}}

\newtheorem{theorem}{Theorem}[section]
\newtheorem{corollary}[theorem]{Corollary}
\newtheorem{lemma}[theorem]{Lemma}
\newtheorem{proposition}[theorem]{Proposition}

\newtheorem{remark}[theorem]{Remark}

\begin{document}

\title[Control and interpolation]{Robust control, multidimensional systems
and multivariable Nevanlinna-Pick interpolation}

\dedicatory{Dedicated to Israel Gohberg on the occasion of his 80th birthday}

\author{Joseph A. Ball}
\address{Department of Mathematics\\
Virginia Tech\\
Blacksburg VA 24061, USA} \email{joball@math.vt.edu}
\author{Sanne ter Horst}
\address{Department of Mathematics\\
Virginia Tech\\
Blacksburg VA 24061, USA} \email{terhorst@math.vt.edu}


\begin{abstract}
The connection between the standard $H^\infty$-problem in control
theory and Nevanlinna-Pick interpolation in operator theory was
established in the 1980s, and has led to a fruitful
cross-pollination between the two fields since. In the meantime,
research in $H^\infty$-control theory has moved on to the study of
robust control for systems with structured uncertainties and to
various types of multidimensional systems, while Nevanlinna-Pick
interpolation theory has moved on independently to a variety of
multivariable settings.  Here we
review these developments and indicate the precise
connections which survive in the more general
multidimensional/multivariable incarnations of the two theories.
\end{abstract}

\maketitle

%%%%%%%%%%%%%%%%%%%%%%%%%%%%%%%%%%%%%%%%%%%%%%%%%%%%%%%%%%%%%%%%%%%%%%%%%%%%%%%%%%%%%%%%%%%%%%%%
%%%%%%%%%%%%%%%%%%%%%%%%%%%%%%%%%%%%%%%%%%%%%%%%%%%%%%%%%%%%%%%%%%%%%%%%%%%%%%%%%%%%%%%%%%%%%%%%
\section{Introduction}

Starting in the early 1980s with the seminal paper
\cite{Zames} of George Zames, there occurred an active interaction
between operator theorists and control engineers in the development
of the early stages of the emerging theory of $H^{\infty}$-control.
The cornerstone for this interaction was the early recognition by
Francis-Helton-Zames \cite{FHZ} that the simplest case of the central
problem of $H^{\infty}$-control (the sensitivity minimization problem)
is one and the same as a Nevanlinna-Pick interpolation problem which had
already been solved in the early part of the twentieth century (see
\cite{Pick, NevPick}). For the standard problem of $H^{\infty}$-control
it was known early on that it could be brought to the so-called Model-Matching
form (see \cite{Doyle-notes, Francis}). In the simplest cases, the
Model-Matching problem converts easily to a Nevanlinna-Pick
interpolation problem of classical type.  Handling the more general
problems of $H^{\infty}$-control required extensions of the theory of
Nevanlinna-Pick interpolation to tangential (or directional)
interpolation conditions for matrix-valued functions; such extensions
of the interpolation theory were pursued by both engineers and
mathematicians (see e.g.~\cite{BGR, Dym, LA, Kimura87, Kimura89}).
Alternatively, the
Model-Matching problem can be viewed as a Sarason problem
which is suitable for application of Commutant Lifting theory
(see \cite{Sarason, FF}).  The approach of
\cite{Francis} used an additional conversion to a Nehari problem
where existing results on the solution of the Nehari problem in
state-space coordinates were applicable (see \cite{Glover, BR}).
The book of Francis  \cite{Francis} was the first book on $H^{\infty}$-control
and provides a good summary of the state of the subject in 1987.

While there was a lot of work emphasizing the connection of the
$H^{\infty}$-problem with interpolation and the related approach
through $J$-spectral factorization
(\cite{BGR,LA,LH,Kimura87,Kimura89,BR,BallCohen}), we should
point out that the final form of the $H^{\infty}$-theory parted ways
with the connection with Nevanlinna-Pick interpolation.  When
calculations were carried out in state-space coordinates, the
reduction to Model-Matching form via the Youla-Ku\v cera
parametrization of stabilizing controllers led to inflation of
state-space dimension; elimination
of non-minimal state-space nodes by finding pole-zero cancellations
demanded tedious brute-force calculations (see \cite{LA, LH}).
A direct solution in state-space coordinates (without reduction to
Model-Matching form and any explicit connection with Nevanlinna-Pick interpolation)
was finally obtained by Ball-Cohen \cite{BallCohen} (via a $J$-spectral
factorization approach) and in the more definitive coupled-Riccati-equation
form of Doyle-Glover-Khargonekar-Francis
\cite{DGKF}.
This latter paper emphasizes the parallels with older control
paradigms (e.g., the Linear-Quadratic-Gaussian and
Linear-Quadratic-Regulator problems) and obtained parallel formulas for the related
$H^{2}$-problem.
The $J$-spectral factorization approach was further
developed in the work of Kimura, Green, Glover, Limebeer, and Doyle
\cite{Kimura89, Green, GGLD}.  A good review of the state of the
theory to this point can be found in the books of Zhou-Doyle-Glover
\cite{ZDG} and Green-Limebeer \cite{GL}.

The coupled-Riccati-equation solution however has now been superseded
by  the Linear-Matrix-Inequality (LMI) solution which came shortly
thereafter; we mention specifically the papers of Iwasaki-Skelton \cite{IS} and
Gahinet-Apkarian \cite{GA}. This solution does not
require any boundary rank conditions entailed in all the earlier
approaches and generalizes in a straightforward way to
more general settings (to be discussed in more detail below).
The LMI form of the solution is particularly appealing from a
computational point of view due to the recent advances in semidefinite
programming (see \cite{GN}).
The book of Dullerud-Paganini \cite{DP} gives an up-to-date account of
these latest developments.

Research in $H^\infty$-control
has moved on in a number of different new directions, e.g., extensions of the
$H^{\infty}$-paradigm to sampled-data systems \cite{ChenFran},
nonlinear systems
\cite{vdS},
hybrid systems \cite{BallChu}, stochastic systems \cite{HinPrit},
quantum stochastic systems \cite{James}, linear repetitive processes
\cite{RogersGalOw},
as well as behavioral frameworks \cite{TrenWil}.  Our focus here will be on
the extensions to robust control for systems with structured
uncertainties and related $H^{\infty}$-control problems for
multidimensional ($N$-D) systems---both frequency-domain and state-space settings.
In the meantime, Nevanlinna-Pick interpolation theory has moved on to
a variety of multivariable settings (polydisk, ball, noncommutative
polydisk/ball); we mention in particular the papers
\cite{Agler-unpublished, DavP, Po1,
AgMcC99, BT, ArPo, BallBoloJFA, BallBoloNYJ, BallBoloJOT, BtH}.

As the transfer function for a multidimensional
system is a function of several variables, one would expect that the
same connections familiar from the 1-D/single-variable case should also
occur in these more general settings; however,  while there had been some
interaction between control theory and several-variable complex
function theory in the older area of systems over rings (see
\cite{KKT,KharSon, BST}), to this point, with a few exceptions
\cite{HeltonACC, HeltonIEEE, BM},
there has not been such an interaction in connection with
$H^{\infty}$-control for $N$-D systems and related such topics.
With this paper we wish to make precise the interconnections which do exist between
the $H^{\infty}$-theory and the interpolation theory in these more
general settings.  As we shall see, some aspects which are taken for
granted in the 1-D/single-variable case become much more subtle in
the $N$-D/multivariable case.  Along the way we
shall encounter a variety of topics that have gained attention
recently, and sometimes less recently, in the engineering
literature.

Besides the present Introduction, the paper consists of five sections
which we now describe:

\textbf{(1)} In Section \ref{S:1D} we lay out four specific results for the
classical 1-D case; these serve as models for the type of results
which we wish to generalize to the $N$-D/multivariable settings.

\textbf{(2)} In Section \ref{S:ring} we survey the recent results of Quadrat
\cite{Q-LN, Q-YKI, Q-Leuven, Q-elementary, Q-Lattice, Q-YKII} on
internal stabilization and parametrization of stabilizing controllers
in an abstract ring setting.  The main point here is that it is
possible to parametrize the set of all stabilizing controllers in
terms of a given stabilizing controller even in settings where the
given plant may not have a double coprime factorization---resolving
some issues left open in the book of Vidyasagar \cite{Vid}.
In the case where a double-coprime factorization is available, the
parametrization formula is more efficient.
Our modest new contribution here is to extend the ideas to the
setting of the standard problem of $H^{\infty}$-control (in the sense
of the book of Francis \cite{Francis}) where the
given plant is assumed to have distinct disturbance and control inputs
and distinct error and measurement outputs.

\textbf{(3)} In Section \ref{S:ND} we look at the
internal-stabilization/$H^{\infty}$-control problem for
multidimensional systems.   These
problems have been studied in a purely frequency-domain
framework (see \cite{Lin1, Lin2}) as well as in a state-space framework (see
\cite{Kaczorek, DuXie, DXZ}).  In Subsection \ref{S:ND-freq}, we
give the frequency-domain formulation of the problem.
When one takes the stable plants to consist of the ring of
structurally stable rational matrix functions, the general results of
Quadrat apply.  In particular,  for this setting
stabilizability of a given plant implies the existence of a double
coprime factorization (see \cite{Q-Leuven}). Application of the
Youla-Ku\v cera parametrization then leads to a Model-Matching form and, in
the presence of some boundary rank conditions, the
$H^{\infty}$-problem converts to a polydisk version of the
Nevanlinna-Pick interpolation problem.  Unlike the situation in the
classical single-variable case, this interpolation problem has no practical
necessary-and-sufficient solution criterion and in practice one
is satisfied with necessary and sufficient conditions for the
existence of a solution in the more restrictive Schur-Agler class
(see \cite{Agler-unpublished, AgMcC99, BT}).

In Subsection \ref{S:ND-statespace} we formulate the
internal-stabilization/$H^{\infty}$-control problem in
Givone-Roesser state-space coordinates.   We indicate the various
subtleties involved in implementing the state-space version
\cite{NJB,KharSon} of the
double-coprime factorization and associated Youla-Ku\v cera parametrization
of the set of stabilizing controllers.  With regard to the
$H^{\infty}$-control problem,
unlike the situation in the
classical 1-D case, there is no useable necessary and sufficient
analysis for solution of the problem; instead what is done (see e.g.
\cite{DuXie, DXZ}) is the use of an
LMI/Bounded-Real-Lemma analysis which provides a convenient set of
sufficient conditions for solution of the problem. This sufficiency
analysis in turn amounts to an $N$-D extension of the LMI solution
\cite{IS, GA} of the 1-D $H^{\infty}$-control
problem and can be viewed as a necessary and sufficient analysis of a
compromise problem (the ``scaled'' $H^{\infty}$-problem).

While stabilization and $H^{\infty}$-control problems have been
studied in the state-space setting \cite{Kaczorek, DuXie, DXZ} and in
the frequency-domain setting \cite{Lin1, Lin2} separately, there does
not seem to have been much work on the precise connections between
these two settings.  The main point of
Subsection \ref{S:ND-equiv} is to study this relationship; while
solving the state-space problem implies a solution of the
frequency-domain problem, the reverse direction is more subtle and
it seems that only partial results are known.  Here we introduce a
notion of {\em modal stabilizability} and {\em modal detectability}
(a modification of the notions of {\em modal controllability} and
{\em modal observability} introduced by Kung-Levy-Morf-Kailath
\cite{2D-II}) to obtain a
partial result on relating a solution of the frequency-domain problem
to a solution of the associated state-space problem.  This result
suffers from the same weakness as a corresponding result in
\cite{2D-II}: just as the authors in \cite{2D-II} were unable to
prove that {\em minimal} (i.e., simultaneously modally controllable
and modally observable) realizations for a given transfer matrix
exist, so also we are unable to prove that a simultaneously modally
stabilizable and modally detectable realization exists.
A basic difficulty in translating from frequency-domain to
state-space coordinates is the failure of the State-Space-Similarity
theorem and related Kalman state-space reduction for $N$-D systems.
Nevertheless, the result is a natural analogue of the corresponding
1-D result.

There is a parallel between the control-theory side and the
interpolation-theory side in that in both cases one is forced to be
satisfied with a compromise solution: the scaled-$H^{\infty}$ problem
on the control-theory side, and the Schur-Agler class (rather than
the Schur class) on the interpolation-theory side.  We include some
discussion on the extent to which these compromises are equivalent.

\textbf{(4)} In Section \ref{S:com-robust} we discuss several 1-D
variations on the internal-stabilization and $H^{\infty}$-control problem
which lead to versions of the $N$-D/multivariable problems discussed
in Section \ref{S:ND}.  It was observed early on that an
$H^{\infty}$-controller has good robustness properties, i.e., an
$H^{\infty}$-controller not only provides stability of the closed-loop
system associated with the given (or {\em nominal}) plant for which
the control was designed, but also for a whole neighborhood of plants
around the nominal plant.  This idea was refined in a number of
directions, e.g., robustness with respect to additive or
multiplicative plant uncertainty, or with respect to uncertainty in a
normalized coprime factorization of the plant (see \cite{McFGl}).
Another model for an uncertainty structure is the
Linear-Fractional-Transformation (LFT) model used by
Doyle and coworkers (see \cite{LZDProc, LZD}). Here a key concept is
the notion of structured singular value $\mu(A)$ for a finite square
matrix $A$ introduced by Doyle and Safonov \cite{DoyleMu,
SafonovMu}
which simultaneously generalizes the norm and the spectral radius
depending on the choice of uncertainty structure (a $C^{*}$-algebra
of matrices with a prescribed block-diagonal structure); we refer to
\cite{PackardDoyle} for a comprehensive survey.
If one assumes that
the controller has on-line access to the uncertainty parameters one
is led to a gain-scheduling problem which can be identified as the
type of multidimensional control problem discussed in Section
\ref{S:ND-statespace}---see \cite{Packard, AG}; we survey this
material in Subsection \ref{S:GS-SSC}. In Subsection \ref{S:GS-freq}
we review the purely frequency-domain approach of Helton
\cite{HeltonACC, HeltonIEEE} toward gain-scheduling which leads to
the frequency-domain internal-stabilization/$H^{\infty}$-control
problem discussed in Section \ref{S:ND-freq}.  Finally, in Section
\ref{S:hybrid} we discuss a hybrid frequency-domain/state-space
model for structured uncertainty which leads to a generalization of
Nevanlinna-Pick interpolation for single-variable functions
where the constraint that the norm be uniformly bounded by 1 is replaced by
the constraint that the $\mu$-singular value be uniformly bounded by 1;
this approach has only been analyzed for very special cases of the
control problem but does lead to interesting new results for operator
theory and complex geometry in the work of Bercovici-Foias-Tannenbaum
\cite{BFT1, BFT2, BFT3, BFT4},
Agler-Young \cite{AY1, AY2, AY3, AY4,
AY5, AY6, AY7, AY8, AY9}, Huang-Marcantognini-Young \cite{HMY}, and Popescu \cite{Po2}.

\textbf{(5)}  The final Section \ref{S:NC} discusses an enhancement of
the LFT-model for structured uncertainty to allow dynamic
time-varying uncertainties.  If the controller is allowed to have on-line
access to these more general uncertainties, then the solution of the
internal-stabilization/$H^{\infty}$-control problem has a form
completely analogous to the classical 1-D case.  Roughly, this result
corresponds to the fact that, with this noncommutative enhanced
uncertainty structure, the a priori upper bound $\widehat \mu(\bA)$ for the
structured singular value $\mu(\bA)$ is actually equal to $\mu(\bA)$,
despite the fact that for non-enhanced structures, the gap
between $\mu$ and $\widehat \mu$ can  be arbitrarily large (see
\cite{Treil}).  In this precise
form, the result appears for the first time in the thesis of
Paganini \cite{Paganini}
but various versions of this type of result have also appeared
elsewhere (see \cite{BFKT, BFT5, FM00, MT, Shamma}).  We discuss this
enhanced noncommutative LFT-model in Subsection \ref{S:NC-statespace}.  In
Subsection \ref{S:NC-FD} we introduce a noncommutative
frequency-domain control problem in the spirit of Chapter 4 of the
thesis of Lu \cite{LuThesis}, where the underlying polydisk occurring in
Section \ref{S:ND-freq} is now replaced by the noncommutative
polydisk consisting of all $d$-tuples of contraction operators on a
fixed separable infinite-dimensional Hilbert space $\cK$ and the
space of $H^{\infty}$-functions is replaced by the space of scalar
multiples of the noncommutative Schur-Agler class introduced in
\cite{BGM2}.   Via an
adaptation of the Youla-Ku\v cera parametrization of stabilizing
controllers, the internal-stabilization/$H^{\infty}$-control problem
can be reduced to a Model-Matching form which has the interpretation
as a noncommutative Sarason interpolation problem. In the final Subsection
\ref{S:NCstate-freq}, we show how the noncommutative state-space
problem is exactly equivalent to the noncommutative frequency-domain
problem  and thereby obtain an analogue of the classical case which
is much more complete than for the commutative-variable case given in
Section \ref{S:ND-equiv}.  In particular, if the problem data are
given in terms of state-space coordinates, the noncommutative Sarason
problem can be solved as an application of the LMI solution of the
$H^{\infty}$-problem.
While there has been quite a bit
of recent activity on this kind of noncommutative function theory
(see e.g.~ \cite{A-KV, BallBoloJOT, HMcCV, KVV09, Po3, Po4}),
the noncommutative Sarason problem has to this point escaped
attention; in particular, it is not clear how the noncommutative
Nevanlinna-Pick interpolation problem studied in \cite{BallBoloJOT}
is connected with the noncommutative Sarason problem.

Finally we mention that each section ends with a ``Notes'' subsection
which discusses more specialized points and makes some additional connections with
existing literature.

\paragraph{Acknowledgement}
The authors thank Quanlei Fang and Gilbert Groenewald for the useful
discussions in an early stage of preparation of the present paper.
We also thank the two anonymous reviewers for their thorough readings of
the first version and constructive suggestions for the preparation of the
final version of this paper.

%%%%%%%%%%%%%%%%%%%%%%%%%%%%%%%%%%%%%%%%%%%%%%%%%%%%%%%%%%%%%%%%%%%%%%%%%%%%%%%%%%%%%%%%%%%%%%%%
%%%%%%%%%%%%%%%%%%%%%%%%%%%%%%%%%%%%%%%%%%%%%%%%%%%%%%%%%%%%%%%%%%%%%%%%%%%%%%%%%%%%%%%%%%%%%%%%
\section{The 1-D systems/single-variable case}\label{S:1D}\setcounter{equation}{0}

\begin{figure}[h]
\setlength{\unitlength}{0.09in}
\centering
\begin{picture}(20,14)
%%boxes
\put(8,2){\framebox(5,4){$K$}}
\put(8,8){\framebox(5,4){$G$}}
%%G to K
\put(14,4.2){$y$}
\put(13,9){\line(1,0){3}}
\put(16,9){\line(0,-1){5}}
\put(16,4){\vector(-1,0){3}}
%%K to G
\put(6,9.2){$u$}
\put(8,4){\line(-1,0){3}}
\put(5,4){\line(0,1){5}}
\put(5,9){\vector(1,0){3}}
%%output
\put(14,11.2){$z$}
\put(13,11){\vector(1,0){5}}
%%input
\put(6,11.2){$w$}
\put(3,11){\vector(1,0){5}}
%%tap1
\put(17,4.2){$v_2$}
\put(19,4){\vector(-1,0){3}}
%%tap2
\put(3,9.2){$v_1$}
\put(2,9){\vector(1,0){3}}
\end{picture}
\caption{Feedback with tap signals}
\label{fig:taps0}
\end{figure}
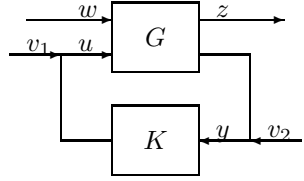

Let ${\mathbb C}[z]$ be the space of polynomials with complex
coefficients and ${\mathbb C}(z)$ the quotient field consisting of
rational functions in the variable $z$.  Let $\cR H^{\infty}$ be
the subring of {\em stable} elements of ${\mathbb C}(z)$ consisting
of those rational functions which are analytic and bounded on the
unit disk ${\mathbb D}$, i.e., with no poles in the closed unit disk
$\overline{\mathbb D}$. We assume to be given a plant
$G= \sbm{ G_{11} & G_{12} \\ G_{21} & G_{22}} \colon \cW \oplus \cU \to \cZ \oplus \cY$
which is given as a block matrix of appropriate size with entries from
${\mathbb C}(z)$. Here the spaces $\cU$, $\cW$, $\cZ$ and $\cY$ have the
interpretation of {\em control-signal space}, {\em disturbance-signal space}, {\em
error-signal space} and {\em measurement-signal space}, respectively, and
consist of column vectors of given sizes $n_{\cU}$, $n_{\cW}$, $n_{\cZ}$ and $n_{\cY}$,
respectively, with entries from ${\mathbb C}(z)$. For this plant $G$ we seek to design
a controller $K \colon \cY \to \cU$, also given as a matrix over ${\mathbb C}(z)$,
that stabilizes the feedback system $\Sigma(G,K)$ obtained from the signal-flow diagram
in Figure \ref{fig:taps0} in a sense to be defined precisely below.

\noindent
Note that the various matrix entries $G_{ij}$ of $G$ are themselves matrices
with entries from ${\mathbb C}(z)$ of compatible sizes (e.g., $G_{11}$ has size
$n_{\cZ} \times n_{\cW}$) and $K$ is a matrix over ${\mathbb C}(z)$ of size $n_{\cU} \times
n_{\cY}$.

The system equations associated with the signal-flow diagram of Figure \ref{fig:taps0}
can be written as
\begin{equation}  \label{standard-eqs}
    \begin{bmatrix} I & - G_{12} & 0 \\ 0 & I & -K \\ 0 & -G_{22} & I
    \end{bmatrix} \begin{bmatrix} z \\ u \\ y \end{bmatrix} =
    \begin{bmatrix} G_{11} & 0 & 0 \\ 0 & I & 0 \\ G_{21} & 0 & I
    \end{bmatrix} \begin{bmatrix} w \\ v_{1} \\ v_{2}
    \end{bmatrix}.
\end{equation}
Here $v_1$ and $v_2$ are tap signals used to detect stability properties of the internal
signals $u$ and $y$.
We say that the system $ \Sigma(G,K)$ is {\em well-posed} if there is a well-defined
map from $ \sbm{ w \\ v_{1} \\ v_{2} }$ to $\sbm{ z \\ u \\ y }$.
  It follows from a standard Schur complement computation that the
  system is well-posed if and only if $\det( I - G_{22}K) \ne 0$, and
  that in that case the map from $\sbm{ w \\ v_{1} \\ v_{2}}$ to
  $\sbm{z \\ u \\ y }$ is given by
$$
 \begin{bmatrix} z \\ u \\ y \end{bmatrix} = \Theta(G,K)
     \begin{bmatrix} w \\ v_{1} \\ v_{2} \end{bmatrix}
$$
where
\begin{align} & \Theta(G,K)
:= \begin{bmatrix} I & -G_{12} & 0 \\ 0 & I & -K \\ 0 & -G_{22} & I
  \end{bmatrix}^{-1} \begin{bmatrix} G_{11} & 0 & 0 \\ 0 & I & 0 \\
  G_{21} & 0 & I \end{bmatrix} = \notag \\
 &\small\begin{bmatrix} G_{11} + G_{12}K (I -G_{22}K)^{-1} G_{21} &
  G_{12}[ I + K (I - G_{22}K)^{-1}G_{22}]   & G_{12} K (I - G_{22} K)^{-1} \\
  K (I - G_{22}K)^{-1} G_{21} & I + K(I - G_{22}K)^{-1} G_{22} &
  K(I  - G_{22}K)^{-1}   \\
  (I - G_{22}K)^{-1} G_{21} & (I - G_{22}K)^{-1} G_{22} &
  (I - G_{22}K)^{-1} \end{bmatrix} \notag \\
  &=\small\begin{bmatrix}
  G_{11} + G_{12} (I - K G_{22})^{-1} K G_{21} &
  G_{12} (I - K  G_{22})^{-1}  & G_{12} (I - K G_{22})^{-1} K \\
(I - K G_{22})^{-1} K G_{21} & (I - K G_{22})^{-1} & (I - K G_{22})^{-1} K \\
[I + G_{22}(I - K G_{22})^{-1} K]G_{21} & G_{22}(I - K G_{22})^{-1} &
I + G_{22}(I - K G_{22})^{-1} K \end{bmatrix}.
\label{ThetaGK}
\end{align}
We say that the system $\Sigma(G, K)$ is {\em internally stable} if
$\Sigma(G,K)$ is well-posed and, in addition, if the map $\Theta(G, K)$
maps $\cR H^{\infty}_{\cW}\oplus \cR H^\infty_{\cU}
 \oplus  \cR H^{\infty}_{\cY}$ into
$ \cR H^{\infty}_{ \cZ} \oplus \cR H^{\infty}_{\cU}
  \oplus \cR H^{\infty}_{\cY}$, i.e., stable
inputs $w,v_{1}, v_{2}$ are mapped to stable outputs $z,u,y$. Note that
this is the same as the condition that the entries of $\Sigma(G,K)$
be in $ \cR H^{\infty}$.

We say that the system $\Sigma(G, K)$ {\em has performance} if
$\Sigma(G,K)$ is internally stable and in addition the transfer
function $T_{zw}$ from $w$ to $z$ has supremum-norm over the unit disk bounded
by some tolerance which we normalize to be equal to 1:
$$  \| T_{zw} \|_{\infty}: = \sup \{ \|T_{zw}(\lambda)\| \colon
\lambda \in {\mathbb D}\} \le 1.
$$
Here $\|T_{zw}(\lambda)\|$   refers to the induced operator norm, i.e., the
largest singular value for the matrix $T_{zw}(\lambda)$.
We say that the system $\Sigma(G,K)$ has {\em strict performance} if
in addition $\|T_{zw}\|_{\infty} < 1$.
The {\em stabilization problem} then is to describe all (if any
exist) internally stabilizing controllers $K$ for the given plant
$G$, i.e., all $K \in \BC(z)^{n_{\cU} \times n_{\cY}}$ so that the
associated closed-loop system $\Sigma(G, K)$ is internally stable.
The {\em standard $H^{\infty}$-problem} is to find all internally
stabilizing controllers which in addition achieve performance
$\|T_{zw}\|_{\infty} \le 1$. The strictly suboptimal $H^{\infty}$-problem
is to describe all internally stabilizing controllers which also achieve
strict performance $\| T_{zw}\|_{\infty} < 1$.

%%%%%%%%%%%%%%%%%%%%%%%%%%%%%%%%%%%%%%%%%%%%%%%%%%%%%%%%%%%%%%%%%%%%%%%%%%%%%%%%%%%%%%%%%%%%%%%%
\subsection{The model-matching problem}\label{S:MMproblem}

Let us now consider the special case where $G_{22}=0$, so that $G$ has
the form $G=\sbm{G_{11}&G_{12}\\G_{21}&0}$. In this case well-posedness
is automatic and $\Theta(G,K)$ simplifies to
\[
\Theta(G,K)=\mat{ccc}{G_{11}+G_{12}KG_{21}&G_{12}&G_{12}K\\KG_{21}&I&K\\G_{21}&0&I}.
\]
Thus internal stability for the closed-loop system $\Sigma(G,K)$ is equivalent to stability
of the four transfer matrices $G_{11}$, $G_{12}$, $G_{21}$ and $K$. Hence internal
stabilizability of $G$ is equivalent to stability of $G_{11}$, $G_{12}$ and $G_{21}$; when the
latter holds a given $K$ internally stabilizes $G$ if and only if $K$ itself is stable.

Now assume that $G_{11}$, $G_{12}$ and $G_{21}$ are stable. Then the $H^\infty$-performance
problem for $G$ consists of finding stable $K$ so that
$\|G_{11}+G_{12}KG_{21}\|_{\infty}\leq 1$.
Following the terminology of \cite{Francis}, the problem is called the Model-Matching Problem.
Due to the influence of the paper \cite{Sarason}, this problem is usually referred to as the
Sarason problem in the operator theory community; in \cite{Sarason} it is shown explicitly how
the problem can be reduced to an interpolation problem.

In general control problems the assumption that $G_{22}=0$ is an unnatural assumption. However,
after making a change of coordinates using the Youla-Ku\v cera parametrization or the
Quadrat parametrization, discussed below, it turns out that the general $H^\infty$-problem
can be reduced to a model-matching problem.

%%%%%%%%%%%%%%%%%%%%%%%%%%%%%%%%%%%%%%%%%%%%%%%%%%%%%%%%%%%%%%%%%%%%%%%%%%%%%%%%%%%%%%%%%%%%%%%%
\subsection{The frequency-domain stabilization and $H^\infty$ problem}

The following result on characterization of stabilizing controllers
is well known (see e.g.~\cite{Francis} or \cite{Vid, VSF} for a more
general setting).

\begin{theorem}  \label{T:1.1}  Suppose that we are given a rational
    matrix function $G
    = \sbm{G_{11} & G_{12} \\ G_{21} & G_{22}}$ of size $(n_{\cZ} +
    n_{\cY}) \times (n_{\cW} + n_{\cU})$ with entries in ${\mathbb
    C}(z)$ as above.  Assume that $G$ is stabilizable, i.e., there
    exists a rational matrix function $K$ of size $n_{\cU} \times
    n_{\cY}$ so that the nine transfer functions in \eqref{ThetaGK}
    are all stable.  Then a given rational matrix function $K$
    stabilizes $G$ if and only if $K$ stabilizes $G_{22}$, i.e.,
    $\Theta(G,K)$ in \eqref{ThetaGK} is stable if and only if
    \begin{align*}
     \Theta(G_{22},K): & = \begin{bmatrix} I + K (I - G_{22}K)^{-1} & K
    (I - G_{2}K)^{-1} \\ (I - G_{22}K)^{-1} & (I - G_{22}K)^{-1}
\end{bmatrix}  \\
       &  =
     \begin{bmatrix} (I - K G_{22})^{-1} & (I - K G_{22})^{-1} K \\
     G_{22}(I - K G_{22})^{-1} & I + G_{22}(I - K G_{22})^{-1} K
     \end{bmatrix}
   \end{align*}
   is stable.
Moreover, if we are given a double coprime factorization for
   $G_{22}$, i.e., stable transfer matrices
   $D$, $N$, $X$, $Y$, $\widetilde D$, $\widetilde N$, $\widetilde X$
   and $\widetilde Y$ so that the determinants of $D$, $\widetilde D$, $X$ and
   $\widetilde X$ are all nonzero (in $\cR H^{\infty}$) and
   \begin{equation}\label{1coprime}
   G_{22} = D^{-1} N = \widetilde N \widetilde D^{-1}, \quad
   \begin{bmatrix} D & -N \\ -\widetilde Y & \widetilde X
   \end{bmatrix} \begin{bmatrix} X & \widetilde N \\ Y &
   \widetilde D \end{bmatrix} =\mat{cc}{I_{n_{\cY}}&0\\0&I_{n_{\cU}}}
   \end{equation}
   (such double coprime factorizations always exists since $\cR H^{\infty}$ is
   a Principal Ideal Domain), then the set of all stabilizing controllers $K$ is given
   by either of the formulas
   \[
   K=(Y+\wtilD\La)(X+\wtilN\La)^{-1}=(\wtilX+\La N)^{-1}(\wtilY+\La D),
   \]
   where $\La$ is a free stable parameter from $\cR
   H^{\infty}_{\cL(\cU, \cY) }$ such that
   $\det (X+\wtilN\La)\not=0$ or equivalently  $\det(\wtilX+\La N)\not=0$.
   \end{theorem}

Through the characterization of the stabilizing controllers, those controllers that,
in addition, achieve performance can be obtained from the solutions of a Model-Matching/Sarason
interpolation problem.

   \begin{theorem}\label{T:1.2}
   Assume that $G\in {\mathbb C}(z)^{(n_{\cZ} +n_{\cY}) \times (n_{\cW} + n_{\cU})}$ is
   stabilizable and that
   $G_{22}$ admits a double coprime factorization \eqref{coprime}.
   Let $K\in {\mathbb C}(z)^{n_\cU \times n_\cY}$. Then $K$ is a solution to the standard
   $H^{\infty}$ problem for $G$ if and only if
   \[
   K=(Y+\wtilD\La)(X+\wtilN\La)^{-1}=(\wtilX+\La N)^{-1}(\wtilY+\La D),
   \]
   where $\La\in \cR H^{\infty}_{\cL(\cU, \cY)}$ so that $\det(X+\wtilN\La)\not=0$,
   or equivalently $\det(\wtilX+\La N)\not=0$, is any solution to the
   Model-Matching/Sarason interpolation
   problem for $\wtilG_{11}$, $\wtilG_{12}$ and $\wtilG_{21}$ defined by
   \[
   \wtilG_{11}:=G_{11}+G_{12} Y D G_{21},\quad
   \wtilG_{12}:=G_{12}\wtilD,\quad
   \wtilG_{21}:=DG_{21},
   \]
   i.e., so that
   \[
   \|\wtilG_{11}+\wtilG_{12}\La \wtilG_{21}\|_{\infty}\leq 1.
   \]
   \end{theorem}

   We note that in case $\widetilde G_{12}$ is injective and
   $\widetilde G_{21}$ is surjective on the unit circle, by absorbing
   outer factors into the free parameter $\Lambda$ we may assume
   without loss of generality that $\widetilde G_{12}$ is inner (i.e., $\widetilde G_{12}(z)$
   is isometric for $z$ on unit circle) and $\widetilde G_{21}$ is
   co-inner (i.e., $\widetilde G_{21}(z)$ is coisometric for $z$ on
   the unit circle).
   Let $\Gamma \colon L^{2}_{\cW} \ominus \widetilde G_{21}^{*} H^{2
   \perp}_{\cU} \to L^{2}_{{\mathcal Z}} \ominus \widetilde G_{12}
   H^{2}_{\cU}$ be the compression of multiplication by $\widetilde
   G_{11}$ to the spaces $L^{2}_{\cW} \ominus \widetilde G_{21}^{*} H^{2 \perp}_{\cU}$
   and  $L^{2}_{\cZ} \ominus \widetilde G_{12} H^{2}_{\cU}$, i.e.,
   $\Gamma = P_{L^{2}_{\cZ} \ominus\widetilde G_{12} H^{2}_{\cU}}
   \widetilde G_{11}|_{L^2_{\cW}\ominus\widetilde G_{21}^{*} H^{2 \perp}_{\cY}}$.
   %We note that in case $\widetilde G_{12}$ is injective and
   %$\widetilde G_{21}$ is surjective on the unit circle, without loss
   %$\widetilde G_{12} H^{2}_{\cU}$ is a closed subspace of
   %$H^{2}_{\cZ}$ and $\widetilde G_{21}^{\dagger} H^{2}_{\cY} \supset
   %H^{2}_{\cW}$ is a closed subspace of $L^{2}_{\cW}$, where
   %$\widetilde G_{21}^{\dagger}$
   %denotes the pointwise Moore-Penrose pseudo-inverse of $\widetilde G_{12}$.
   %Let $\Gamma \colon L^2_{\cW}\ominus\widetilde G_{21}^{\dagger} H^{2}_{\cY} \to
   %L^{2}_{\cZ} \ominus \widetilde G_{12} H^{2}_{\cU}$ be the compression of $G_{11}$
   %to the spaces $L^2_{\cW}\ominus\widetilde G_{21}^{\dagger} H^{2}_{\cY}$ and
   %$L^{2}_{\cZ} \ominus \widetilde G_{12} H^{2}_{\cU}$, i.e.,
   %$\Gamma=P_{L^{2}_{\cZ} \ominus\widetilde G_{12} H^{2}_{\cU}}G_{11}
   %P_{L^2_{\cW}\ominus\widetilde G_{21}^{\dagger} H^{2}_{\cY}}$.
   %Define an operator $\Gamma \colon L^2_{\cW}\ominus\widetilde G_{21}^{\dagger} H^{2}_{\cY} \to
%   L^{2}_{\cZ} \ominus \widetilde G_{12} H^{2}_{\cU}$ by
%   $$
%    \Gamma \colon \widetilde G_{21}^{\dagger} f \mapsto P_{L^{2}_{\cZ} \ominus
%    \widetilde G_{12} H^{2}_{\cU}}
%      (\widetilde G_{11} \widetilde G_{21}^{\dagger} f).
%   $$
   Then, as a consequence of the
   Commutant Lifting theorem (see \cite[Corollary 10.2 pages
   40--41]{FFGK}), one can see that
   the strict Model-Matching/Sarason interpolation problem
   posed in Theorem \ref{T:1.2} has a
   solution if and only if $\| \Gamma\|_{op} < 1$.
   Alternatively, in case $\widetilde G_{12}$ and $\widetilde G_{21}$
   are square and invertible on the unit circle, one can convert this
   Model-Matching/Commutant-Lifting problem to a bitangential
   Nevanlinna-Pick interpolation problem (see \cite[Theorem
   16.9.3]{BGR}), a direct generalization of the connection between
   a model-matching/Sarason interpolation problem with
   Nevanlinna-Pick interpolation as given in \cite{Sarason, FHZ} for
   the scalar case,  but we will not go into the details of this here.

%%%%%%%%%%%%%%%%%%%%%%%%%%%%%%%%%%%%%%%%%%%%%%%%%%%%%%%%%%%%%%%%%%%%%%%%%%%%%%%%%%%%%%%%%%%%%%%%
\subsection{The state-space approach}

We now restrict the classes of admissible plants and controllers to the
transfer matrices whose entries are in $\BC(z)_0$, the space of rational
functions without a pole at 0 (i.e., analytic in a neighborhood of 0).
In that case, a transfer matrix $F:\cU\to\cY$ with entries in $\BC(z)_0$
admits a {\em state-space realization}: There exists a quadruple $\{A,B,C,D\}$
consisting of matrices whose sizes are given by
\begin{equation}\label{sysmat}
\mat{cc}{A&B\\C&D}:\mat{c}{\cX\\\cU}\to\mat{c}{\cX\\\cY},
\end{equation}
where the {\em state-space} $\cX$ is finite dimensional, so that
\[
F(z)=D+zC(I-zA)^{-1}B
\]
for $z$ in a neighborhood of 0. Sometimes we consider quadruples
$\{A,B,C,D\}$ of operators, of compatible size as above, without
any explicit connection to a transfer matrix, in which case we just
speak of a realization.

Associated with the realization $\{A,B,C,D\}$ is the linear discrete-time
system of equations
\[
\Sigma:=\left\{\begin{array}{ccc}
x(n+1)&=&Ax(n)+Bu(n),\\
y(n)&=&Cx(n)+Du(n).
\end{array}\right.\qquad(n\in\BZ_+)
\]
The system $\Si$ and function $F$ are related through the fact that $F$ is the
{\em transfer-function} of $\Si$. The two-by-two matrix \eqref{sysmat} is called
the {\em system matrix} of the system $\Si$.

For the rest of this section we shall say that
   an operator $A$ on a finite-dimensional state space $\cX$ is {\em
   stable} if all its eigenvalues are in the open unit disk, or,
   equivalently, $\|A^{n}x\| \to 0$ as $n \to \infty$ for each $x \in
   \cX$. The following result deals with two key notions for the stabilizability problem
on the state-space level.

%
%   The internal stabilization/$H^{\infty}$-control problem calls for
%   a dynamic output controller. Before discussing the solution of
%   this problem in state-space coordinates, it makes sense to discuss
%   the easier problem of stabilization via static state-feedback or
%   output-injection (where the controller has full access to the
%   internal state).  For the rest of this section we shall say that
%   an operator $A$ on a finite-dimensional state space $\cS$ is {\em
%   stable} if all its eigenvalues are in the open unit disk, or,
%   equivalently, $\|A^{n}x\| \to 0$ as $n \to \infty$ for each $x \in
%   \cX$.
%
%   The basic results on this theme are as follows.

   \begin{theorem}  \label{T:1.3}
       (I) Suppose that $\{A,B\}$ is an {\em input pair}, i.e.,
       $A,B$ are operators with $A \colon \cX \to \cX$ and $B
       \colon \cU \to \cX$ for a finite-dimensional state space
       $\cX$ and a finite-dimensional input space $\cU$.  Then the
       following are equivalent:
       \begin{enumerate}
           \item $\{A, B\}$ is {\em operator-stabilizable}, i.e., there
           exists a state-feedback operator $F \colon \cX \to
           \cU$ so that the operator $A + BF$ is stable.

           \item $\{A, B\}$ is {\em Hautus-stabilizable}, i.e., the
           matrix
           pencil $\begin{bmatrix} I - z A & B \end{bmatrix}$ is
           surjective for each $z$ in the closed unit disk
           $\overline{\mathbb D}$.

           \item The Stein inequality
           $$
             A X A^{*} - X - B B^{*} < 0
           $$
           has a positive-definite solution $X$. Here $\Gamma < 0$
	   for a square matrix $\Gamma$ means that $-\Gamma$ is
	   positive definite.
       \end{enumerate}

      (II)  Dually, if $\{C,A\}$ is an output pair, i.e., $C,A$ are
       operators with $A \colon \cX \to \cX$ and $C \colon \cX
       \to \cY$ for a finite-dimensional state space $\cX$ and a
       finite-dimensional output space $\cY$, then the following
       are equivalent:
       \begin{enumerate}
           \item $\{C,A\}$ is {\em operator-detectable}, i.e.,
           there exists an output-injection operator $L \colon
           \cY \to \cX$ so that $A + LC$ is stable.

           \item $\{C,A\}$ is {\em Hautus-detectable}, i.e., the
           matrix pencil $\sbm{ I - zA \\ C }$ is injective for
           all $z$ in the closed disk $\overline{\mathbb D}$.

           \item The Stein inequality
           $$
             A^{*} Y A - Y - C^{*} C < 0
         $$
         has a positive definite solution $Y$.
         \end{enumerate}
       \end{theorem}

       When the input pair $\{A,B\}$ satisfies any one (and hence
       all) of the three equivalent conditions in part (I) of
       Theorem \ref{T:1.3}, we shall say simply that $\{A, B\}$ is
       {\em stabilizable}.  Similarly, if $(C,A)$ satisfies any
       one of the three equivalent conditions in part (II), we
       shall say simply that $\{C,A\}$ is detectable.  Given a
       realization $\{A,B,C,D\}$, we shall say that $\{A,B,C,D\}$ is
       stabilizable and detectable if $\{A,B\}$ is stabilizable and
       $\{C,A\}$ is detectable.

   In the state-space formulation of the internal
   stabilization/$H^{\infty}$-control problem, one assumes to be
   given a state-space realization for the plant $G$:
   \begin{equation}  \label{1.G}
   G(z) = \begin{bmatrix} D_{11} & D_{12} \\ D_{21} & D_{22}
\end{bmatrix} +  z  \begin{bmatrix} C_{1} \\ C_{2}
\end{bmatrix} (I - zA)^{-1} \begin{bmatrix} B_{1} & B_{2} \end{bmatrix}
    \end{equation}
    where the system matrix has the form
    \begin{equation}  \label{1.sysmat}
    \begin{bmatrix} A & B_{1} & B_{2} \\ C_{1} & D_{11} & D_{12} \\
    C_{2} & D_{21} & D_{22} \end{bmatrix} \begin{bmatrix} \cX \\
    \cW \\ \cU \end{bmatrix} \to \begin{bmatrix} \cX \\ \cZ \\
    \cY \end{bmatrix}.
     \end{equation}
    One then seeks a
     controller $K$ which is also given in terms of a state-space
     realization
     $$
     K(z) = D_{K} + z C_{K}(I - zA_{K})^{-1} B_{K}
     $$
     which provides internal stability (in the state-space sense to
     de defined below)
     and/or $H^{\infty}$-performance for the closed-loop system.
     Well-posedness of the closed-loop system is equivalent to invertibility of $I -
     D_{22}D_{K}$.  To keep various formulas affine in the design
     parameters $A_{K}, B_{K}, C_{K}, D_{K}$, it is natural to assume
     that $D_{22} = 0$; this is considered not unduly restrictive
     since under the assumption of well-posedness this can always be arranged
     via a change of variables (see
     \cite{IS}). Then the closed loop system $\Th(G,K)$ admits a state space
     realization $\{A_{cl},B_{cl},C_{cl},D_{cl}\}$ given by its system matrix
     \begin{equation}  \label{1clsysmat}
    \begin{bmatrix}  A_{cl} & B_{cl} \\ C_{cl} & D_{cl} \end{bmatrix} =
        \left[ \begin{array}{cc|c}
        A+B_{2} D_{K}C_{2} & B_{2} C_{K} & B_{1} + B_{2} D_{K} D_{21} \\
        B_{K} C_{2} & A_{K} & B_{K} D_{21} \\
        \hline C_{1} + D_{12} D_{K} C_{2} & D_{12} C_{K} & D_{11} +
        D_{12} D_{K} D_{21} \end{array} \right]
     \end{equation}
     and {\em internal stability (in the state-space sense)} is taken to
     mean that $A_{cl} = \sbm{ A + B_{2} D_{K} C_{2} & B_{2} C_{K} \\
     B_{K} C_{2} & A_{K} }$ should be stable, i.e., all eigenvalues
     are in the open unit disk.

    The following result characterizes when a given $G$ is internally
    stabilizable in the state-space sense.

     \begin{theorem}  \label{T:1.4} (See Proposition 5.2 in \cite{DP}.)
    Suppose that we are given a system matrix as in
    \eqref{1.sysmat} with $D_{22} = 0$ with associated transfer matrix $G$ as in
    \eqref{1.G}.  Then there exists a $K(z) = D_{K} + z
    C_{K} (I - z A_{K})^{-1} B_{K}$ which internally stabilizes
    $G$ (in the state-spaces sense) if and only if $\{A,B_2\}$ is stabilizable
    and $\{C_{2},A\}$ is detectable. In this case one such
    controller is given by the realization $\{A_K,B_K,C_K,D_K\}$ with system matrix
    $$ \mat{cc}{A_K&B_K\\C_K&D_K}= \begin{bmatrix} A+B_{2}F + LC_{2} &-L \\ F & 0
    \end{bmatrix}
    $$
    where $F$ and $L$ are state-feedback and output-injection
    operators chosen so that $A + B_{2}F$ and $A + L C_{2}$ are
    stable.
       \end{theorem}

In addition to the state-space version of the stabilizability problem we also
consider a (strict) state-space $H^\infty$ problem, namely to find a controller
$K$ given by a state-space realization $\{A_K,B_K,C_K,D_K\}$ of compatible size
so that the transfer-function $T_{zw}$ of the closed loop system, given by the
system matrix \eqref{1clsysmat}, is stable (in the state-space sense) and has
a supremum norm $\|T_{zw}\|_\infty$ of at most $1$ (less than 1).

  The definitive solution of the $H^{\infty}$-control problem in
  state-space coordinates for a time was the coupled-Riccati-equation
  solution due to Doyle-Glover-Khargonekar-Francis \cite{DGKF}.  This solution
  has now been superseded by the LMI solution of Gahinet-Apkarian \cite{GA}
  which can be stated as follows.  Note that the problem can be
  solved directly without first processing the data to the
  Model-Matching form.

  \begin{theorem}  \label{T:1.5}
      Let $\{A,B,C,D\} = \left\{ A, \sbm{ B_{1} & B_{2}}, \sbm{ C_{1} & C_{2}},
      \sbm{ D_{11} & D_{12} \\
  D_{21} & 0} \right\}$ be a given realization. Then there exists a solution for
      the strict state-space $H^\infty$-control problem associated with $\{A,B,C,D\}$ if and only if
      there exist positive-definite matrices $X,Y$ satisfying the LMIs
      \begin{align}
         & \begin{bmatrix} N_{c} & 0 \\ 0 & I \end{bmatrix} ^{*}
          \begin{bmatrix} AYA^{*} - Y & AYC_{1}^{*} & B_{1} \\
       C_{1}Y A^{*} & C_{1} Y C_{1}^{*} - I & D_{11} \\
         B_{1}^{*} & D_{11}^{*} & -I \end{bmatrix} \label{1.Y-LMI}
       \begin{bmatrix} N_{c} & 0 \\ 0 & I \end{bmatrix} <  0,\quad Y>0, \\
           & \begin{bmatrix} N_{o} & 0 \\ 0 & I \end{bmatrix} ^{*}
          \begin{bmatrix} A^{*} X A - X & A^{*} X B_{1} & C_{1} ^{*} \\
         B_{1}^{*} X A & B_{1}^{*} X B_{1} - I & D_{11}^{*} \\
         C_{1} & D_{11} & -I \end{bmatrix}
         \begin{bmatrix} N_{o} & 0 \\ 0 & I \end{bmatrix} < 0,\quad X>0,
         \label{1.X-LMI}
        \end{align}
      and the coupling condition
      \begin{equation}\label{1.coupling}
      \mat{cc}{X&I\\I&Y}\geq 0.
      \end{equation}
      Here $N_{c}$ and $N_{o}$ are matrices chosen so that
      \begin{align*}
          & N_{c} \text{ is injective and }
       \operatorname{Im } N_{c} = \operatorname{Ker } \begin{bmatrix}
       B_{2}^{*} & D_{12}^{*} \end{bmatrix} \text{ and } \\
       & N_{o} \text{ is injective and } \operatorname{Im } N_{o} =
       \operatorname{Ker } \begin{bmatrix} C_{2} & D_{21} \end{bmatrix}.
       \end{align*}
      \end{theorem}

    We shall discuss the proof of Theorem \ref{T:1.5} in Section
    \ref{S:ND-statespace} below in the context of a more general
    multidimensional-system $H^{\infty}$-control problem.

     The next result is the key to transferring from the
     frequency-domain version of the
     internal-stabilization/$H^{\infty}$-control problem to the
     state-space version.

     \begin{theorem}  \label{T:1.6}
     (See Lemma 5.5 in \cite{DP}.)  Suppose that the
     realization $\{A, B_{2},
     C_{2},0\}$ for the plant $G_{22}$ and the realization
     $\{A_{K}, B_{K}, C_{K}, D_{K}\}$ for the controller $K$ are both
     stabilizable and detectable.   Then $K$ internally
     stabilizes $G_{22}$ in the state-space sense if and only if
     $K$ stabilizes $G_{22}$ in the frequency-domain sense, i.e.,
     the closed-loop matrix $A_{cl} = \sbm{A + B_{2}D_{K} C_{2} &
     B_{2}C_{K} \\ B_{K} C_{2} & A_{K}}$ is stable if and only if
     the associated transfer matrix
     $$ \Theta(G_{22},K) = \begin{bmatrix} I & D_{K} \\ 0 & I
     \end{bmatrix} + z \begin{bmatrix} D_{K}C_{2} & C_{K} \\ C_{2} & 0
     \end{bmatrix} (I - z A_{cl})^{-1} \begin{bmatrix} B_{2} & B_{2}
     D_{K} \\ 0 & B_{K} \end{bmatrix}
     $$
     has all matrix entries in $\cR H^{\infty}$.
     \end{theorem}

   \subsection{Notes}  \label{S:1.notes}
   In the context of the discussion immediately after the statement
   of Theorem \ref{T:1.2},
   in case $\widetilde G_{12}$ and/or $\widetilde G_{21}$ drop rank at
   points on the unit circle, the Model-Matching problem in Theorem
   \ref{T:1.2} may convert
   to a boundary Nevanlinna-Pick interpolation problem for which
   there is an elaborate specialized theory (see  e.g.~Chapter 21 of
   \cite{BGR} and the more recent \cite{BD06}).
   However, if one sticks with the strictly suboptimal
   version of the problem, one can solve the problem with the
   boundary interpolation conditions if and only if one can solve the
   problem without the boundary interpolation conditions, i.e.,
   boundary interpolation conditions are irrelevant as far as
   existence criteria are concerned.  This is the route taken in the
   LMI solution of the $H^{\infty}$-problem and provides one
   explanation for the disappearance of any rank conditions in the
   formulation of the solution of the problem.  For a complete
   analysis of the relation between the coupled-Riccati-equation of
   \cite{DGKF} versus the LMI solution of \cite{GA}, we refer to
   \cite{Scherer}.

%%%%%%%%%%%%%%%%%%%%%%%%%%%%%%%%%%%%%%%%%%%%%%%%%%%%%%%%%%%%%%%%%%%%%%%%%%%%%%%%%%%%%%%%%%%%%%%%
%%%%%%%%%%%%%%%%%%%%%%%%%%%%%%%%%%%%%%%%%%%%%%%%%%%%%%%%%%%%%%%%%%%%%%%%%%%%%%%%%%%%%%%%%%%%%%%%
\section{The fractional representation approach to stabilizability and performance}\label{S:ring}
\setcounter{equation}{0}

In this section we work in the general framework of the fractional
representation approach to stabilization of linear systems as
introduced originally  by Desoer, Vidyasagar and coauthors
\cite{DLMS, VSF} in the 1980s and refined only recently in the work
of Quadrat \cite{Q-YKI, Q-Lattice, Q-YKII}.  For an overview of the
more recent developments we recommend the survey article \cite{Q-LN}
and for a completely elementary account of the generalized
Youla-Ku\v{c}era parametrization with all the algebro-geometric
interpretations stripped out we recommend \cite{Q-elementary}.

The set of {\em stable} single-input single-output (SISO) transfer functions is
assumed to be given by
a general ring ${\mathbb A}$ in place of the ring $\cR H^{\infty}$
used for the classical case as discussed in Section \ref{S:1D}; the only assumption
which we shall impose on
${\mathbb A}$ is that it be a commutative integral domain.  It
therefore has a quotient field ${\mathbb K}: = Q({\mathbb A}) = \{ n/d \colon
d, n \in {\mathbb A}, d\not=0\}$ which shall be
considered as the set of all possible SISO transfer functions (or plants).
Examples of ${\mathbb A}$ which come up include the ring ${\mathbb
R}_{s}(z)$ of real rational functions of the complex variable $z$
with no poles in the closed right half plane, the Banach algebra
$R H^{\infty}({\mathbb C}_{+})$ of all bounded analytic functions on the
right half plane ${\mathbb C}_{+}$ which are real on the positive real axis,
and their discrete-time analogues: (1) real rational functions with no poles in the
closed unit disk (or closed exterior of the unit disk depending on
how one sets conventions), and (2) the Banach algebra $R
H^{\infty}({\mathbb D})$ of all bounded holomorphic functions on the
unit disk ${\mathbb D}$ with real values on the real interval $(-1,
1)$.  There are also Banach subalgebras of $R H^{\infty}({\mathbb
C}_{+})$ or $R H^{\infty}({\mathbb D})$ (e.g., the Wiener algebra and
its relatives such as the Callier-Desoer class---see \cite{CZ}) which
are of interest.  In addition to these examples there are
multivariable analogues, some of which we shall discuss in the next
section.

We now introduce some notation.
We assume that the {\em control-signal space } $\cU$, the {\em
disturbance-signal space} ${\cW}$, the {\em error-signal space}
${\cZ}$ and the {\em measurement signal space} $\cY$ consist of
column vectors of given sizes $n_{\cU}$, $n_{\cW}$, $n_{\cZ}$ and
$n_{\cY}$, respectively, with entries from the quotient field
${\mathbb K}$ of ${\mathbb A}$:
$$
\cU = \BK^{n_{\cU}}, \quad \cW = \BK^{n_{\cW}}, \quad
\cZ = \BK^{n_{\cZ}}, \quad \cY = \BK^{n_{\cY}}.
$$
We are given a plant $G= \sbm{ G_{11} & G_{12} \\ G_{21} & G_{22}}
\colon \cW \oplus \cU \to \cZ \oplus \cY$ and seek to design a
controller $K \colon \cY \to \cU$ that stabilizes the system
$\Sigma(G,K)$ of Figure \ref{fig:taps0} as given in Section \ref{S:1D}.
The various matrix entries $G_{ij}$ of $G$ are now  matrices
with entries from $\BK$ (rather than $\cR H^{\infty}$ as in the
classical case) of compatible sizes (e.g., $G_{11}$ has size
$n_{\cW} \times n_{\cU}$) and $K$ is a matrix over $\BK$ of size
$n_{\cU} \times n_{\cY}$. Again $v_1$ and $v_2$ are tap signals
used to detect stability properties of the internal signals $u$ and $y$.

Just as was explained in Section \ref{S:1D} for the classical case,
the system $ \Sigma(G,K)$ is {\em well-posed} if there is a well-defined
map from $ \sbm{ w \\ v_{1} \\ v_{2} }$ to $\sbm{ z \\ u \\ y }$ and
this happens exactly when $\det( I - G_{22}K) \ne 0$ (where the
determinant now is an element of ${\mathbb A}$); when this is the
case, the map from $\sbm{ w \\ v_{1} \\ v_{2}}$ to
  $\sbm{z \\ u \\ y }$ is given by
$$
 \begin{bmatrix} z \\ u \\ y \end{bmatrix} = \Theta(G,K)
     \begin{bmatrix} w \\ v_{1} \\ v_{2} \end{bmatrix}
$$
where $\Theta(G,K)$ is given by \eqref{ThetaGK}.
We say that the system $\Sigma(G, K)$ is {\em internally stable} if
$\Sigma(G,K)$ is well-posed and, in addition, if the map $\Theta(G, K)$
maps $\BA^{n_{\cW}} \oplus \BA^{n_{\cU}} \oplus  \BA^{n_{\cY}}$ into
$\BA^{n_{\cZ}} \oplus \BA^{n_{\cU}} \oplus \BA^{n_{\cY}}$, i.e., stable
inputs $w,v_{1}, v_{2}$ are mapped to stable outputs $z,u,y$. Note that
this is the same as the entries of $\Sigma(G,K)$ being in $\BA$.

To formulate the standard problem of $H^{\infty}$-control, we assume
that $\BA$ is equipped with a positive-definite inner product making
$\BA$ at least a pre-Hilbert space with norm $\| \cdot \|_{\BA}$; in
the classical case, one takes this norm to be the $L^{2}$-norm over
the unit circle.  Then
we say that the system $\Sigma(G, K)$ {\em has performance} if
$\Sigma(G,K)$ is internally stable and in addition the transfer
function $T_{zw}$ from $w$ to $z$ has induced operator norm bounded
by some tolerance which we normalize to be equal to 1:
$$  \| T_{zw} \|_{op}: =\sup  \{ \| z \|_{\BA^{n_{\cZ}}} \colon
\|w\|_{\BA^{n_{\cW}}} \le 1, v_{1} = 0,\, v_{2} = 0\} \le 1.
$$
We say that the system $\Sigma(G,K)$ has {\em strict performance} if
in fact $\|T_{zw}\|_{op} < 1$.
The {\em stabilization problem} then is to describe all (if any
exist) internally stabilizing controllers $K$ for the given plant
$G$, i.e., all $K \in \BK^{n_{\cU} \times n_{\cY}}$ so that the
associated closed-loop system $\Sigma(G, K)$ is internally stable.
The {\em standard $H^{\infty}$-problem} is to find all internally
stabilizing controllers which in addition achieve performance
$\|T_{zw}\|_{op} \le 1$. The strictly suboptimal $H^{\infty}$-problem
is to describe all internally stabilizing controllers which achieve
strict performance $\| T_{zw}\|_{op} < 1$.

The $H^{\infty}$-control problem for the special case where $G_{22}=0$
is the Model-Matching problem for this setup. With the same arguments
as in Subsection \ref{S:MMproblem} it follows that stabilizability forces
$G_{11}$, $G_{12}$ and $G_{21}$ all to be stable (i.e., to have all
matrix entries in ${\mathbb A}$) and then $K$ stabilizes exactly when
also $K$ is stable.

%%%%%%%%%%%%%%%%%%%%%%%%%%%%%%%%%%%%%%%%%%%%%%%%%%%%%%%%%%%%%%%%%%%%%%%%%%%%%%%%%%%%%%%%%%%%%%%%
\subsection{Parametrization of stabilizing controllers in terms of a given
stabilizing controller} \label{S:stab-w/oCP}
We return to the general case i.e., $G= \sbm{ G_{11} & G_{12} \\ G_{21} & G_{22}}
\colon \cW \oplus \cU \to \cZ \oplus \cY$.
Now suppose we have a stabilizing controller $K\in \BK^{n_{\cU}\times
n_{\cY}}$. Set
\begin{equation}\label{UV}
U=(I-G_{22}K)^{-1} \quad\text{and}\quad V=K(I-G_{22}K)^{-1}.
\end{equation}
Then $U\in \BA^{n_{\cY}\times n_{\cY}}$, $V\in \BA^{n_{\cU}\times n_{\cY}}$,
$\det U\not=0\in \BA$, $K=VU^{-1}$ and $U-G_{22}V=I$. Furthermore,
$\Th(G,K)$ can then be
written as
\begin{equation}\label{Theta2}
\Th(G,K) = \Th(G; U,V):=
\mat{ccc}{G_{11}+G_{12}VG_{21}&G_{12}+G_{12}VG_{22}&G_{12}V\\
VG_{21}&I+VG_{22}&V\\
UG_{21}&UG_{22}&U}.
\end{equation}
It is not hard to see that if $U\in \BA^{n_{\cY}\times n_{\cY}}$
and $V\in \BA^{n_{\cU}\times n_{\cY}}$ are such that $\det U\not=0$,
$U-G_{22}V=I$ and \eqref{Theta2} is stable, i.e., in
$\BA^{(n_{\cZ} + n_{\cU} + n_{\cY})\times (n_{\cW} + n_{\cU }+ n_{\cY})}$,
then $K=VU^{-1}$ is a stabilizing controller.  A dual result holds if we set
\begin{equation}\label{tilUtilV}
\widetilde U =
(I - K G_{22})^{-1} \quad\text{and}\quad \widetilde V = (I - K G_{22})^{-1}K.
\end{equation}
In that case $\widetilde U\in \BA^{n_{\cU}\times n_{\cU}}$,
$\widetilde V\in \BA^{n_{\cU}\times n_{\cY}}$, $\det \widetilde U\not=0\in \BA$,
$K=\widetilde U^{-1}\widetilde V$, $\widetilde U-\widetilde VG_{22}=I$ and we
can write $\Th(G,K)$ as
\begin{equation}   \label{Theta2'}
       \Th(G,K)=\Theta(G; \widetilde U,\widetilde V) =
       \begin{bmatrix} G_{11} + G_{12} \widetilde V G_{21} &
           G_{12} \widetilde U & G_{12} \widetilde V \\
           \widetilde V G_{21} & \widetilde U & \widetilde V \\
           (I + G_{22} \widetilde V) G_{21} & G_{22} \widetilde U
           & I + G_{22} \widetilde V \end{bmatrix},
       \end{equation}
while conversely, for any $\widetilde U\in \BA^{n_{\cU}\times n_{\cU}}$ and
$\widetilde V\in \BA^{n_{\cU}\times n_{\cY}}$ with $\det \widetilde U\not=0$ and
$\widetilde U-\widetilde VG_{22}=I$ and such that
\eqref{Theta2'} is stable, we have that $K=\widetilde U^{-1}\widetilde V$
is a stabilizing controller.

This leads to the following first-step more linear reformulation of
the definition of internal stabilization.

\begin{theorem}\label{T:stab1}
A plant $G$ defined by a transfer matrix $G \in \BK^{ (n_{\cZ} + n_{\cY}) \times (n_{\cW} +
n_{\cU}) }$ is internally stabilizable if
and only if one of the following equivalent assertions holds:
\begin{enumerate}
    \item There exists $L = \sbm{ V \\ U } \in \BA^{(n_{\cU} +
    n_{\cY}) + n_{\cY}}$ with $\det U \ne 0$ such that:
    \begin{enumerate}
    \item The block matrix \eqref{Theta2} is stable (i.e., has all
    matrix entries  in $\BA$), and
    \item $\begin{bmatrix} -G_{22} & I \end{bmatrix} L = I$.
    \end{enumerate}
  Then the controller $K = V U^{-1}$ internally stabilizes the plant
  $G$ and we have:
  $$
    U = (I - G_{22} K)^{-1}, \quad V = K (I - G_{22}K)^{-1}.
   $$

   \item There exists $\widetilde L =\sbm{ \widetilde U & -
 \widetilde V} \in \BA^{n_{\cU} \times (n_{\cU} +
   n_{\cY})}$ with $\det \widetilde U \ne 0$ such that:
   \begin{enumerate}
       \item The block matrix \eqref{Theta2'} is stable (i.e., has
       all matrix entries in $\BA$), and
 \item  $\widetilde L  \sbm{ I \\  G_{22}} : = \sbm{ \widetilde U & -
 \widetilde V} \sbm{ I \\ G_{22}} = I$.
 \end{enumerate}
 If this is the case, then the controller $K = \widetilde U^{-1}
 \widetilde V$ internally stabilizes the plant $G$ and we have:
 $$
 \widetilde U = (I - K G_{22})^{-1}, \quad \widetilde V = (I - K
 G_{22})^{-1} K.
 $$
  \end{enumerate}
\end{theorem}

With this result in hand, we are able to get a parametrization for
the set of all stabilizing controllers in terms of an assumed
particular stabilizing controller.

\begin{theorem}\label{T:stab2}$ $\\[-.5cm]
    \begin{enumerate}
    \item
Let $K_*\in \BK^{n_{\cU}\times n_{\cY}}$ be a stabilizing controller
for $G\in \BK^{(n_{\cZ}+n_{\cY})\times (n_{\cW}+n_{\cU})}$.
Define $U_*=(I-G_{22}K_*)^{-1}$ and $V_*=K(I-G_{22}K_*)^{-1}$. Then the set of all
stabilizing controllers is given by
\begin{equation}\label{QtoK1}
K=(V_*+Q)(U_*+G_{22}Q)^{-1},
\end{equation}
where $Q \in \BK^{n_{\cU} \times n_{\cY}}$ is an element of the set
\begin{equation}  \label{Omega}
\Om:=\left\{Q\in \BK^{n_{\cU} \times n_{\cY}} \colon
\begin{bmatrix} G_{12} \\ I \\ G_{22} \end{bmatrix} Q \begin{bmatrix}
G_{21} & G_{22} & I \end{bmatrix} \in
\BA^{(n_{\cZ}+n_{\cU}+n_{\cY}) \times  (n_{\cW}+n_{\cU}+n_{\cY})}\right\}
\end{equation}
such that in addition $\det(U_*+G_{22}Q)\not=0$.

\item Let $K_{*} \in \bK^{n_{\cU} \times n_{\cY}}$ be a
stabilizing controller for $G \in \BK^{(n_{\cZ}+n_{\cY}) \times (n_{\cW} + n_{\cU})}$.
Define $\widetilde U_{*} = (I - K_{*}
G_{22})^{-1}$ and $\widetilde V_{*} = (I - K_{*}G_{22})^{-1} K_{*}$.
Then the set of all controllers is given by
\begin{equation}\label{QtoK2}
K=(\widetilde U_{*}  + Q G_{22})^{-1} (\widetilde V_{*} + Q),
\end{equation}
where $Q \in \BK^{n_{\cU} \times n_{\cY}}$ is an element of the set
$\Omega$ \eqref{Omega} such that in addition
$\det (\widetilde U_{*} + Q G_{22}) \ne 0$.
\end{enumerate}
\noindent Moreover, if $Q\in\Omega$, that $\det(U_*+G_{22}Q)\not=0$ if and only if
$\det (\widetilde U_{*} + Q G_{22}) \ne 0$, and the formulas \eqref{QtoK1}
and \eqref{QtoK2} give rise to the same controller $K$.
\end{theorem}

\begin{proof}  By Theorem \ref{T:stab1}, if $K$ is a stabilizing
    controller for $G$,
then $K$ has the form $K = V U^{-1}$ with $L = \sbm{ U \\ V}$ as in part (1) of Theorem
    \ref{T:stab1} and then $\Theta(G,K)$ is as in \eqref{Theta2}.
    Similarly $\Theta(G,K_{*})$ is given as $\Theta(G; U_{*},V_{*})$
    in \eqref{Theta2} with
    $U_{*},V_{*}$ in place of $U,V$.  As by assumption
    $\Theta(G; U_{*},V_{*})$ is stable, it follows that $\Theta(G; U,V)$
    is stable if and only if $\Theta(G;  U,V) - \Theta(G;
    U_{*},V_{*})$ is stable. Let $Q=V-V_*$; as $U=I+G_{22}V$ and $U_*=I+G_{22}V_*$,
it follows that $U-U_*=G_{22}Q$. {}From \eqref{Theta2} we then see that the stable
quantity $\Th(G; U,V)-\Th(G; U_*,V_*)$ is given by
 $$
    \Theta(G; U,V) - \Theta(G; U_{*},V_{*}) =
     \begin{bmatrix} G_{12} \\ I \\ G_{22}
    \end{bmatrix} Q \begin{bmatrix} G_{21} & G_{22} & I \end{bmatrix}.
   $$
Thus
\[
K=VU^{-1}=(V_*+(V-V_*))(U_*+(U-U_*))^{-1}=(V_*+Q)(U_* + G_{22}Q)^{-1},
\]
where $Q$ is an element of $\Om$ such that $\det(U_*+G_{22}Q)\not=0$.

Conversely, suppose $K$ has the form $K=(V_*+Q)(U_* + G_{22}Q)^{-1}$ where
$Q\in\Om$ and $\det(U_*+G_{22}Q)\not=0$. Define $V=V_*+Q$, $U=U_*+G_{22}Q$.
Then one easily checks that
$$
    \Theta(G; U,V) = \Theta(G; U_{*},V_{*})+
     \begin{bmatrix} G_{12} \\ I \\ G_{22}
    \end{bmatrix} Q \begin{bmatrix} G_{21} & G_{22} & I \end{bmatrix}
   $$
is stable and
\[
\mat{cc}{-G_{22}&I}\mat{c}{V\\U}=\mat{cc}{-G_{22}&I}\mat{c}{V_*\\U_*}
+\mat{cc}{-G_{22}&I}\mat{c}{Q\\G_{22}Q}=I+0=I.
\]
So $K=VU^{-1}$ stabilizes $G$ by part (1) of Theorem \ref{T:stab1}.
   This completes the proof of the first statement of the theorem.
   The second part follows in  a similar way by using the second statement
   in Theorem \ref{T:stab1} and $Q=\wtilV-\wtilV_*$. Finally, since
   $V=\wtilV$ and $V_*=\wtilV_*$, we find that indeed $\det(U_*+G_{22}Q)\not=0$
   if and only if $\det (\widetilde U_{*} + Q G_{22}) \ne 0$, and the formulas
   \eqref{QtoK1} and \eqref{QtoK2} give rise to the same controller $K$.
   \end{proof}

   The drawback of the parametrization of the stabilizing controllers
   in Theorem \ref{T:stab2} is that the set $\Omega$ is not really a
   free-parameter set.  By definition, $Q \in \Omega$ if $Q$ itself
   is stable (from the (1,3) entry in the defining matrix for the
   $\Omega$ in \eqref{Omega}), but, in addition, the eight additional
   transfer matrices
   \begin{align*}  G_{12}Q G_{21}, \quad G_{12} Q G_{22},  \quad G_{12} Q, \quad
   Q G_{21},  \\  Q G_{22}, \quad G_{22}Q G_{21}, \quad G_{22} Q
   G_{22}, \quad G_{22} Q
   \end{align*}
   should all be stable as well. The next lemma shows how the parameter set
   $\Omega$ can in turn be parametrized by a free stable parameter
   $\Lambda$ of size $(n_{\cU} + n_{\cY}) \times (n_{\cU} + n_{\cY})$.

   \begin{lemma}  \label{L:free-param} Assume that $G$ is
       stabilizable and that $K_{*}$ is a particular stabilizing
       controller for $G$. Let $Q\in\BK^{n_\cU\times n_\cY}$. Then the following are equivalent:
\begin{itemize}

\item[(i)] $Q$ is an element of the set $\Omega$ in \eqref{Omega},

\item[(ii)] $\mat{c}{I\\G_{22}}Q\mat{cc}{G_{22}&I}$ is stable,

\item[(iii)] $Q$ has the form $Q = \widetilde
       L \Lambda L$ for a stable free-parameter $\Lambda \in
       \BA^{(n_{\cU} +n_{\cY}) \times (n_{\cU} + n_{\cY})}$, where
       $\widetilde L \in {\mathbb A}^{n_{\cU} \times (n_{\cU} +
       n_{\cY})}$ and $L \in \BA^{(n_{\cU} + n_{\cY}) \times n_{\cY}}$
       are given by
       \begin{equation}  \label{tildeLL}
\widetilde L = \begin{bmatrix} (I - K_{*} G_{22})^{-1} & -
       (I - K_{*} G_{22})^{-1} K_{*} \end{bmatrix}, \quad
       L = \begin{bmatrix} -K_{*} (I - G_{22}K_{*})^{-1} \\ (I - G_{22}
       K_{*})^{-1} \end{bmatrix}.
       \end{equation}

\end{itemize}
       \end{lemma}

\begin{proof} The implication (i) $\Longrightarrow$ (ii) is obvious. Suppose that
       $\La=\sbm{I\\G_{22}}Q\sbm{G_{22}&I}$ is stable. Note that
\begin{align*}
    \widetilde L \Lambda L & =
    \begin{bmatrix} (I - K_{*} G_{22})^{-1} & -
       (I - K_{*} G_{22})^{-1} K_{*} \end{bmatrix}
       \begin{bmatrix} I \\ G_{22} \end{bmatrix} Q\times\\
&           \qquad\qquad\qquad\qquad\qquad\qquad\times\begin{bmatrix} G_{22} & I \end{bmatrix}
           \begin{bmatrix} -K_{*} (I - G_{22}K_{*})^{-1} \\ (I - G_{22}
              K_{*})^{-1} \end{bmatrix} \\
              & = Q.
\end{align*}
Hence (ii) implies (iii). Finally assume $Q=\widetilde L \Lambda L$ for a
stable $\La$. To show that $Q\in\Om$, as $\Lambda$ is stable, it suffices to show that
$$
L_{1}: = \begin{bmatrix} G_{12} \\ I \\ G_{22} \end{bmatrix}
\widetilde L \text{ is stable, and }
L_{2}: = L \begin{bmatrix} G_{21} & G_{22} & I \end{bmatrix} \text{
is stable.}
$$
Spelling out $L_{1}$, using the definition of $\widetilde L$ from
\eqref{tildeLL}, gives
$$
 L_{1} = \begin{bmatrix} G_{12} \\ I \\ G_{22} \end{bmatrix}
 \begin{bmatrix} (I - K_{*} G_{22})^{-1}  & - (I - K_{*} G_{22})^{-1}
     K_{*} \end{bmatrix}.
 $$
 We note that each of the six matrix entries of $L_{1}$ are stable,
 since they all occur among the matrix entries of $\Theta(G,K_{*})$
 (see \eqref{ThetaGK}) and $K_{*}$ stabilizes $G$ by assumption.
 Similarly, each of the six matrix entries of $L_{2}$ given by
 $$
 L_{2} = \begin{bmatrix} -K_{*} (I - G_{22}K_{*})^{-1} \\ (I - G_{22}
 K_{*})^{-1} \end{bmatrix} \begin{bmatrix} G_{21} & G_{22} & I
\end{bmatrix}
$$
is stable since $K_{*}$ stabilizes $G$.  It therefore follows that
$Q\in \Omega$ as wanted.
\end{proof}

We say that $K$ stabilizes $G_{22}$ if the map $\sbm{v_1\\v_2}\mapsto\sbm{u\\y}$
in Figure \ref{fig:taps0} is stable, i.e., the usual stability holds with
$w=0$ and $z$ ignored. This amounts to the stability of the lower right $2\times 2$
block in $\Th(G,K)$:
\[
\mat{cc}{(I - K G_{22})^{-1} & (I - K G_{22})^{-1} K\\
G_{22}(I - K G_{22})^{-1} & I + G_{22}(I - K G_{22})^{-1} K}.
\]

The equivalence of (i) and (ii) in Lemma \ref{L:free-param} implies
the following result.

\begin{corollary}\label{C:G22G}  Assume that $G$ is stabilizable.  Then $K$
    stabilizes $G$ if and only if $K$ stabilizes $G_{22}$.
 \end{corollary}

\begin{proof}
Assume $K_*\in\BK^{n_\cU\times n_\cY}$ stabilizes $G$. Then in particular
the lower left $2\times 2$ block in $\Theta(G,K_*)$ is stable. Thus $K_*$
stabilizes $G_{22}$. Moreover, $K$ stabilizes $G_{22}$ if and only if $K$
stabilizes $G$ when we impose $G_{11}=0$, $G_{12}=0$ and $G_{21}=0$, that
is, $K$ is of the form \eqref{QtoK1} with $U_*$ and $V_*$ as in Theorem
\ref{T:stab2} and $Q\in \BK^{n_{\cU} \times n_{\cY}}$  is such that
$\sbm{I\\G_{22}}Q\sbm{G_{22}&I}$ is stable. But then it follows from the implication
(ii) $\Longrightarrow$ (i) in Lemma \ref{L:free-param} that $Q$ is in $\Omega$,
and thus, by Theorem \ref{T:stab2}, $K$ stabilizes $G$ (without $G_{11}=0$,
$G_{12}=0$, $G_{21}=0$).
\end{proof}

 Combining Lemma \ref{L:free-param} with Theorem \ref{T:stab2} leads
 to the following generalization of Theorem \ref{T:1.1} giving a
 parametrization of stabilizing controllers without
 the assumption of any coprime factorization.

 \begin{theorem} \label{T:stab3}  Assume that $G \in \BK^{(n_{\cZ} +
    n_{\cY}) \times (n_{\cW} + n_{\cU})}$ is stabilizable and that
    $K_{*}$ is one stabilizing controller for $G$.  Define $U_{*} =
    (I - G_{22} K_{*})^{-1}$, $V_{*} = K_{*}(I - G_{22} K_{*})^{-1}$,
    $\widetilde U_{*} = (I - K_{*} G_{22})^{-1}$ and $\widetilde V_{*} = (I -
    K_{*} G_{22})^{-1} K_{*}$.
    Then the set of all stabilizing controllers for $G$ are given by
    \begin{align*}
    & K = (V_{*} + Q) (U_{*} + G_{22} Q)^{-1}
    =(\widetilde U_{*} + Q G_{22})^{-1} (\widetilde V_{*} + Q),
  \end{align*}
  where $Q = \widetilde L \Lambda L$  where $\widetilde L$ and $L$ are
  given by \eqref{tildeLL} and $\Lambda$ is a free stable parameter
  of size $(n_{\cU} + n_{\cY}) \times (n_{\cU} + n_{\cY})$ so that
  $\det(U_{*} + G_{22} Q)\not=0$ or equivalently
  $\det(\widetilde U_{*} + Q G_{22})\not=0$.
  \end{theorem}

\subsection{The Youla-Ku\v cera parametrization}

There are two drawbacks to the para\-metrization of the stabilizing controllers
obtained in Theorem \ref{T:stab3}, namely, to find all stabilizing controllers one first
has to find a particular stabilizing controller, and secondly, the map $\Lambda\mapsto Q$
given in Part (iii) of Lemma \ref{L:free-param} is in general not one-to-one.
We now show that, under the
additional hypothesis that $G_{22}$ admits a double coprime factorization, both issues
can be remedied, and we are thereby led to the well known Youla-Ku\v cera parametrization
for the stabilizing controllers.

Recall that $G_{22}$ has a double coprime factorization in case there exist stable
transfer matrices $D$, $N$, $X$, $Y$, $\widetilde D$, $\widetilde N$, $\widetilde X$
and $\widetilde Y$ so that the determinants of $D$, $\widetilde D$, $X$ and
$\widetilde X$ are all nonzero (in $\BA$) and
\begin{equation}\label{coprime}
G_{22} = D^{-1} N = \widetilde N \widetilde D^{-1}, \quad
\begin{bmatrix} D & -N \\ -\widetilde Y & \widetilde X
\end{bmatrix} \begin{bmatrix} X & \widetilde N \\ Y &
\widetilde D \end{bmatrix} =\mat{cc}{I_{n_{\cY}}&0\\0&I_{n_{\cU}}}.
\end{equation}

According to Corollary \ref{C:G22G} it suffices to focus on describing the
stabilizing controllers of $G_{22}$. Note that $K$ stabilizes $G_{22}$
means that
\[
\mat{cc}{(I - K G_{22})^{-1} & (I - K G_{22})^{-1} K\\
G_{22}(I - K G_{22})^{-1} & I + G_{22}(I - K G_{22})^{-1} K}
\]
is stable, or, by Theorem \ref{T:stab2}, that $K$ is given by \eqref{QtoK1} or \eqref{QtoK2}
for some $Q\in\BK^{n_\cU\times n_\cY}$ so that $\sbm{I\\G_{22}}Q\sbm{G_{22}&I}$ is stable.

In case $G_{22}$ has a double coprime factorization Quadrat shows in
\cite[Proposition 4]{Q-elementary} that the equivalence of (ii) and (iii) in
Lemma \ref{L:free-param} has the following refinement. We provide a proof
for completeness.

\begin{lemma}  \label{L:free-paramYK}
Suppose that $G_{22}$ has a double coprime factorization \eqref{coprime}.
Let $Q\in\BK^{n_\cU\times n_\cY}$. Then $\sbm{I\\G_{22}}Q\sbm{G_{22}&I}$
is stable if and only if $Q=\widetilde D \Lambda D$ for some
$\La\in\BA^{n_\cU\times n_\cY}$.
\end{lemma}

\begin{proof}
Let $Q=\wtilD\La D$ for some $\La\in\BA^{n_\cU\times n_\cY}$. Then
\[
\mat{c}{I\\G_{22}}Q\mat{cc}{G_{22}&I}=\mat{cc}{QG_{22}&Q\\G_{22}QG_{22}&G_{22}Q}=
\mat{cc}{\wtilD\La N&\wtilD\La D\\
\wtilN\La N&\wtilN\La D}.
\]
Hence $\sbm{I\\G_{22}}Q\sbm{G_{22}&I}$ is stable.

Conversely, assume that $\sbm{I\\G_{22}}Q\sbm{G_{22}&I}$ is stable.
Set $\La=\wtilD^{-1}QD^{-1}$. Then with $X$, $Y$, $\wtilX$ and $\wtilY$ the transfer
matrices from the coprime factorization \eqref{coprime} we have
\begin{eqnarray*}
\La&=&\mat{cc}{\wtilX&-\wtilY}\mat{c}{\wtilD\\\wtilN}\La\mat{cc}{N&D}\mat{c}{-Y\\X}\\
&=&\mat{cc}{\wtilX&-\wtilY}\mat{cc}{\wtilD\La N&\wtilD\La
D\\\wtilN\La N&\wtilN\La D}\mat{c}{-Y\\X}\\
&=&\mat{cc}{\wtilX&-\wtilY}\mat{cc}{QG_{22}&Q\\G_{22}QG_{22}&G_{22}Q}\mat{c}{-Y\\X}.
\end{eqnarray*}
Thus $\La$ is stable.
\end{proof}

%%%%%%%%%%%%%%%%%%%%%%%%%%%%%%%%%%%%%%%%%%%%%%%%%%%%%%%%%%%%%%%%%%%%%%%%
\begin{lemma}\label{L:coprimestab}
Assume that $G_{22}$ admits a double coprime factorization \eqref{coprime}. Then
$K_{0}$ is a stabilizing controller for $G_{22}$ if and only if there exist
$X_0\in \BA^{n_\cY\times n_\cY}$, $Y_0\in \BA^{n_\cU\times n_\cY}$, $\wtilX_0\in \BA^{n_\cU\times n_\cU}$
and $\wtilY_0\in \BA^{n_\cU\times n_\cY}$ with $\det(X_{0})\not=0$,
$\det(\wtilX_{0})\not=0$ so
that $K_{0}=Y_0X_0^{-1}=\wtilX_0^{-1}\wtilY_0$ and
\begin{equation*}\label{coprime2}
\mat{cc}{D&-N\\-\wtilY_0&\wtilX_0}\mat{cc}{X_0&\wtilN\\Y_0&\wtilD}
=\mat{cc}{I_{n_\cY}&0\\0&I_{n_\cU}}.
\end{equation*}
In particular, $K=YX^{-1}=\wtilX^{-1}\wtilY$ is a stabilizing
controller for $G_{22}$, where $X,Y, \widetilde X, \widetilde Y$ come
from the double coprime factorization \eqref{coprime} for $G_{22}$.
\end{lemma}

\begin{proof}
Note that if $K$ is a stabilizing controller for $G_{22}$, then, in particular,
\begin{equation}\label{G22stab}
\mat{cc}{I&-K\\ -G_{22}&I}^{-1}
=\mat{cc}{(I-KG_{22})^{-1}&K(I-G_{22}K)^{-1}\\
(I-G_{22}K)^{-1}G_{22}&(I-G_{22}K)^{-1}}
\end{equation}
is stable.
The above identity makes sense, irrespectively of $K$ being a stabilizing
controller, as long as the left hand side is invertible.
Let $X$, $Y$, $\wtilX$ and $\wtilY$ be the transfer matrices from the
double coprime factorization. Set $K=\wtilX^{-1}\wtilY=YX^{-1}$. Then we have
\begin{eqnarray*}
\mat{cc}{\wtilX&-\wtilY\\-N&D}^{-1}\mat{cc}{\wtilX&0\\0&D}
&=&\left(\mat{cc}{\wtilX^{-1}&0\\0&D^{-1}}\mat{cc}{\wtilX&-\wtilY\\-N&D}\right)^{-1}\\
&=&\mat{cc}{I&-\wtilX^{-1}\wtilY\\-D^{-1}N&I}=\mat{cc}{I&-K\\-G_{22}&I}^{-1}.
\end{eqnarray*}
Since $\wtilX$, $D$ and
$\left[ \begin{smallmatrix} \wtilX&-\wtilY\\-N&D \end{smallmatrix}
\right]^{-1}=\left[ \begin{smallmatrix} \wtilD&Y\\\wtilN&X
\end{smallmatrix} \right]$
are stable, it follows that the right-hand side of \eqref{G22stab} is stable as well.
We conclude that $K = \widetilde X^{-1} \widetilde Y = Y X^{-1}$
stabilizes $G_{22}$.

Now let $K_0$ be any stabilizing controller for $G_{22}$. It follows from the
first part of the proof that $K=YX^{-1}=\wtilX^{-1}\wtilY$ is  stabilizing
 for $G_{22}$. Define $V$ and $U$ by \eqref{UV} and $\wtilV$ and $\wtilU$ by
\eqref{tilUtilV}. Then, using Theorem \ref{T:stab2} and Lemma \ref{L:free-paramYK},
there exists a $\La\in\BA^{n_\cU\times n_\cY}$ so that
\[
K_0=(V+Q)(U+G_{22}Q)^{-1}=(\wtilU+QG_{22})^{-1}(\wtilV+Q),
\]
where $Q=\widetilde D \Lambda D$. We compute that
\begin{equation}\label{comp4}
(I-G_{22}K)^{-1}=(I-D^{-1}NYX^{-1})^{-1}=X(DX-NY)^{-1}D=XD
\end{equation}
and
\begin{equation}\label{comp5}
(I-KG_{22})^{-1}=(I-\wtilX^{-1}\wtilY\wtilN\wtilD^{-1})^{-1}
=\wtilD(\wtilX\wtilD-\wtilY\wtilN)^{-1}\wtilX=\wtilD\wtilX.
\end{equation}
Thus
\[
V=YD,\quad U=XD,\quad \wtilV=\wtilD\wtilY,\quad \wtilU=\wtilD\wtilX.
\]
Therefore
\begin{eqnarray}
\label{comp1}
K_0&=&(V+Q)(U+G_{22}Q)^{-1}=(YD+\wtilD\La D)(XD+\wtilN\La D)^{-1}\\
&=&(Y+\wtilD\La)(X+\wtilN\La)^{-1}\notag
\end{eqnarray}
and
\begin{eqnarray}
\label{comp2}
K_0&=&(\wtilU+QG_{22})^{-1}(\wtilV+Q)
=(\wtilD\wtilX+\wtilD\La N)^{-1}(\wtilD\wtilY+\wtilD\La D)\\
&=&(\wtilX+\La N)^{-1}(\wtilY+\La D).\notag
\end{eqnarray}
Set
\[
Y_0=(Y+\wtilD\La),\quad X_0=(X+\wtilN\La),\quad
\wtilY_0=(\wtilY+\La D),\quad \wtilX_0=(\wtilX+\La N).
\]
Then certainly $\det X_0\not=0$ and $\det \wtilX_0\not=0$,
and we have
\begin{align*}
& \mat{cc}{D&-N\\-\wtilY_0&\wtilX_0}\mat{cc}{X_0&\wtilN\\Y_0&\wtilD}
 = \mat{cc}{D&-N\\-\wtilY-\La D&\wtilX+\La N}
\mat{cc}{X+\wtilN\La&\wtilN\\Y+\wtilD\La&\wtilD}\\
& \qquad = \mat{cc}{I&0\\-\La&I}\begin{bmatrix} D & -N \\ -\widetilde Y & \widetilde X
\end{bmatrix} \begin{bmatrix} X & \widetilde N \\ Y &
\widetilde D \end{bmatrix}\mat{cc}{I&0\\\La&I}
  =\mat{cc}{I&0\\-\La&I}\mat{cc}{I&0\\\La&I}\\
  &\qquad=\mat{cc}{I&0\\0&I}.
\end{align*}
\end{proof}

Since any stabilizing controller for $G$ is also a stabilizing controller for $G_{22}$,
the following corollary is immediate.

%%%%%%%%%%%%%%%%%%%%%%%%%%%%%%%%%%%%%%%%%%%%%%%%%%%%%%%%%%%%%%%%%%%%%%%%
\begin{corollary}\label{C:coprimestab}
Assume that $G\in \BK^{(n_{\cZ}+n_{\cY})\times (n_{\cW}+n_{\cU})}$ is a stabilizable and that $G_{22}$
admits a double coprime factorization. Then any stabilizing controller $K$ of $G$
admits a double coprime factorization.
\end{corollary}

%%%%%%%%%%%%%%%%%%%%%%%%%%%%%%%%%%%%%%%%%%%%%%%%%%%%%%%%%%%%%%%%%%%%%%%%
\begin{lemma}\label{L:coprime}
Assume that $G$ is stabilizable and that $G_{22}$ admits a double coprime
factorization. Then there exists a double coprime factorization
\eqref{coprime} for $G_{22}$ so that $DG_{21}$ and $G_{12}\wtilD$ are stable.
\end{lemma}

\begin{proof}
Let $K$ be a stabilizing controller for $G$. Then $K$ is also a stabilizing
controller for $G_{22}$. Thus, according to Lemma \ref{L:coprimestab}, there
exists a double coprime factorization \eqref{coprime} for $G_{22}$ so that
$K=YX^{-1}=\wtilX\wtilY^{-1}$. Note that \eqref{coprime} implies that
$\sbm{X&\widetilde N\\Y&\widetilde D}\sbm{D&-N\\-\widetilde Y&\widetilde X}=I$.
In particular, $\wtilD\wtilY=YD$ and $\wtilN\wtilX=XN$. Moreover,
from the computations \eqref{comp4} and \eqref{comp5} we see that
\[
(I-G_{22}K)^{-1}=XD\ands (I-KG_{22})^{-1}=\wtilD\wtilX.
\]
Inserting these identities into the formula for $\Th(G,K)$, and using
that $K$ stabilizes $G$, we find that
\[
\Th(G,K)=
\mat{ccc}{G_{11}+G_{12}YDG_{21}&G_{12}\wtilD\wtilX&G_{12}\wtilD\wtilY\\
YDG_{21}&\wtilD\wtilX&\wtilD\wtilY\\XDG_{21}&\wtilN\wtilX&I+\wtilN\wtilY}\text{ is stable}.
\]
In particular $\mat{cc}{G_{12}\wtilD\wtilX&G_{12}\wtilD\wtilY}$ is stable, and thus
\[
\mat{cc}{G_{12}\wtilD\wtilX&G_{12}\wtilD\wtilY}\mat{c}{\wtilN\\-\wtilD}
=G_{12}\wtilD(\wtilX\wtilN-\wtilY\wtilD)=G_{12}\wtilD
\]
is stable. Similarly, since $\sbm{YDG_{21}\\XDG_{21}}$ is stable, we find that
\[
\mat{cc}{-N & D}\mat{c}{YDG_{21}\\XDG_{21}}=(-NY + DX)DG_{21}=DG_{21}
\]
is stable.
\end{proof}

We now present an alternative proof of Corollary \ref{C:G22G} for the
case that $G_{22}$ admits a double coprime factorization.

%%%%%%%%%%%%%%%%%%%%%%%%%%%%%%%%%%%%%%%%%%%%%%%%%%%%%%%%%%%%%%%%%%%%%%%%
\begin{lemma}\label{L:coprimeG22G}
Assume that $G$ is stabilizable and $G_{22}$ admits a double coprime
factorization. Then $K$ stabilizes $G$ if and only if $K$
stabilizes $G_{22}$.
\end{lemma}

\begin{proof}
It was already noted that in case $K$ stabilizes $G$, then $K$ also
stabilizes $G_{22}$. Now assume that $K$ stabilizes $G_{22}$. Let
$Q\in\BK^{n_\cU\times n_\cY}$ so that $K$ is given by \eqref{QtoK1}.
It suffices to show that $Q\in\Omega$, with $\Omega$ defined by \eqref{Omega}.
Since $G$ is stabilizable, it follows from Lemma \ref{L:coprime} that there
exists a double coprime factorization \eqref{coprime} of $G_{22}$ so that
$DG_{21}$ and $G_{12}\wtilD$ are stable. According to Lemma
\ref{L:free-paramYK}, $Q=\widetilde D\La D$ for some $\La\in \BA^{n_\cU\times n_\cY}$.
It follows that
\begin{eqnarray*}
\mat{c}{G_{12}\\I\\G_{22}}Q\mat{ccc}{G_{21}&G_{22}&I}
&=&\mat{c}{G_{12}\widetilde D\\\widetilde D\\G_{22}\widetilde D}\La\mat{ccc}{DG_{21}&DG_{22}&D}\\
&=&\mat{c}{G_{12}\widetilde D\\\widetilde D\\\widetilde N}\La\mat{ccc}{DG_{21}&N&D}
\end{eqnarray*}
is stable. Hence $Q\in\Omega$.
\end{proof}

Combining the results from the Lemmas \ref{L:free-paramYK}, \ref{L:coprimestab}
and \ref{L:coprimeG22G} with Theorem \ref{T:stab2} and the
computations \eqref{comp1} and \eqref{comp2} from the proof of
Lemma \ref{L:coprimestab} we obtain the Youla-Ku\v cera
parametrization of all stabilizing controllers.

%%%%%%%%%%%%%%%%%%%%%%%%%%%%%%%%%%%%%%%%%%%%%%%%%%%%%%%%%%%%%%%%%%%%%%%%
\begin{theorem}\label{T:stab4}
Assume that $G\in \BK^{(n_{\cZ}+n_{\cY})\times (n_{\cW}+n_{\cU})}$ is stabilizable and that $G_{22}$
admits a double coprime factorization \eqref{coprime}. Then the set of all
stabilizing controllers is given by
\[
K=(Y+\wtilD\La)(X+\wtilN\La)^{-1}=(\wtilX+\La N)^{-1}(\wtilY+\La D),
\]
where $\La$ is a free stable parameter from $\BA^{n_\cU\times n_\cY}$ such that
$\det (X + \wtilN\La)\not=0$ or equivalently  $\det(\wtilX+\La N)\not=0$.
\end{theorem}

\subsection{The standard $H^\infty$-problem reduced to model matching.}
\label{S:HinfToMM}

We now consider the $H^\infty$-problem for a plant $G= \sbm{ G_{11} & G_{12} \\ G_{21} & G_{22}}
\colon \cW \oplus \cU \to \cZ \oplus \cY$, i.e., we seek a controller $K \colon \cY \to \cU$
so that not only $\Theta(G,K)$ in \eqref{ThetaGK} is stable, but also
\[
\|G_{11} + G_{12}K (I -G_{22}K)^{-1} G_{21}\|_{op}\leq 1.
\]

Assume that the plant $G$ is stabilizable, and that $K_* \colon \cY \to \cU$
stabilizes $G$. Define $U_*$, $V_*$, $\wtilU_*$ and $\wtilV_*$ as in Theorem \ref{T:stab2}.
We then know that all stabilizing controllers of $G$ are given by
\[
K=(V_* + Q)(U_* + G_{22}Q)^{-1}=(\widetilde U_{*}  + Q G_{22})^{-1}(\wtilV_*+Q)
\]
where $Q \in \BK^{n_{\cU} \times n_{\cY}}$ is any element of $\Omega$ in \eqref{Omega}.
We can then express the transfer matrices $U$ and $V$ in \eqref{UV}
in terms of $Q$ as follows:
\begin{eqnarray*}
U&=&(I-G_{22}K)^{-1}
=(I-G_{22}(V_*-Q)(U_*-G_{22}Q)^{-1})^{-1}\\
&=&(U_*-G_{22}Q)(U_*-G_{22}Q-G_{22}(V_*-Q))^{-1}\\
&=&(U_*-G_{22}Q)(U_*-G_{22}V_*)^{-1}\\
&=&(U_*-G_{22}Q),
\end{eqnarray*}
where we used that $U_*-G_{22}V_*=I$, and
\[
V=KU=V_*-Q.
\]
Similar computations provide the formulas
\[
\wtilU=\wtilU_*+QG_{22}\ands \wtilV=\wtilV_*+Q
\]
for the transfer matrices $\wtilU$ and $\wtilV$ in \eqref{tilUtilV}.
Now recall that $\Theta(G,K)$ can be expressed in terms of $U$ and $V$
as in \eqref{Theta2}. It then follows that left upper block in $\Theta(G,K)$
is equal to
\begin{eqnarray}\label{comp3}
G_{11} + G_{12}K (I -G_{22}K)^{-1} G_{21}&=&G_{11}+G_{12}VG_{21}\\
&=&G_{11}+G_{12}V_*G_{21}-G_{12}QG_{21}.\notag
\end{eqnarray}
The fact that $K_*$ stabilizes $G$ implies that
$\wtilG_{11}:=G_{11}+G_{12}V_*G_{21}$ is stable, and thus $G_{12}QG_{21}$
is stable as well. We are now close to a reformulation of the $H^\infty$-problem
as a model matching problem. However, to really formulate it as a
model matching problem, we need to apply the change of design
parameter $Q\mapsto \La$
defined in Lemma \ref{L:free-param}, or Lemma \ref{L:free-paramYK} in case $G_{22}$
admits a double coprime factorization.  The next two results extend the idea
of Theorem \ref{T:1.2} to this more general setting.

%%%%%%%%%%%%%%%%%%%%%%%%%%%%%%%%%%%%%%%%%%%%%%%%%%%%%%%%%%%%%%%%%%%%%%%%
\begin{theorem}\label{T:HinftoMM1}
Assume that $G\in \BK^{(n_{\cZ}+n_{\cY})\times (n_{\cW}+n_{\cU})}$ is stabilizable and
let $K\in \BK^{n_\cU\times n_\cY}$. Then $K$ is a solution to the standard
$H^{\infty}$ problem for $G$ if and only if
\[
K=(V_* + Q)(U_* + G_{22}Q)^{-1}=(\widetilde U_{*}  + Q G_{22})^{-1}(\wtilV_*+Q)
\]
with $Q=\wtilL\La L$, where $\wtilL$ and $L$ are defined by \eqref{tildeLL},
so that $\det(U_* + G_{22}Q)\not=0$, or equivalently
$\det(\widetilde U_{*}  + Q G_{22})\not=0$,
and $\Lambda \in\BA^{(n_{\cU} +n_{\cY}) \times (n_{\cU} + n_{\cY})}$  is any
solution to the model matching problem for $\wtilG_{11}$, $\wtilG_{12}$ and
$\wtilG_{21}$ defined by
\[
\wtilG_{11}:=G_{11}+G_{12}V_*G_{21},\quad
\wtilG_{12}:=G_{12}\wtilL,\quad
\wtilG_{21}:=LG_{21},
\]
i.e., so that
\[
\|\wtilG_{11}+\wtilG_{12}\La \wtilG_{21}\|_{op}\leq 1.
\]
\end{theorem}

%%%%%%%%%%%%%%%%%%%%%%%%%%%%%%%%%%%%
\begin{proof}
The statement essentially follows from Theorem \ref{T:stab3} and the
computation \eqref{comp3} except that we need to verify that the
functions $\wtilG_{11}$, $\wtilG_{12}$ and $\wtilG_{21}$ satisfy the
conditions to be data for a model matching problem, that is, they
should be stable. It was already observed that $\wtilG_{11}$ is
stable. The fact that $\wtilG_{12}$ and $\wtilG_{21}$ are stable
was shown in the proof of Lemma \ref{L:free-param}.
\end{proof}
%%%%%%%%%%%%%%%%%%%%%%%%%%%%%%%%%%%%

We have a similar result in case $G_{22}$ admits a double coprime factorization.

%%%%%%%%%%%%%%%%%%%%%%%%%%%%%%%%%%%%%%%%%%%%%%%%%%%%%%%%%%%%%%%%%%%%%%%%
\begin{theorem}\label{T:HinftoMM2}
Assume that $G\in \BK^{(n_{\cZ}+n_{\cY})\times (n_{\cW}+n_{\cU})}$ is stabilizable and that
$G_{22}$ admits a double coprime factorization \eqref{coprime}.
Let $K\in \BK^{n_\cY\times n_\cU}$. Then $K$ is a solution to the standard
$H^{\infty}$ problem for $G$ if and only if
\[
K=(Y + \wtilD\La)(X + \wtilN\La)^{-1}=(\wtilX+\La N)^{-1}(\wtilY+\La D),
\]
where $\La\in\BA^{n_\cU\times n_\cY}$ so that $\det(X + \wtilN\La)\not=0$,
or equivalently $\det(\wtilX+\La N)\not=0$, is any solution to the model matching
problem for $\wtilG_{11}$, $\wtilG_{12}$ and $\wtilG_{21}$ defined by
\[
\wtilG_{11}:=G_{11}+G_{12} Y D G_{21},\quad
\wtilG_{12}:=G_{12}\wtilD,\quad
\wtilG_{21}:=DG_{21},
\]
i.e., so that
\[
\|\wtilG_{11}+\wtilG_{12}\La \wtilG_{21}\|_{op}\leq 1.
\]
\end{theorem}

%%%%%%%%%%%%%%%%%%%%%%%%%%%%%%%%%%%%
\begin{proof}
The same arguments apply as in the proof of Theorem \ref{T:HinftoMM1}, except that
in this case Lemma \ref{L:coprime} should be used to show that $\wtilG_{12}$ and
$\wtilG_{21}$ are stable.
\end{proof}
%%%%%%%%%%%%%%%%%%%%%%%%%%%%%%%%%%%%

\subsection{Notes}  \label{S:Notes1}
The development in Section \ref{S:stab-w/oCP} on the parametrization
of stabilizing controllers without recourse to a double coprime
factorization of $G_{22}$ is based on the exposition of Quadrat
\cite{Q-elementary}. It was already observed by Zames-Francis
\cite{ZF} that $Q = K (I - G_{22}K)^{-1}$ can be used as a free
stable design parameter in case $G_{22}$ is itself already stable; in
case $G_{22}$ is not stable, $Q$ is subject to some additional
interpolation conditions.  The results of \cite{Q-elementary} is an
adaptation of this observation to the general ring-theoretic setup.
The more theoretical papers \cite{Q-YKI, Q-YKII} give  module-theoretic
interpretations for the structure associated with internal
stabilizability.  In particular, it comes out that every matrix
transfer function $G_{22}$ with entries in ${\mathbb K}$ has a
double-coprime factorization if and only if ${\mathbb A}$ is a {\em
Bezout domain}, i.e., every finitely generated ideal in ${\mathbb
A}$ is principal; this recovers a result already appearing in the
book of Vidyasagar \cite{Vid}.  A new result which came out of this
module-theoretic interpretation was that internal stabilizability of a plant
$G_{22}$ is equivalent to the existence of a double-coprime
factorization for $G_{22}$ exactly when the ring ${\mathbb A}$ is
{\em projective-free}, i.e., every submodule of a finitely generated
free module over ${\mathbb A}$ must itself be free.  This gives an
explanation for the earlier result of Smith \cite{Smith} that this
phenomenon holds for the case where ${\mathbb A}$ is equal
$H^{\infty}$ over the unit disk or right-half plane.
%The result of Corollary \ref{C:G22G}, see also Lemma \ref{L:coprimeG22G},
%appears \cite{Q-Leuven} and proves a conjecture of Lin posed in \cite{Lin3}.

Earlier less complete results concerning parametrization of the set of
stabilizing controllers without the assumption of a coprime factorization
were obtained by Mori \cite{Mori02} and Sule \cite{Sule}.  Mori
\cite{Mori04} also showed that the internal-stabilization problem can
be reduced to model matching form for the general case where the
plant has the full $2 \times 2$-block structure $G = \sbm{ G_{11} &
G_{12} \\ G_{21} & G_{22}}$.

Lemma \ref{L:coprimeG22G} for the classical case is Theorem 2 on page 35 in
\cite{Francis}.
The proof there relies in a careful analysis of signal-flow diagrams; we believe
that our proof is more direct.

  \section{Feedback control for linear time-invariant multidimensional
  systems} \label{S:ND}\setcounter{equation}{0}

  \subsection{Multivariable frequency-domain formulation}
  \label{S:ND-freq}

  The most obvious multivariable analogue of the classical
  single-variable setting considered in the book of Francis
  \cite{Francis} is as follows. We take the underlying field
   to be the complex numbers ${\mathbb C}$; in the engineering applications, one
  usually requires that the underlying field be the reals ${\mathbb
  R}$, but this can often be incorporated at the end by using
  the characterization of real rational functions as being those
  complex rational functions which are invariant under the
  conjugation operator $s(z) \mapsto \overline{s(\overline{z})}$.
   We let ${\mathbb D}^{d} = \{ z =
  (z_{1}, \dots, z_d) \colon |z_{k}| < 1\}$ be the unit polydisk in the
  $d$-dimensional complex space ${\mathbb C}^{d}$ and we take our
  ring ${\mathbb A}$ of stable plants to be the
  ring ${\mathbb C}(z)_{s}$ of all rational functions $s(z) =
  \frac{p(z)}{q(z)}$ in $d$ variables (thus, $p$ and $q$ are
  polynomials in the $d$ variables $z_{1}, \dots, z_{d}$   where we
  set $z = (z_{1}, \dots, z_{d})$) such that $s(z)$ is bounded on the
  polydisk ${\mathbb D}^{d}$. The ring ${\mathbb C}[z]$ of polynomials in $d$ variables is a
  unique factorization domain so we may assume that $p$ and $q$ have
  no common factor (i.e., that $p$ and $q$ are {\em relatively coprime})
  in the fractional representation $s = \frac{p}{q}$
  for any element of ${\mathbb C}(z_{1}, \dots, z_{d})$. Unlike in
  the single-variable case, for the case $d>1$ it can happen that $p$
  and $q$ have common zeros in ${\mathbb C}^{d}$ even when they are
  coprime in ${\mathbb C}[z]$ (see \cite{Youla} for an early analysis
  of the resulting distinct notions of coprimeness).  It turns out
  that for $d \ge 3$, the ring ${\mathbb C}(z)_{s}$ is difficult to work with
  since the denominator $q$ for a stable ring
  element depends in a tricky way on the numerator $p$: if $s
  \in {\mathbb C}(z)_{s}$ has coprime fractional representation $s =
  \frac{p}{q}$, while it is the case that
  necessarily $q$ has no zeros in the open polydisk ${\mathbb
  D}^{d}$, it can happen that the zero variety of $q$ touches the boundary
  $\partial {\mathbb D}^{d}$ as long as the zero variety of $p$ also
  touches the same points on the boundary in such a way that the
  quotient $s = \frac{p}{q}$ remains bounded on ${\mathbb D}^{d}$.
  Note that at such a boundary point $\zeta$, the quotient $s = p/q$
  has no well-defined value.  In the engineering literature (see
  e.g.~\cite{Bose, SRP, KharMun}), such a point is known as a {\em
  nonessential singularity of the second kind}.

  To avoid this difficulty, Lin \cite{Lin1, Lin2} introduced the ring
  ${\mathbb C}(z)_{ss}$ of {\em structured stable} rational
  functions, i.e., rational functions $s \in {\mathbb C}(z)$ so that
  the denominator $q$ in any coprime fractional representation $s =
  \frac{p}{q}$ for $s$ has no zeros in the closed polydisk
  $\overline{\mathbb D}^{d}$.  According to the result of
  Kharitonov--Torres-Mu\~noz \cite{KharMun},  whenever $s =
  \frac{p}{q} \in {\mathbb C}(z)_{s}$ is stable in the first
  (non-structured) sense, an arbitrarily small perturbation of
  the coefficients of $q$ may lead to the perturbed $q$ having zeros in
  the open polydisk ${\mathbb D}^{d}$ resulting in the perturbed
  version $s = \frac{p}{q}$ of $s$ being unstable; this phenomenon does
does not occur for $s\in \BC(z)_{ss}$, and thus {\em
  structured stable} can be viewed just as a robust version of {\em
  stable} (in the unstructured sense).  Hence one can argue that
  structured stability is the more desirable property from an engineering
perspective. In the
  application to delay systems using the systems-over-rings approach
  \cite{BST, KharSon, KKT}, on the other hand, it is the collection
  ${\mathbb C}(z)_{ss}$ of structurally
  stable rational functions  which comes up in the first place.

  As the ring ${\mathbb A} = {\mathbb C}(z)_{ss}$ is a commutative
  integral domain, we can apply the results of Section \ref{S:ring}
  to this particular setting.  It was proved in connection with work
  on systems-over-rings rather than multidimensional systems (see
  \cite{BST, KKT}) that the
  ring ${\mathbb C}(z)_{ss}$ is {\em projective-free}.  As pointed
  out in the notes of Section \ref{S:ring} above, it follows that
  stabilizability of $G_{22}$ is
  equivalent to the existence of a double coprime factorization for
  the plant $G_{22}$ (see \cite{Q-Leuven}), thereby settling a
  conjecture of Lin \cite{Lin1, Lin2, Lin3}.  We summarize these results as follows.

  \begin{theorem} \label{T:Lin} Suppose that we are given a system $G = \left[
      \begin{smallmatrix} G_{11} & G_{12} \\ G_{21} & G_{22}
      \end{smallmatrix} \right]$ over the quotient field
      $Q({\mathbb C}(z)_{ss})$ of the ring ${\mathbb C}(z)_{ss}$ of
      structurally stable rational functions in $d$ variables.  If
      there exists a controller $K = Y X^{-1} = \widetilde X^{-1}
      \widetilde Y$ which internally stabilizes $G$, then $G_{22}$
      has a double coprime factorization and all internally
      stabilizing controllers $K$ for $G$ are given by the
      Youla-Ku\v cera parametrization.
      \end{theorem}

Following Subsection \ref{S:HinfToMM}, the Youla-Ku\v cera
parametrization can then be used to rewrite the $H^{\infty}$-problem in
the form of a model-matching problem:
      {\em Given $T_{1}, T_{2}, T_{3}$ equal to matrices over
      ${\mathbb C}(z)_{ss}$ of respective sizes $n_{\cZ} \times
      n_{\cW}$, $n_{\cW} \times n_{\cU}$ and $n_{\cY} \times
      n_{\cW}$, find a matrix $\Lambda$ over ${\mathbb C}(z)_{ss}$ of
      size $n_{\cU} \times n_{\cY}$ so that the affine expression $S$
      given by
      \begin{equation}   \label{NDmod-match}
        S =  T_{1} + T_{2} \Lambda T_{3}
      \end{equation}
      has supremum norm at most 1, i.e., $\|S\|_{\infty}=\max\{\|S(z)\|\colon z\in\ov{\BD}^d\}
      \leq1$.}

  For mathematical convenience we shall now widen the class of admissible
       solutions and allow $\Lambda_{1}, \dots, \Lambda_{J}$ to be
       in the   algebra $H^{\infty}({\mathbb D}^{d})$ of bounded analytic
       functions on ${\mathbb D}^{d}$. The unit ball of
       $H^{\infty}({\mathbb D}^{d})$ is the set of  all holomorphic functions $S$
       mapping the polydisk ${\mathbb D}^{d}$ into the closed
       single-variable unit disk ${\mathbb D} \subset {\mathbb C}$; we
       denote this space by ${\mathcal S}_{d}$, the {\em $d$-variable
       Schur class}. While $T_{1}$, $T_{2}$ and $T_{3}$ are assumed
       to be in ${\mathbb C}(z)_{ss}$, we allow $\Lambda$ in
       \eqref{NDmod-match} to be in $H^{\infty}({\mathbb D}^{d})$.

      Just as in the classical one-variable case, it is possible to
      give the model-matching form \eqref{NDmod-match} an
      interpolation  interpretation, at least for special cases (see
      \cite{HeltonACC, HeltonIEEE, BM}).  One such case is where
     $n_{\cW} = n_{\cZ} = n_{\cY} = 1$ while $n_{\cU}
      = J$.  Then $T_{1}$ and $T_{3}$ are scalar while $T_{2} = \left[
      \begin{smallmatrix} T_{2,1} & \cdots & T_{2,J} \end{smallmatrix}
      \right]$ is a row.  Assume in addition that $T_{3} = 1$.
      Then the model-matching form \eqref{NDmod-match} collapses
      to
\begin{equation}  \label{scalarNP-modmatch}
      S = T_{1} + T_{21} \Lambda_{1} + \cdots + T_{2J} \Lambda_{J}
  \end{equation}
  where $\Lambda_{1}, \dots \Lambda_{J}$ are $J$ free stable scalar
  functions.  Under the assumption that the intersection of the zero
  varieties of $T_{2,1}, \dots, T_{2,J}$ within the closed polydisk
  $\overline{\mathbb D}^{d}$ consists of finitely many (say $N$)  points
  $$
  z_{1} = (z_{1,1}, \dots, z_{1,d}), \cdots, z_{N} = (z_{N,1}, \dots,
  z_{N,d})
  $$
  and if we let $w_{1}, \dots, w_{N}$ be the values of $T_{1}$ at
  these points
  $$
    w_{1} = T_{1}(z_{1}), \dots, w_{N} = T_{1}(z_{N}),
  $$
  then it is not hard to see that a function $S \in  {\mathbb
  C}(z)_{ss}$ has the form \eqref{scalarNP-modmatch} if and only if
  it satisfies the interpolation conditions
  \begin{equation}   \label{ND-NevPick}
  S(z_{i}) = w_{i} \text{ for } i = 1, \dots, N.
  \end{equation}
  In this case the model-matching problem thus becomes the following
  finite-point Nevanlinna-Pick interpolation problem over ${\mathbb
  D}^{d}$:  {\em find $S \in {\mathbb
  C}(z)_{ss}$ subject to $|S(z)| \le 1$ for all $z \in {\mathbb
  D}^{d}$ which satisfies the interpolation conditions
  \eqref{ND-NevPick}.} Then the $d$-variable $H^{\infty}$-Model-Matching problem
  becomes: {\em find $S \in {\mathcal S}_{d}$ so that $S(z_{1}) =
  w_{1}$ for $i = 1, \dots, N$.}

  A second case (see \cite{BM}) where the polydisk Model-Matching
  Problem can be reduced to an interpolation problem is the case where $T_{2}$ and $T_{3}$
  are square (so $n_{\cZ} = n_{\cU}$ and $n_{\cY} = n_{\cW}$) with invertible values
  on the distinguished boundary of the polydisk; under these
  assumptions it is shown in \cite{BM} (see Theorem 3.5 there) how the model-matching problem
  is equivalent to a {\em bitangential Nevanlinna-Pick interpolation problem
  along a subvariety}, i.e.,
  bitangential interpolation conditions are specified along all
  points of a codimension-1 subvariety of ${\mathbb D}^{d}$ (namely, the union of
  the zero sets of $\det T_{2}$ and $\det T_{3}$ intersected with
  ${\mathbb D}^{d}$).  For $d=1$, codimension-1 subvarieties are
  isolated points in the unit disk; thus the codimension-1
  interpolation problem is a direct
  generalization of the bitangential Nevanlinna-Pick interpolation problem
  studied in \cite{BGR, Dym, FF}.
  However for the
  case where the number of variables $d$ is at least 3, there
  is no theory with results parallel to those of the classical case.

  Nevertheless, if we change the problem somewhat there is a theory
  parallel to the classical case.  To formulate this adjustment, we
  define the {\em $d$-variable Schur-Agler class} $\mathcal{S A}_{d}$
to consist of those functions $S$ analytic on the polydisk
  for which the operator $S(X_{1}, \dots, X_{d})$ has norm at most 1
  for any collection $X_{1}, \dots, X_{d}$ of $d$ commuting strict contraction operators
  on a separable Hilbert space ${\mathcal K}$; here $S(X_{1}, \dots,
  X_{d})$ can be defined via the formal power series for $S$:
  $$ S(X_{1}, \dots, X_{d}) = \sum_{n \in {\mathbb Z}^{d}_{+}} s_{n}
  X^{n},\quad \text{if }S(z) = \sum_{n \in {\mathbb Z}^{d}_{+}} s_{n}
  z^{n}
  $$
  where we use the standard multivariable notation
  $$
  n = (n_{1}, \dots, n_{d}) \in {\mathbb Z}^{d}_{+}, \quad
  X^{n} = X_{1}^{n_{1}} \cdots X_{d}^{n_{d}} \text{ and }
   z^{n} = z_{1}^{n_{1}} \cdots z_{d}^{n_{d}}.
   $$
  For the cases $d=1,2$, it turns out, as a consequence of the von
  Neumann inequality or the Sz.-Nagy dilation theorem for $d=1$ and
  of  the And\^o dilation theorem \cite{Ando} for $d=2$ (see \cite{Paulsen,
  BSV} for a full discussion), that the Schur-Agler class
  $\mathcal{S A}_{d}$ and the Schur class $\mathcal{S}_{d}$ coincide,
  while, due to an explicit example of Varopoulos, the inclusion
  $\mathcal {S A}_{d} \subset \mathcal{S}_{d}$  is strict for $d \ge
  3$.

  There is a result due originally to Agler \cite{Agler-unpublished}
  and developed and refined in a number of directions since (see
  \cite{AgMcC99, BT} and \cite{AgMcCbook} for an overview)
  which parallels the one-variable case; for the case of a simple
  set of interpolation conditions  \eqref{ND-NevPick} the result is
  as follows: {\em there exists a function $S$ in the Schur-Agler
  class $\mathcal{ S A}_{d}$ which satisfies the set of interpolation
  conditions $S(z_{i}) = w_{i}$ for $i = 1, \dots, N$ if and only if
  there exist $d$ positive semidefinite matrices ${\mathbb P}^{(1)},
  \dots, {\mathbb P}^{(d)}$ of size $N \times N$ so that}
  $$
    1 - w_{i} \overline{w_{j}} = \sum_{k=1}^{d} (1 - z_{i,k}
    \overline{z_{j,k}} ) {\mathbb P}^{(k)}_{i,j}.
  $$
For the case $d=1$, the Pick matrix
$\BP=\sbm{\frac{1-w_i\overline{w}_j}{1-z_i\overline{z}_j}}_{i,j=1}^N$
is the unique solution of this equation, and we recover the classical criterion
$\BP\geq0$ for the existence of solutions to the Nevanlinna-Pick problem.
There is a later realization result of Agler \cite{Agler-Hellinger}
(see also \cite{AgMcC99, BT}):  a given holomorphic function $S$ is in the
Schur-Agler class $\mathcal{ S A}_{d}(\cL(\cU, \cY))$ if and only if
$S$ has a contractive Givone-Roesser realization:
$S(z) = D + C (I - Z(z) A)^{-1} Z(z) B$ where $\sbm{ A & B \\ C & D }
\colon (\oplus_{k=1}^{d} \cX_{k} \oplus \cU) \to (\oplus_{k=1}^{d}
\cX_{k} \oplus \cY)$ is contractive with $Z(z) = \sbm{
z_{1}I_{\cX_{1}} & & \\ & \ddots & \\ & & z_{d} I_{\cX_{d}} }$.

  Direct application of the Agler result to the
  bitangential Nevanlinna-Pick interpolation problem along a subvariety, however, gives a
    solution criterion involving an infinite
    Linear Matrix Inequality (where the unknown matrices have
    infinitely many rows and columns
    indexed by the points of the interpolation-node
    subvariety)---see \cite[Theorem 4.1]{BM}.
    Alternatively, one can use the polydisk Commutant Lifting Theorem
    from \cite{BLTT} to get a solution criterion involving a Linear
    Operator Inequality \cite[Theorem 5.2]{BM}.  Without further
    massaging, either approach
    is computationally unattractive; this is in contrast with the
    state-space approach discussed below. In that setting there exists
computable sufficient conditions, in terms of a pair of LMIs and a coupling
condition, that in general are only sufficient, unless one works with a more
conservative notion of stability and performance.

  \subsection{Multidimensional state-space formulation}
  \label{S:ND-statespace}

The starting point in this subsection is a quadruple $\{A,B,C,D\}$
consisting of operators $A$, $B$, $C$ and $D$  so that
$\sbm{A&B\\C&D}:\sbm{\cX\\\cW\oplus\cU}\to\sbm{\cX\\\cZ\oplus\cY}$
and a partitioning $\cX=\cX_1\oplus\cdots\oplus\cX_d$ of the space $\cX$.
Associate with such a quadruple $\{A,B,C,D\}$ is a linear state-space
system $\Si$ of Givone-Roesser type (see \cite{GR}) that evolves over $\BZ_+^d$ and is given
by the system of equations
\begin{equation}\label{GRmodel}
\Si:=\left\{\begin{array}{c}\sbm{x_1(n+e_1)\\\vdots\\x_d(n+e_d)}
=A\sbm{x_1(n)\\\vdots\\x_d(n)}+Bu(n)\\[.2cm]
y(n)=Cx(n)+Du(n)\end{array}\right.\qquad (n\in\BZ_+^d),
\end{equation}
with initial conditions a specification of the state values
$x_k(\sum_{j\not= k}t_j e_j)$ for $t = (t_{1}, \dots, t_{d}) \in
{\mathbb Z}^{d}_{+}$ subject to $t_{k}=0$ where $k=1,\ldots,d$.
Here $e_k$ stands for the $k$-th unit vector in $\BC^d$ and
$x(n)=\sbm{x_1(n)\\\vdots\\x_d(n)}$. We call $\cX$ the {\em state-space} and $A$ the
{\em state operator}. Moreover, the block operator matrix $\sbm{A&B\\C&D}$ is referred
to as the {\em system matrix}.

Following \cite{Kaczorek}, the Givone-Roesser system \eqref{GRmodel} is said to be
{\em asymptotically stable} in case, for zero input $u(n)=0$ for $n\in\BZ_+^d$ and initial
conditions with the property
\[
\sup_{t \in {\mathbb Z}^{d}_{+} \colon t_{k}=0}
\|x_k(\sum_{j=1}^{d} t_j e_j)\|<\infty\text{ for }k=1,\ldots,d,
\]
the state sequence $x$ satisfies
\[
\sup_{n\in\BZ_+^d}\|x(n)\|<\infty\ands
\lim_{n\to\infty}\|x(n)\|=0,
\]
where $n\to\infty$ is to be interpreted as $\min\{n_1,\ldots n_d\}\to\infty$ when
$n=(n_1,\ldots,n_d)\in\BZ_+^d$.

With the Givone-Roesser system \eqref{GRmodel} we associate the transfer
function $G(z)$ given by
\begin{equation}\label{transfer}
G(z)=D+C(I-Z(z)A)^{-1}Z(z)B,
\end{equation}
defined al least for $z\in\BC^d$ with $\|z\|$ sufficiently small, where
\begin{equation}\label{Z(z)}
Z(z)=\mat{ccc}{z_1I_{\cX_1}&&\\&\ddots&\\&&z_dI_{\cX_d}}\quad (z\in\BC^d).
\end{equation}
We then say that $\{A,B,C,D\}$ is a {\em (state-space) realization}
for the function $G$, or if $G$ is not specified, just refer to
$\{A,B,C,D\}$ as a realization. The realization $\{A,B,C,D\}$, or
just the operator $A$, is said to be {\em Hautus-stable} in case the
pencil $I-Z(z)A$ is invertible on the closed polydisk $\ov{\BD}^d$.

Here we only consider the case that $\cX$ is finite dimensional;
then the entries of the transfer function $G$ are in the quotient
field $Q(\BC(z)_{ss})$ of ${\mathbb C}(z)_{ss}$ and are analytic at 0, and it is
straightforward to see that $G$ is
structurally stable in case $G$ admits a Hautus-stable realization.
For the case $d=2$ it has been asserted in the literature
\cite[Theorem 4.8]{Kaczorek} that asymptotic stability and Hautus
stability are equivalent; presumably this assertion continues to hold
for general $d \ge 1$ but we do not go into details here.

Given a realization $\{A, B, C, D\}$ where the decomposition $\cX =
\cX_{1} \oplus \cdots \oplus \cX_{d}$ is understood, our main interest
will be in Hautus-stability; hence we shall say
simply that $A$ is {\em stable} rather than {\em Hautus-stable}.

As before we consider controllers $K$ in $Q(\BC(z)_{ss})$ of size $n_\cY\times n_\cU$
that we also assume to be given by a state-space realization
\begin{equation}\label{statecontroller}
K(z)=D_K+C_K(I-Z_K(z)A_K)^{-1}Z_K(z)B_K
\end{equation}
with system matrix $\sbm{A_K&B_K\\C_K&D_K}:\sbm{\cX_K\\\cY}\to\sbm{\cX_K\\\cU}$,
a decomposition of the state-space $\cX_K=\cX_{1,K}\oplus\cdots\oplus\cX_{d,K}$
and $Z_K(z)$ defined analogous to $Z(z)$ but with respect to the decomposition
of $\cX_K$.  We now further specify the matrices $B$, $C$ and $D$ from the
realization $\{A,B,C,D\}$ as
\begin{equation}\label{decomps}
B=\mat{cc}{B_1&B_2},\quad C=\mat{c}{C_1\\C_2},\quad D=\mat{cc}{D_{11}&D_{12}\\D_{21}&D_{22}}
\end{equation}
compatible with the decompositions $\cZ\oplus\cY$ and $\cW\oplus\cU$.
We can then form the closed loop system $G_{cl}=\Si(G,K)$ of
the two transfer functions. The closed loop system $G_{cl}=\Si(G,K)$
corresponds to the feedback connection
\[
\mat{ccc}{A&B_1&B_2\\C_1&D_{11}&D_{12}\\C_2&D_{21}&D_{22}}
\mat{c}{x\\w\\u}\to\mat{c}{\wtil{x}\\z\\y},\quad
\mat{cc}{A_K&B_K\\C_K&D_K}:\mat{c}{x_K\\u_K}\to\mat{c}{\wtil{x}_K\\y_K}
\]
subject to
\[
x=Z(z)\wtil{x},\quad x_K=Z_K(z)\wtil{x}_K,\quad u_K=y\ands y_K=u.
\]
This feedback loop is well-posed exactly when $I-D_{22}D_K$ is invertible. Since,
under the assumption of well posedness, one can always arrange via a change of
variables that $D_{22}=0$ (cf., \cite{IS}), we shall assume that
$D_{22}=0$ for the remainder of the paper. In that case
well-posedness is automatic and the closed loop system $G_{cl}$ admits a state-space realization
\begin{equation}\label{cltransfer}
G_{cl}(z)=D_{cl}+C_{cl}(I-Z_{cl}(z)A_{cl})^{-1}Z_{cl}(z)C_{cl}
\end{equation}
with system matrix
\begin{equation}  \label{clsysmat}
   \begin{bmatrix}  A_{cl} & B_{cl} \\ C_{cl} & D_{cl} \end{bmatrix} =
       \left[ \begin{array}{cc|c}
       A+B_{2} D_{K}C_{2} & B_{2} C_{K} & B_{1} + B_{2} D_{K} D_{21} \\
       B_{K} C_{2} & A_{K} & B_{K} D_{21} \\
       \hline C_{1} + D_{12} D_{K} C_{2} & D_{12} C_{K} & D_{11} +
       D_{12} D_{K} D_{21} \end{array} \right]
\end{equation}
and
\[
Z_{cl}(z)=\mat{cc}{Z(z)&0\\0&Z_{K}(z)}\quad(z\in\BC^d).
\]

The {\em state-space (internal) stabilizability problem} then is: Given the
realization $\{A,B,C,D\}$ find a compatible controller $K$ with
realization $\{A_K,B_K,C_K,D_K\}$ so that the closed-loop
realization $\{A_{cl},B_{cl},C_{cl},D_{cl}\}$ is stable,
i.e., so that $I-Z_{cl}(z)A_{cl}$ is invertible on the closed
polydisk $\ov{\BD}^d$. We also consider the {\em strict state-space
$H^\infty$-problem}: Given the realization $\{A,B,C,D\}$, find a
compatible controller $K$ with realization $\{A_K,B_K,C_K,D_K\}$ so
that the closed loop realization $\{A_{cl},B_{cl},C_{cl},D_{cl}\}$
is stable and the closed-loop system $G_{cl}$ satisfies
$\|G_{cl}(z)\|< 1$ for all $z\in\BD^d$.\bigskip

\paragraph{\bf State-space stabilizability}

In the fractional representation setting of Section \ref{S:ring} it
took quite some effort to derive the result: ``If $G$ is
stabilizable, then $K$ stabilizes $G$ if and only if $K$ stabilizes
$G_{22}$'' (see Corollary \ref{C:G22G} and Lemma
\ref{L:coprimeG22G}). For the state-space stabilizability problem
this result is obvious, and what is more, one can drop the
assumption that $G$ needs to be stabilizable. Indeed, $G_{22}$
admits the realization $\{A,B_2,C_2,0\}$ (assuming $D_{22}=0$), so
that the closed-loop realization for $\Si(G_{22},K)$ is equal to
$\{A_{cl},0,0,0\}$. In particular, both closed-loop realizations
have the same state operator $A_{cl}$, and thus $K$ with realization
$\{A_K,B_K,C_K,D_K\}$ stabilizes $G_{22}$ if and only if $K$
stabilizes $G$, without any assumption on the stabilizability
of $G$.

The state-space stabilizability problem does not have a clean solution;
To discuss the partial results which exist, we first introduce some terminology.

Let $\{A,B,C,D\}$ be a given realization as above with decomposition
of $B$, $C$ and $D$ as in \eqref{decomps}. The {\em Givone-Roesser
output pair} $\{C_2,A\}$ is said to be {\em Hautus-detectable} if
the block-column matrix $\sbm{I-Z(z) A\\C_2}$ is of maximal rank
$n_{\cX}$ (i.e., is left invertible) for all $z$ in the closed
polydisk $\overline{\mathbb D}^{d}$. We say that $\{C_2,A\}$ is {\em
operator-detectable} in case there exists an output-injection
operator $L \colon \cY \to \cX$ so that $A + LC_2$ is stable.
Dually, the {\em Givone-Roesser input pair} $\{A, B_2\}$ is called
{\em Hautus-stabilizable} if it is the case that the block-row
matrix $\begin{bmatrix} I - A Z(z) & B_2 \end{bmatrix}$ has maximal
rank $n_{\cX}$ (i.e., is right invertible) for all $z \in
\overline{\mathbb D}^{d}$, and {\em operator-stabilizable} if there
is a state-feedback operator $F \colon \cX \to \cU$ so that $A + B_2
F$ is stable.  Notice that both Hautus-detectability and
operator-detectability for the pair $(C,A)$ reduce to stability of
$A$ in case $C=0$.  A similar remark applies to stabilizability for
an input pair $(A, B)$.

We will introduce yet another notion of detectability and stabilizability shortly, but
in order to do this we need a stronger notion of stability. We first define $\cD$ to be
the set
\begin{equation}\label{Zcom}
\cD=\left\{\sbm{X_1&&\\&\ddots&\\&&X_d}\colon X_i:\cX_i\to\cX_i,\ i=1,,\ldots,d\right\},
\end{equation}
which is also equal to the commutant of $\{Z(z)\colon z\in\BZ^d\}$
in the $C^*$-algebra of bounded operators on $\cX$. We then say that
the realization $\{A,B,C,D\}$, or just $A$, is {\em scaled stable}
in case there exists an invertible operator $Q\in\cD$ so that
$\|Q^{-1}AQ\|<1$, or, equivalently, if there exists a positive
definite operator $X$ (notation $X>0$) in $\cD$ so that $AXA^*-X<0$.
To see that the two definitions coincide, take either
$X=QQ^*\in\cD$, or, when starting with $X>0$, factor $X$ as $X=QQ^*$
for some $Q\in\cD$. It is not hard to see that scaled stability
implies stability. Indeed, assume there exists an invertible
$Q\in\cD$ so that $\|Q^{-1}AQ\|<1$. Then $Z(z)Q^{-1}AQ=Q^{-1}Z(z)AQ$
is a strict contraction for each $z\in\ov{\BD}^d$, and thus
$Q^{-1}(I-Z(z)A)Q=I-Z(z)Q^{-1}AQ$ is invertible on $\ov{\BD}^d$. But
then $I-Z(z)A$ is invertible on $\ov{\BD}^d$ as well, and $A$ is
stable. The converse direction, even though asserted in
\cite{Piekarski,LF81}, turns out not to be true in general, as shown
in \cite{Anderson} via a concrete example. The output pair
$\{C_2,A\}$ is then said to be {\em scaled-detectable} if there
exists an output-injection operator $L \colon \cY \to \cX$ so that
$A + LC_2$ is scaled stable, and the input pair $\{A, B_2\}$ is
called {\em scaled-stabilizable} if there exists a state-feedback
operator $F \colon \cX \to \cU$ so that $A + B_2 F$ is scaled
stable.

While a classical result for the 1-D case states that operator, Hautus and scaled
detectability, as well as operator, Hautus and scaled stabilizability, are equivalent,
in the multidimensional setting considered here only one direction is clear.

\begin{proposition}  \label{P:stable/Hautus}
Let $\{A,B,C,D\}$ be a given realization as above with decomposition of $B$, $C$
and $D$ as in \eqref{decomps}.
 \begin{enumerate}
\item If the output pair $\{C_2, A\}$ is scaled-detectable, then $\{C_2,A\}$ is also
operator-detectable. If the output pair $\{C_2, A\}$ is
operator-detectable, then $\{C_2, A\}$ is also Hautus-detectable.

\item If the input pair $\{A,B_2\}$ is scaled-stabilizable, then $\{A,B_2\}$ is also
operator-stabilizable. If the input pair $\{A,B_2\}$ is
operator-stabilizable, then $\{A,B_2\}$ is also Hautus-stabilizable.
\end{enumerate}
\end{proposition}

\begin{proof}
Since scaled stability is a stronger notion than stability,
the first implications of both (1) and (2) are obvious. Suppose that
$L \colon \cY \to \cX$ is such that $A + LC_2$ is stable.
Then
$$ \begin{bmatrix} I & - Z(z) L \end{bmatrix} \begin{bmatrix} I
- Z(z)A \\ C_{2} \end{bmatrix} = I - Z(z) (A + LC_{2})
$$
is invertible for all $z \in \overline{\mathbb D}^{d}$ from which it
follows that $\{C_2, A\}$ is Hautus-detectable. The last assertion
concerning stabilizability follows in a similar way by making use of
the identity
$$
\begin{bmatrix} I - A Z(z) & B_{2} \end{bmatrix} \begin{bmatrix} I
\\ - F Z(z) \end{bmatrix} = I - (A + B_{2}F) Z(z).
$$
\end{proof}

The combination of operator-detectability together with
operator-stabiliz\-ability is strong enough for stabilizability
of the realization $\{A,B,C,D\}$ and we have the following weak
analogue of Theorem \ref{T:1.4}

\begin{theorem}\label{T:NDstate-stab}
Let $\{A,B,C,D\}$ be a given realization as above with decomposition
of $B$, $C$ and $D$ as in \eqref{decomps} (with $D_{22} = 0$). Assume that $\{C_{2},A\}$
is operator-detectable and $\{A,B_{2}\}$ is operator-stabilizable.
Then $\{A,B,C,D\}$ is stabilizable. Moreover, in this case one
stabilizing controller is $K \sim \{A_{K}, B_{K}, C_{K}, D_{K}\}$
where
\begin{equation}   \label{stabcon}
\begin{bmatrix} A_{K} & B_{K} \\ C_{K} & D_{K} \end{bmatrix} =
    \begin{bmatrix} A + B_{2} F + L C_{2} & -L
    \\ F & 0 \end{bmatrix}
\end{equation}
where $L:\cY\to\cX$ and $F:\cX\to\cU$ are any operators chosen
such that $A+LC_2$ and $A+FB_2$ are stable.
\end{theorem}

\begin{proof} It is possible to motivate these formulas with some
    observability theory (see \cite{DP}) but, once one has the
    formulas, it is a simple direct check that
    \begin{align*}
\begin{bmatrix} A_{cl} & B_{cl} \\ C_{cl} & D_{cl} \end{bmatrix} & =
    \begin{bmatrix} A + B_{2} D_{K} C_{2} & B_{2} C_{K} \\ B_{K}
    C_{2} & A_{K} \end{bmatrix}  \\
    & = \begin{bmatrix} A & B_{2} F \\ - L C_{2} & A + B_{2}F + L
    C_{2} \end{bmatrix}.
 \end{align*}
 It is now a straightforward exercise to check that this last matrix
 can be put in the triangular form $\sbm{ A + L C_{2} & 0 \\ -LC_{2}
 & A + B_{2}F}$ via a sequence of block-row/block-column similarity
 transformations, from which we conclude that $A_{cl}$ is stable as
 required.
\end{proof}

\begin{remark}  \label{R:NDcoprime}  {\em A result for the
      systems-over-rings setting that is analogous to that of Theorem
      \ref{T:NDstate-stab} is given in \cite{KharSon}.  There the
      result is given in terms of a Hautus-type
      stabilizable/detectable condition; in the systems-over-rings setting,
      Hautus-detectability/stabilizability is equivalent to
      operator-detectability/stabilizability (see Theorem 3.2 in \cite{KKT})
      rather than merely sufficient as
      in the present setting (see Proposition \ref{P:stable/Hautus}
      above).  The Hautus-type notions of detectability and stabilizability in
      principle are checkable using methods from \cite{Jury}: see the
      discussion in \cite[page 161]{KKT}.  The
      weakness of Theorem \ref{T:NDstate-stab} for our
      multidimensional setting is that there are no checkable
      criteria for when $\{C_2,A\}$ and $\{A,B_2\}$ are
      operator-detectable and operator-stabilizable since the Hautus test is
      in general only necessary.

      An additional weakness of Theorem \ref{T:NDstate-stab} is that
      it goes in only one direction: we do not assert that
      operator-detectability of $\{C_{2},A\}$ and operator-stabilizability for
      $\{A, B_{2}\}$ is necessary for stabilizability
      of $\{A,B,C,D\}$.  These weaknesses probably explain why
      apparently this result does not appear explicitly in the
      control literature.
      }\end{remark}

While there are no tractable necessary and sufficient
conditions for solving the state-space stabilizability
problem available, the situation turns out quite differently
when working with the more conservative notion of scaled stability.
The following is a more complete analogue of Theorem \ref{T:1.4}
combined with Theorem \ref{T:1.3}.

%%%%%%%%%%%%%%%%%%%%%%%%%%%%%%%%%%%%%%%%%%%%%%%%%%%%%%%%%%%%%%%%%%%%%%%%
\begin{theorem}\label{T:stab}
Let $\{A,B,C,D\}$ be a given realization. Then $\{A, B, C, D\}$ is
scaled-stabilizable, i.e., there exists a controller
$K$ with realization $\{A_K,B_K,C_K,D_K\}$ so that the closed loop state operator
$A_{cl}$ is scaled stable, if and only if the
input pair $\{A,B_2\}$ is scaled operator-stabilizable and the
output pair $\{C_2,A\}$ is scaled operator-detectable, i.e.,
there exist matrices $F$ and $L$ so that $A+B_2F$ and $A+LC_2$ are
scaled stable.  In this case the controller $K$ given by
\eqref{stabcon} solves the scaled-stabilization problem for $\{A, B,
C, D\}$.  Moreover:
\begin{enumerate}
    \item The following conditions concerning the input pair
    are equivalent:
    \begin{enumerate}
    \item $\{A, B_{2}\}$ is scaled operator-stabilizable.

    \item There exists $Y \in \cD$ satisfying the LMIs:
    \begin{equation} \label{stableLMI1Y}
    B_{2,\perp}(AYA^*-Y)B_{2,\perp}^*<0, \quad Y > 0
    \end{equation}
    where $B_{2, \perp}$ any injective operator with range equal
    to $ \operatorname{Ker} B_{2}$.

    \item There exists $Y \in \cD$ satisfying the LMIs
    \begin{equation}   \label{stableLMI2Y}
        AYA^{*} - Y - B_{2} B_{2}^{*} < 0, \quad Y > 0.
     \end{equation}
   \end{enumerate}
   \item The following conditions concerning the output pair are
   equivalent:
   \begin{enumerate}
       \item $\{C_{2} A,\}$ is scaled operator-detectable.

       \item There exists $X \in \cD$ satisfying the LMIs:
       \begin{equation} \label{stableLMI1X}
       C_{2,\perp}^{*} (A^{*} X A - X) C_{2, \perp} < 0, \quad X>0.
       \end{equation}
       where $C_{2, \perp}$ any injective operator with range equal
       to $\operatorname{Ker} C_{2}$.

       \item There exists $X \in \cD$ satisfying the LMIs
       \begin{equation}   \label{stableLMI2X}
          A^{*}X A - X - C_{2}^{*} C_{2} < 0, \quad X > 0.
        \end{equation}
      \end{enumerate}
    \end{enumerate}
   \end{theorem}

One of the results we shall use in the proof of Theorem \ref{T:stab} is known as
Finsler's lemma \cite{Finsler}, which also plays a key role in \cite{LZD,IS}.
This result can be interpreted as a refinement of the Douglas lemma \cite{D66}
which is well known in the operator theory community.

%%%%%%%%%%%%%%%%%%%%%%%%%%%%%%%%%%%%%%%%%%%%%%%%%%%%%%%%%%%%%%%%%%%%%%%%
\begin{lemma}[Finsler's lemma]\label{L:FinslerI}
Assume $R$ and $H$ are given matrices of appropriate size with
$H=H^*$. Then there exists a $\mu>0$ so that $\mu R^{*}R>H$ if and
only if $R_\perp^* H R_\perp<0$ where $R_\perp$ is any injective
operator with range equal to $\ker R$.
\end{lemma}

Finsler's lemma can be seen as a special case of another important result, which we
shall refer to as Finsler's lemma II. This is one of the main underlying tools
in the proof of the solution to the $H^\infty$-problem obtained in \cite{GA,AG}.

%%%%%%%%%%%%%%%%%%%%%%%%%%%%%%%%%%%%%%%%%%%%%%%%%%%%%%%%%%%%%%%%%%%%%%%%
\begin{lemma}[Finsler's lemma II]\label{L:FinslerII}
Given matrices $R$, $S$ and $H$ of appropriate sizes with $H=H^*$,
the following are equivalent:
\begin{itemize}
\item[(i)] There exists a matrix $J$ so that
$H+\mat{cc}{R^*&S^*}\mat{cc}{0&J^*\\J&0}\mat{c}{R\\S}<0$,

\item[(ii)] $R_\perp^*HR_\perp<0$ and  $S_\perp^*HS_\perp<0$, where $R_\perp$ and $S_\perp$
are injective operators with ranges equal to $\ker R$ and $\ker S$, respectively.
\end{itemize}
\end{lemma}

The proof of Finsler's Lemma II given in \cite{GA} uses only basic linear algebra and
is based on a careful administration of the kernels and ranges from the various matrices.
In particular, the matrices $J$ in statement (i) can actually be constructed from the data.
We show here how Finsler's lemma follows from the extended version.

%%%%%%%%%%%%%%%%%%%%%%%%%%%%%%%%%%%%
\begin{proof}[Proof of lemma \ref{L:FinslerI} using Lemma \ref{L:FinslerII}]
Apply Lemma \ref{L:FinslerII} with $R=S$. Then (ii) reduces to
$R_\perp^* H R_\perp<0$, which is equivalent to the existence of a matrix $J$ so
that $K=-(J^*+J)$ satisfies $R^{*}KR >H$. Since for such a matrix $K$ we have $K^*=K$,
it follows that $R^{*}\wtil{K}R >H$ holds for $\wtil{K}=\mu I$ as long as $\mu I>K$.
\end{proof}
%%%%%%%%%%%%%%%%%%%%%%%%%%%%%%%%%%%%

With these results in hand we can proof Theorem \ref{T:stab}.

%%%%%%%%%%%%%%%%%%%%%%%%%%%%%%%%%%%%
\begin{proof}[Proof of Theorem \ref{T:stab}]
We shall first prove that scaled stabilizability of $\{A,B,C,D\}$ is equivalent to
the existence of solutions $X$ and $Y$ in $\cD$ for the LMIs \eqref{stableLMI1X}
and \eqref{stableLMI1Y}. Note that $A_{cl}$ can be written in the following
affine way:
\begin{equation}   \label{Acl}
A_{cl}=\mat{cc}{A&0\\0&0}+\mat{cc}{0&B_2\\I&0}\mat{cc}{A_K&B_K\\C_K&D_K}\mat{cc}{0&I\\C_2&0}.
\end{equation}
Now let $X_{cl}:\cX\oplus\cX_K$ be an invertible matrix in $\cD_{cl}$, where $\cD_{cl}$
stands for the commutant of $\{Z_{cl}(z)\colon z\in\BZ^d\}$. Let $X$ be the compression
of $X_{cl}$ to $\cX$ and $Y$ the compression of $X_{cl}^{-1}$ to $\cX$. Then $X,Y\in\cD$.
Assume that $X_{cl}>0$. Thus, in particular, $X>0$ and $Y>0$. Then
$A_{cl}X_{cl}A_{cl}-X_{cl}<0$ if and only if
\begin{equation}\label{SC1}
\mat{cc}{-X_{cl}^{-1}&A_{cl}\\A_{cl}^*&-X_{cl}}<0.
\end{equation}
Now define
\[
H=\mat{cc}{-X_{cl}^{-1}&\mat{cc}{A&0\\0&0}\\\mat{cc}{A^*&0\\0&0}&-X_{cl}},\quad
R^*=\mat{cc}{0&0\\0&0\\0&C_2^*\\I&0},\quad S^*=\mat{cc}{0&I\\B_2&0\\0&0\\0&0}
\]
and
\[
J=\mat{cc}{A_K&B_K\\C_K&D_K}.
\]
Note that $H$, $R$ and $S$ are determined by the problem data, while $J$ amounts to the
system matrix of the controller to be designed. Then
\begin{equation}\label{}
\mat{cc}{-X_{cl}^{-1}&A_{cl}\\A_{cl}^*&-X_{cl}}
=H+\mat{cc}{R^*&S^*}\mat{cc}{0&J^*\\J&0}\mat{c}{R\\S}.
\end{equation}
Thus, by Finsler's lemma II, the inequality \eqref{SC1} holds for some $J=\sbm{A_K&B_K\\C_K&D_K}$
if and only if $R_\perp^*HR_\perp<0$ and $S_\perp^*HS_\perp<0$, where without loss
of generality we can take
\[
R_\perp=\mat{ccc}{I&0&0\\0&I&0\\0&0&C_{2,\perp}\\0&0&0}\ands
S_\perp=\mat{ccc}{0&0&0\\B_{2,\perp}&0&0\\0&I&0\\0&0&I}
\]
with $C_{2,\perp}$ and $B_{2,\perp}$ as described in part (b) of statements 1 and 2.
Writing out $R_\perp^*HR_\perp$ we find that $R_\perp^*HR_\perp<0$ if and only if
\[
\mat{cc}{-X_{cl}^{-1}&\mat{c}{AC_{2,\perp}\\0}\\
\mat{cc}{C_{2,\perp}^*A^*&0}&-C_{2,\perp}^*XC_{2,\perp}}<0
\]
which, after taking a Schur complement, turns out to be equivalent to
\[
C_{2,\perp}^*(A^*XA-X)C_{2,\perp}
=\mat{cc}{C_{2,\perp}^*A^*&0}X_{cl}\mat{c}{AC_{2,\perp}\\0}-C_{2,\perp}^*XC_{2,\perp}<0.
\]
A similar computation shows that $S_\perp^*HS_\perp<0$ is equivalent to
$B_{2,\perp}(AYA^*-Y)B_{2,\perp}^*<0$. This proves the first part of our claim.

For the converse direction assume we have $X$ and $Y$ in $\cD$ satisfying
\eqref{stableLMI1X}--\eqref{stableLMI1Y}.
Most of the implications in the above argumentation go both ways,
and it suffices to prove that there exists an operator $X_{cl}$ on
$\cX\oplus\cX_K$ in $\cD_{cl}$, with $\cX_K$ an arbitrary finite dimensional
Hilbert space with some partitioning
$\cX_K=\cX_{K,1}\oplus\cdots\oplus\cX_{K,d}$, so that $X_{cl}>0$
and $X$ and $Y$ are the compressions to $\cX$ of $X_{cl}$ and $X_{cl}^{-1}$,
respectively. Since \eqref{stableLMI1X}--\eqref{stableLMI1Y} hold with $X$ and $Y$ replaced by
$\rho X$ and $\rho Y$ for any positive number $\rho$, we may assume without loss of
generality that $\sbm{X&I\\I&Y}>0$. The existence of the required matrix $X_{cl}$ can
then be derived from Lemma 7.9 in \cite{DP} (with $n_K=n$). To enforce the fact
that $X_{cl}$ be in $\cD_{cl}$ we decompose $X=\textup{diag}(X_1,\ldots,X_d)$ and
$Y=\textup{diag}(Y_1,\ldots,Y_d)$ as in \eqref{Zcom} and complete $X_i$ and $Y_i$ to
positive definite matrices so that $\sbm{X_i&*\\ *&*}^{-1}=\sbm{Y_i&*\\ *&*}$.

To complete the proof it remains to show the equivalences of parts (a), (b) and (c)
in both statements 1 and 2. The equivalences of the parts (b) and (c) follows immediately
from Finsler's lemma with $R=B_2$ (respectively, $R=C_2^*$) and $H=AYA^*-Y$ (respectively,
$H=A^*XA-X$), again using that $X$ in \eqref{stableLMI1X} can be replaced with $\mu X$
(respectively, $Y$ in \eqref{stableLMI1Y} can be replaced with $\mu Y$) for any positive
number $\mu$.

We next show that (a) is equivalent to (b) for statement 1; for statement 2 the result
follows with similar arguments. Let $F:\cX\to\cU$, and let $X\in\cD$ be positive definite.
Taking a Schur complement it follows that
\begin{equation}\label{SC2}
(A^*+F^*B_2^*)X(A+B_2F)-X<0
\end{equation}
if and only if
\[
\mat{cc}{-X^{-1}&A+B_2F\\A^*+F^*B_2^*&-X}<0.
\]
Now write
\begin{align*}
&\mat{cc}{-X^{-1}&A+B_2F\\A^*+F^*B_2^*&-X}=\\
&\qquad\qquad\mat{cc}{-X^{-1}&A\\A^*&-X}+\mat{cc}{B_2&0\\0&I}\mat{cc}{0&F\\F^*&0}\mat{cc}{B_2^*&0\\0&I}.
\end{align*}
Thus, applying Finsler's lemma II with
\begin{equation}\label{FinsData}
H=\mat{cc}{-X^{-1}&A\\A^*&-X},\quad R=\mat{cc}{B_2^*&0},\quad
S=\mat{cc}{0&I} \ands J=F,
\end{equation}
we find that there exists an $F$ so that \eqref{SC2} holds if and only if
\[
R_\perp^* HR_\perp<0 \ands S_\perp^* HS_\perp<0
\]
with now $R_\perp=\sbm{B_{2,\perp}&0\\0&I}$ and $S_\perp=\sbm{I\\0}$.
The latter inequality is the same as $-X^{-1}<0$ and thus vacuous.
The first inequality, after writing out $R_\perp^* HR_\perp$, turns out to be
\[
\mat{cc}{-B_{2,\perp}^*X^{-1}B_{2,\perp}&B_{2,\perp}^*A\\A^*B_{2,\perp}& -X}<0,
\]
which, after another Schur complement, is equivalent to
$B_{2,\perp}^*(AX^{-1}A^*-X^{-1})B_{2,\perp}$.
\end{proof}
%%%%%%%%%%%%%%%%%%%%%%%%%%%%%%%%%%%%

Since scaled stability implies stability, it is clear that
finding operators $F$ and $L$ wit $A+B_2F$ and $A+LC_2$
scaled-stable implies that $A+B_2F$ and $A+LC_2$ are also
stable. In particular, having such operators $F$ and $L$ we
find the coprime factorization of $G_{22}$ via the functions in
Theorem \ref{T:NDstate-stab}. While there are no known tractable
necessary and sufficient conditions for
operator-detectability/stabilizability, the LMI criteria in parts
(iii) and (iv) of Theorem \ref{T:stab} for the scaled versions are
considered computationally tractable. Moreover, an
inspection of the last part of the proof shows how operators $F$ and
$L$ so that $A+B_2F$ and $A+LC_2$ are scaled stable can be
constructed from the solutions $X$ and $Y$ from the LMIs in
\eqref{stableLMI1Y}--\eqref{stableLMI2X}: Assume we have $X,Y\in\cD$ satisfying
\eqref{stableLMI1Y}--\eqref{stableLMI2X}. Define $H$, $R$ and $S$ as in \eqref{FinsData},
and determine a $J$ so that
$H+\sbm{R^*&S^*}\sbm{0&J^*\\J&0}\sbm{R\\S}<0$; this is possible as
the proof of Finsler's lemma II is essentially constructive. Then
take $F=J$. In a similar way one can construct $L$ using the LMI
solution $Y$.\bigskip

\paragraph{\bf Stability versus scaled stability, $\mu$ versus $\widehat\mu$}

We observed above that the notion of scaled stability is stronger,
and more conservative than the more intuitive notions of
stability in the Hautus or asymptotic sense. This
remains true in a more general setting that has proved useful in the
study of robust control \cite{LZD,DP,PackardDoyle} and that we will
encounter later in the paper.

Let $A$ be a bounded linear operator on a Hilbert space $\cX$. Assume that in
addition we are given a unital $C^*$-algebra $\DS$ which is realized concretely
as a subalgebra of $\cL(\cX)$, the space of bounded linear operators on $\cX$.
The complex structured singular value $\mu_{\DS}(A)$ of $A$ (with respect to the
structure $\DS$) is defined as
\begin{equation}\label{mu}
\mu_\DS(A)=\frac{1}{\inf\{\si(\De)\colon \De\in\DS,\ I-\De A\text{ is not invertible}\}}.
\end{equation}
Here $\si(M)$ stands for the largest singular value
of the operator $M$. Note that this contains two standard measures for $A$: the operator
norm $\|A\|$ if we take $\DS=\cL(\cX)$, and $\rho(A)$, the spectral radius of $A$,
if we take $\DS=\{\la I_\cX\colon \la\in\BC\}$; it is not hard to see that for
any unital $C^*$-algebra $\DS$ we have $\rho(A) \leq\mu_\DS(A)\leq
\|A\|$. See
\cite{PackardDoyle} for a tutorial introduction on the complex structured singular value
and \cite{FM00} for the generalization to algebras of operators on infinite dimensional spaces.

The $C^*$-algebra that comes up in the context of stability for the $N$-D systems
studied in this section is $\DS=\{Z(z)\colon z\in\BC^d\}$. Indeed, note that for this choice
of $\DS$ we have that $A$ is stable if and only if
$\mu_\DS(A)<1$.

In order to introduce the more conservative measure for $A$ in this context, we
write $\cD_\DS$ for the commutant of the $C^*$-algebra $\DS$ in $\cL(\cX)$.
We then define
\begin{align}\notag
\widehat\mu_\DS(A) & = \inf\{\gamma\colon \|Q^{-1}AQ\|<\ga\text{ for some invertible }Q\in\cD_\DS\}\\
&= \inf\{\gamma\colon  AXA^*-\gamma X<0 \text{ for some }X\in\cD_\DS,X>0\}.
\label{muhat}
\end{align}
The equivalence of the two definitions again goes through the relation between $X$
and $Q$ via $X=Q^*Q$. It is immediate that with $\DS=\{Z(z)\colon z\in\BC^d\}$
we find $\cD_\DS=\cD$ as in \eqref{Zcom}, and that $A$ is scaled stable if and
only if $\widehat\mu_\DS(A)<1$.\bigskip

\paragraph{\bf The state-space $H^\infty$-problem}
The problems of finding tractable necessary and sufficient conditions for the strict
state-space $H^\infty$-problem are similar to that for the state-space
stabilizability problem. Here one also typically resorts to a more conservative
`scaled' version of the problem.

We say that the realization $\{A,B,C,D\}$ with decomposition \eqref{decomps}
has {\em scaled performance} whenever there exists an invertible $Q\in\cD$ so
that
\begin{equation}\label{ScaledPerf1}
\left\|\mat{cc}{Q^{-1}&0\\0&I_{\cZ\oplus\cY}}\mat{cc}{A&B\\C&D}
\mat{cc}{Q&0\\0&I_{\cW\oplus\cU}}\right\|<1,
\end{equation}
or, equivalently, if there exists an $X>0$ in $\cD$ so that
\begin{equation}\label{ScaledPerf2}
\mat{cc}{A&B\\C&D}\mat{cc}{X&0\\0&I_{\cW\oplus\cU}}\mat{cc}{A&B\\C&D}^*
-\mat{cc}{X&0\\0&I_{\cW\oplus\cU}}<0.
\end{equation}
The equivalence of the two definitions goes as for the scaled stability case
through the relation $X=QQ^*$. Looking at the left upper entry in \eqref{ScaledPerf2}
it follows that scaled performance of $\{A,B,C,D\}$ implies scaled stability.
Moreover, if \eqref{ScaledPerf1} holds for $Q\in\cD$, then it is not hard to see
that the transfer function $G(z)$ in \eqref{transfer} is also given by
\[
G(z)=D+C'(I-Z(z)A')^{-1}Z(z)B'
\]
where the system matrix
\[
\mat{cc}{A'&B'\\C'&D}=
\mat{cc}{Q^{-1}&0\\0&I_{\cZ\oplus\cY}}\mat{cc}{A&B\\C&D}\mat{cc}{Q&0\\0&I_{\cW\oplus\cU}}
\]
is equal to a strict contraction. It then follows from a standard fact on feedback connections
(see e.g.~Corollary 1.3 page 434 of \cite{FF} for a very general
formulation)
that $\|G(z)\|<1$ for $z\in\ov{\BD}^d$, i.e., $G$ has strict performance.
The {\em scaled $H^\infty$-problem} is then to find
a controller $K$ with realization $\{A_K,B_K,C_K,D_K\}$ so that
the closed loop system $\{A_{cl},B_{cl},C_{cl},D_{cl}\}$ has scaled performance. The above
analysis shows that solving the scaled $H^\infty$-problem implies solving that state-space
$H^\infty$-problem. The converse is again not true in general. Further elaboration of the
same techniques as used in the proof of Theorem \ref{T:stab} yields the following result
for the scaled $H^\infty$-problem; see \cite{AG,GA}. For the connections between the Theorems
\ref{T:perform} and \ref{T:stab}, in the more general setting of LFT models with structured
uncertainty, we refer to \cite{BFGtH09}.  Note that the result
collapses to Theorem \ref{T:1.5} given in the Introduction when we
specialize to the single-variable case $d=1$.

%%%%%%%%%%%%%%%%%%%%%%%%%%%%%%%%%%%%%%%%%%%%%%%%%%%%%%%%%%%%%%%%%%%%%%%%
\begin{theorem}\label{T:perform}
Let $\{A,B,C,D\}$ be a given realization. Then there exists a solution for
the scaled $H^\infty$-problem associated with $\{A,B,C,D\}$ if and only if
there exist $X,Y\in\cD$ satisfying LMIs:
\begin{align}
       & \begin{bmatrix} N_{c} & 0 \\ 0 & I \end{bmatrix} ^{*}
            \begin{bmatrix} AYA^{*} - Y & AYC_{1}^{*} & B_{1} \\
     C_{1}Y A^{*} & C_{1} Y C_{1}^{*} - I & D_{11} \\
       B_{1}^{*} & D_{11}^{*} & -I \end{bmatrix} \label{Y-LMI}
     \begin{bmatrix} N_{c} & 0 \\ 0 & I \end{bmatrix} <  0,\quad Y>0, \\
         & \begin{bmatrix} N_{o} & 0 \\ 0 & I \end{bmatrix} ^{*}
        \begin{bmatrix} A^{*} X A - X & A^{*} X B_{1} & C_{1} ^{*} \\
       B_{1}^{*} X A & B_{1}^{*} X B_{1} - I & D_{11}^{*} \\
       C_{1} & D_{11} & -I \end{bmatrix}
       \begin{bmatrix} N_{o} & 0 \\ 0 & I \end{bmatrix} < 0,\quad X>0,
           \label{X-LMI}
      \end{align}
and the coupling condition
\begin{equation}\label{coupling}
\mat{cc}{X&I\\I&Y}\geq 0.
\end{equation}
Here $N_{c}$ and $N_{o}$ are matrices chosen so that
    \begin{align*}
        & N_{c} \text{ is injective and }
     \operatorname{Im } N_{c} = \operatorname{Ker } \begin{bmatrix}
     B_{2}^{*} & D_{12}^{*} \end{bmatrix} \text{ and } \\
     & N_{o} \text{ is injective and } \operatorname{Im } N_{o} =
     \operatorname{Ker } \begin{bmatrix} C_{2} & D_{21} \end{bmatrix}.
     \end{align*}
\end{theorem}

Note that Theorem \ref{T:perform} does not require that the problem
be first brought into model-matching form; thus this solution
bypasses the Nevanlinna-Pick-interpolation interpretation of the
$H^\infty$-problem.

\subsection{Equivalence of frequency-domain and state-space
  formulations}   \label{S:ND-equiv}

  In this subsection we suppose that we are given a transfer matrix
   $G$ of size $(n_{\cZ}+n_{\cY}) \times (n_{\cW}+ n_{\cU})$ with
   coefficients in $Q({\mathbb C}(z)_{ss})$ as in Section
   \ref{S:ND-freq} with a given state-space realization as in Subsection
\ref{S:ND-statespace}:
   \begin{equation}  \label{G-real}
     G(z) = \begin{bmatrix} G_{11} & G_{12} \\ G_{21} & G_{22}
 \end{bmatrix} =
     \begin{bmatrix} D_{11} & D_{12} \\ D_{21} & D_{22} \end{bmatrix} +
     \begin{bmatrix} C_{1} \\C_{2} \end{bmatrix}  (I - Z(z)A)^{-1} Z(z)
         \begin{bmatrix} B_{1} & B_{2} \end{bmatrix}
   \end{equation}
   where $Z(z)$ is as in \eqref{Z(z)}. We again consider the problem of
finding stabilizing controllers $K$, also equipped with a state-space
realization
   \begin{equation}   \label{K-real}
     K(z) = D_{K} + C_{K} (I - Z_K(z) A_{K})^{-1} Z_K(z) B_{K},
   \end{equation}
in either the state-space stability or in the frequency-domain stability sense.
A natural question is whether the frequency-domain $H^{\infty}$-problem
with formulation in state-space coordinates is the same as the state-space
$H^{\infty}$-problem formulated in Section \ref{S:ND-statespace}.

   For simplicity in the computations to follow, we shall always
   assume that the plant $G$ has been normalized so that $D_{22} = 0$.
   In one direction the result is clear.  Suppose that $K(z) = D_{K} +
   C_{K}(I - Z(z) A_{K})^{-1} Z(z) B_{K}$ is a stabilizing controller
   for $G(z)$ in the state-space sense.  It follows that the
   closed-loop state matrix
   \begin{equation}   \label{Acl'}
     A_{cl} =\begin{bmatrix}  A + B_{2} D_{K} C_{2} & B_{2}C_{K} \\ B_{K} C_{2} &
     A_{K} \end{bmatrix}
   \end{equation}
   is stable, i.e.,  $I - Z_{cl}(z)A_{cl}$ is invertible for all $z$ in
the closed polydisk $\overline{\mathbb D}^{d}$, with $Z_{cl}(z)$ as defined in
Subsection \ref{S:ND-statespace}. On the other hand
   one can compute that the transfer matrix $\Theta(G_{22},K): = \sbm{
   I & -K(z) \\ -G_{22}(z) & I}^{-1}$ has realization
   \begin{equation}  \label{Wtilde}
       \widetilde W(z) = \begin{bmatrix} I & I_{D_{K}} \\ 0 & I
   \end{bmatrix} +
   \begin{bmatrix} D_{K} C_{2} & C_{K} \\ C_{2} & 0 \end{bmatrix}
       \left( I - Z_{cl}(z) A_{cl} \right)^{-1} Z_{cl}(z)
       \begin{bmatrix} B_{2} & B_{2}D_{K} \\ 0 & B_{K} \end{bmatrix}.
   \end{equation}
   As the resolvent expression $\left(I - Z_{cl}(z)A_{cl}\right)^{-1}$
   has no singularities in the closed polydisk $\overline{\mathbb
   D}^{d}$, it is clear that $\widetilde W(z)$ has matrix entries in
   ${\mathbb C}(z)_{ss}$, and it follows that $K$ stabilizes
   $G_{22}$ in the frequency-domain sense.  Under the assumption that
   $G$ is internally stabilizable (frequency-domain sense), it
   follows from Corollary \ref{C:G22G} that $K$
   also stabilizes $G$ (frequency-domain sense).

   We show that the converse direction holds under an additional
   assumption.  The early paper \cite{2D-II} of
   Kung-L\'evy-Morf-Kailath introduced the notion of modal
   controllability and modal observability for 2-D systems.  We
   extend these notions to $N$-D systems as follows.  Given a
   Givone-Roesser output pair $\{C,A\}$, we say that $\{C,A\}$ is
{\em modally observable} if the block-column matrix $\sbm{ I -
   Z(z)A \\ C }$ has maximal rank $n_{\cX}$ for a generic point $z$
   on each irreducible component of the variety $\det (I - Z(z) A) =
   0$.  Similarly we say that the Givone-Roesser input pair $\{A, B\}$
   is {\em modally controllable} if the block-row matrix
   $\sbm{  I - A Z(z) & B }$ has maximal rank
   $n_{\cX}$ for a generic point on each irreducible component of the
   variety $\det (I - A Z(z)) = \det (I - Z(z) A) = 0$.  Then the
   authors of \cite{2D-II} define the realization $\{A,B,C,D\}$ to be
   {\em minimal} if both $\{C,A\}$ is modally observable and $\{A, B\}$
   is modally controllable.  While this is a natural notion of
   minimality, unfortunately it is not clear that an arbitrary
   realization $\{A, B, C, D\}$ of a given transfer function $S(z) = D
   + C (I - Z(z)A)^{-1} Z(z) B$ can be reduced to a minimal realization
   $\{A_{0}, B_{0}, C_{0}, D_{0}\}$ of the same transfer function
   $S(z) = D_{0} + C_{0} (I - Z(z) A_{0})^{-1} Z(z) B_{0}$.

   As a natural modification of the notions of modally observable and
   modally controllable, we now introduce the notions of modally
   detectable and modally stabilizable as follows.
   For $\{C,A\}$ a Givone-Roesser output pair,
    we say that $\{C, A\}$ is {\em modally detectable} if the
column matrix $\sbm{ I - Z(z)  A \\ C }$
    has maximal rank $n_{\cX}$ for a generic point $z$ on each
    irreducible component of the variety $\det (I - Z(z) A) = 0$ which
    enters into the polydisk $\overline{\mathbb D}^{d}$.  Similarly, we
    say that the Givone-Roesser input pair $\{A, B\}$ is {\em modally
    stabilizable} if the row matrix $\sbm{ I -  A Z(z) & B }$ has
    maximal rank $n_{\cX}$ for a generic point $z$ on each irreducible
    component of the variety $\det (I - Z(z) A)=0$ which has nonzero
    intersection with the closed polydisk $\overline{\mathbb D}^{d}$.
    We then have the following partial converse of the observation
    made above that state-space internal stabilization implies
    frequency-domain internal stabilization; this is an $N$-D version
    of Theorem \ref{T:1.6} in the Introduction.

      \begin{theorem} \label{T:NDfreq-state}
Let \eqref{G-real} and \eqref{K-real} be given realizations for
$G \colon \sbm{ \cW \\ \cU }\to \sbm{ \cZ \\ \cY }$ and $K \colon \cY \to \cU$.
Assume that $\{C_2,A\}$ and $\{C_K,A_K\}$ are modally detectable and
$\{A,B_2\}$ and $\{A_K,B_K\}$ are modally stabilizable. Then $K$
internally stabilizes $G_{22}$ in the state-space sense (and thus state-space
stabilizes $G$) if and only if $K$ stabilizes $G_{22}$ in
the frequency-domain sense (and $G$ if $G$ is stabilizable in the frequency-domain sense).
       \end{theorem}

     \begin{remark}\label{R:NDfreq-state}
     {\em  As it is not clear that a given realization can
     be reduced to a modally observable and modally controllable
     realization for a given transfer function, it is equally not
     clear whether a given transfer function has a modally
     detectable and modally stabilizable realization.  However,
     in the case that $d=1$, such realizations always exists and
     Theorem \ref{T:NDfreq-state} recovers the standard 1-D
     result (Theorem \ref{T:1.6} in the Introduction).
     }\end{remark}

    The proof of Theorem \ref{T:NDfreq-state} will make frequent
    use of the following basic result from the theory of holomorphic
    functions in several complex variables. For the proof we refer
    to \cite[Theorem 4 page 176]{Shabat}; note that if
    the number of variables $d$ is 1, then the only analytic set of
    codimension at least  2 is the empty set and the theorem is
    vacuous; the theorem has content only when the number of
    variables is at least 2.

   \begin{theorem} \label{T:Shabat} {\em \textbf{Principle of Removal
       of Singularities}}  Suppose that the complex-valued function
       $\varphi$ is holomorphic on a set $S$ contained in ${\mathbb
       C}^{d}$ of the form $S = {\mathcal D} - \cE$ where $\cD$ is an
       open set in ${\mathbb C}^{d}$ and $\cE$ is the intersection
       with $\cD$ of an analytic set of codimension at least 2.  Then
       $\varphi$ has analytic continuation to a function holomorphic
       on all of $\cD$.
     \end{theorem}

   We shall also need some preliminary lemmas.

   \begin{lemma}  \label{L:feedback-inv}
       \begin{enumerate}
       \item Modal detectability is invariant under output
       injection, i.e., given a Givone-Roesser output pair $\{C,A\}$
       (where $A \colon \cX \to \cX$ and $C \colon \cX \to \cY$)
       together with
       an output injection operator $L \colon \cY \to \cX$, then
       the pair $\{C,A\}$ is modally detectable if and only if the
       pair $\{C, A + LC\}$ is modally detectable.

       \item Modal stabilizability is invariant under state
       feedback, i.e.,  given a Givone-Roesser input pair $\{A,B\}$
       (where $A \colon \cX \to \cX$ and $B \colon \cU \to \cX$)
       together with a state-feedback operator $F \colon \cX \to
       \cU$,  then the pair $\{A, B\}$ is modally stabilizable if
       and only if the pair $\{A + B F, B\}$ is modally stabilizable.
       \end{enumerate}
 \end{lemma}

 \begin{proof}
     To prove the first statement, note the identity
     $$
     \begin{bmatrix} I & -Z(z) L \\ 0 & I \end{bmatrix}
     \begin{bmatrix} I - Z(z) A \\ C \end{bmatrix} =
         \begin{bmatrix} I - Z(z) (A + LC) \\ C \end{bmatrix}.
    $$
    Since the factor $\sbm{ I & - Z(z) L \\ 0 & I }$ is invertible for
    all $z$, we conclude that, for each $z \in {\mathbb C}^{d}$,
    $\sbm{I - Z(z) A \\ C}$ has maximal rank exactly when $\sbm{I -
    Z(z) (A + LC) \\ C}$ has maximal rank, and hence, in particular,
    the modal detectability for  $\{C,A\}$ holds exactly when modal
    detectability for $\{C, A + LC\}$ holds.

    The second statement follows in a similar way from the identity
    $$
    \begin{bmatrix} I - A Z(z) & B \end{bmatrix} \begin{bmatrix} I & 0
    \\ -F Z(z) & I \end{bmatrix} = \begin{bmatrix} I - (A + BF)
    Z(z) & B \end{bmatrix}.
    $$
    \end{proof}

    \begin{lemma}   \label{L:ND-sing}
    Suppose that the function $W(z)$ is stable (i.e., all
    matrix entries of $W$ are in ${\mathbb C}(z)_{ss}$) and suppose
    that
     \begin{equation}  \label{real-W}
     W(z) = D + C (I - Z(z) A)^{-1} Z(z) B
    \end{equation}
    is a realization for $W$ which is both modally detectable and
    modally  stabilizable.  Then the matrix $A$ is stable, i.e.,
    $(I - Z(z) A)^{-1}$ exists for all $z$ in the closed polydisk
    $\overline{\mathbb D}^{d}$.
    \end{lemma}

    \begin{proof}
    As $W$ is stable and $Z(z)B$ is trivially stable, then
    certainly
    \begin{equation}  \label{identity1}
     \begin{bmatrix} I - Z(z) A \\ C \end{bmatrix} (I - Z(z) A)^{-1}
     Z(z) B = \begin{bmatrix} Z(z) B \\ W(z) - D \end{bmatrix}
    \end{equation}
    is stable (i.e., holomorphic on $\overline{\mathbb D}^{d}$).
    Trivially $\sbm{ I - Z(z) A \\ C }$ has maximal rank $n_{\cX}$
    for all $z \in \overline{\mathbb D}^{d}$ where $\det (I - Z(z) A)
    \ne 0$.  By assumption, $\sbm{ I - Z(z)A \\ C }$ has maximal rank
    generically on each irreducible component of the zero variety of
    $\det (I - Z(z) A)$  which intersects $\overline{\mathbb D}^{d}$.
    We conclude that $\sbm{ I - Z(z) A \\ C }$ has maximal rank
    $n_{\cX}$ at all points of $\overline{\mathbb D}^{d}$ except those
    in an exceptional set $\cE$ which is contained in a subvariety,
    each irreducible component of which has codimension at least 2.
    In a neighborhood of each such point $z \in \overline{\mathbb D}^{d} - \cE$, $\sbm{
    I - Z(z) A \\ C }$ has a holomorphic left inverse; combining this
    fact with the identity \eqref{identity1}, we see that $(I - Z(z)
    A)^{-1} Z(z) B$ is holomorphic on $\overline{\mathbb D}^{d} -
    \cE$.  By Theorem \ref{T:Shabat}, it follows that $(I - Z(z)
    A)^{-1} Z(z) B$ has analytic continuation to all of $\overline{\mathbb
    D}^{d}$.

    We next note the identity
    \begin{equation}   \label{identity2}
    \begin{bmatrix} Z(z) & (I - Z(z) A)^{-1} Z(z) B  \end{bmatrix}
        = Z(z) (I - A Z(z))^{-1} \begin{bmatrix}  I - A Z(z) & B
        \end{bmatrix}
   \end{equation}
   where the quantity on the left-hand side is holomorphic on
   $\overline{\mathbb D}^{d}$ by the result established above.  By
   assumption $\{A, B\}$ is modally stabilizable; by an argument
   analogous to that used above for the modally detectable pair
   $\{C,A\}$, we see that the pencil $\begin{bmatrix} I - A Z(z) & B
 \end{bmatrix}$ has a holomorphic right inverse in the neighborhood of
 each point $z$ in $\overline{\mathbb D}^{d} - \cE'$ where the
 exception set $\cE'$ is contained in a subvariety each irreducible
 component of which has codimension at least 2.   Multiplication of the
 identity \eqref{identity2} on the right by this right inverse then
 tells us that $Z(z) (I - Z(z) A)^{-1}$ is holomorphic on
 $\overline{\mathbb D}^{d} - \cE'$.  Again by Theorem \ref{T:Shabat},
 we conclude that in fact $Z(z) (I - Z(z) A)^{-1}$ is holomorphic on all of
 $\overline{\mathbb D}^{d}$.

 We show that $(I - Z(z) A)^{-1}$ is holomorphic on
 $\overline{\mathbb D}^{d}$ as follows. Let $E_j:\cX\to\cX_j$ be the projection on the
$j$-th component of $\cX=\cX_1\oplus\cdots\oplus\cX_d$.
Note that the first block
 row of $(I - Z(z) A)^{-1}$ is equal to $z_{1} E_{1} (I - Z(z)
 A)^{-1}$.  This is holomorphic on the closed polydisk $\overline
 {\mathbb D}^{d}$.  For $z$ in a sufficiently small polydisk $|z_{i}|
 < \rho$ for $i = 1, \dots, d$, $(I - Z(z) A)^{-1}$ is analytic and
 hence $z_{1} E_{1} (I - Z(z) A)^{-1}|_{z_{1} = 0} = 0$.  By analytic
 continuation, it then must hold that $z_{1}(E_{1}(I - Z(z) A)^{-1} =
 0$ for all $z = (0, z_{2}, \dots, z_{d})$ with $|z_{i}| \le 1$ for $i
 = 2, \dots, d$.  For each fixed $(z_{2}, \dots, z_{d})$, we may
 use the single-variable result that one can divide out zeros to
 conclude that $E_{1} (I - Z(z) A)^{-1}$ is holomorphic in $z_{1}$ at $z_{1}
 = 0$. As the result is obvious for $z_{1} \ne 0$, we conclude that
 $E_{1} (I - Z(z) A)^{-1}$ is holomorphic on the whole closed polydisk
 $\overline{\mathbb D}^{d}$.  In a similar way working with the
 variable $z_{i}$,  one can show that $E_{i} (I - Z(z) A)^{-1}$ is
 holomorphic on the whole closed polydisk, and it follows that $(I -
 Z(z) A)^{-1} = \sbm{ E_{1} \\ \vdots \\ E_{d}} (I - Z(z) A)^{-1}$ is
 holomorphic on the whole closed polydisk as wanted.
 \end{proof}

    We are now ready for the proof of Theorem \ref{T:NDfreq-state}.

    \begin{proof}[Proof of Theorem \ref{T:NDfreq-state}.]
    Suppose that $K$ stabilizes $G_{22}$ in
    the frequency-domain sense.
    This simply means that the transfer
   function $\widetilde W$ given by \eqref{Wtilde} is holomorphic on the closed polydisk
   $\overline{\mathbb D}^{d}$. To show that $A_{cl}$ is stable, by
   Lemma \ref{L:ND-sing} it
   suffices to show that
   $\left\{ \sbm{ D_{K}C_{2} & C_{K} \\ C_{2} & 0}, A_{cl} \right\}$ is
   modally detectable and that $
 \left\{ A_{cl}, \sbm{  B_{2} & B_{2}D_{K} \\ 0 & B_{K} } \right\}$
 is modally stabilizable.

 To prove that $\left\{ \sbm{ D_{K}C_{2} & C_{K} \\ C_{2} & 0}, A_{cl} \right\}$
 is modally detectable,  from the definition \eqref{Acl} we note that
 $$
 A_{cl} = \begin{bmatrix} A & 0 \\ 0 & A_{K} \end{bmatrix} +
 \begin{bmatrix}  B_{2} & 0 \\ 0 & B_{K} \end{bmatrix}
 \begin{bmatrix} D_{K}C_{2} & C_{K} \\ C_{2} & 0 \end{bmatrix}.
 $$
 By Lemma \ref{L:feedback-inv} we see that modal detectability of
 $ \left\{\sbm{ D_{K} C_{2} & C_{K} \\ C_{2} & 0}, A_{cl} \right\}$
  is equivalent to modal detectability of $\left\{ \sbm{ D_{K}C_{2} & C_{K}
  \\ C_{2} & 0 }, \sbm{A & 0 \\ 0 & A_{K}}\right\}$.  As
  $\sbm{D_{K}C_{2} & C_{K} \\ C_{2} & 0} = \sbm{D_{K} & I \\ I & 0}
  \sbm{C_{2}& 0 \\ 0 & C_{K}}$ with $\sbm{D_{K} & I \\ I & 0}$
  invertible, it is easily seen that modal detectability of the input pair
  $\left\{ \sbm{D_{K}C_{2} & C_{K}
    \\ C_{2} & 0 }, \sbm{A & 0 \\ 0 & A_{K}}\right\}$
    is equivalent to modal detectability of
   $\left\{ \sbm{C_{2} & 0 \\ 0 & C_{K}}, \sbm{A & 0 \\ 0 & A_{K} }
   \right\}$.  But the modal detectability of this last pair in turn
   follows from its diagonal form and the assumed modal detectability
   of $\{C_{2}, A\}$ and $\{C_{K}, A_{K}\}$.

   The modal      stabilizability of $\left\{A_{cl}, \sbm{ B_{2} &
   B_{2}D_{K} \\ 0 & B_{K}} \right\}$ follows in a similar way by
   making use of the identities
   $$
   A_{cl}\! =\! \begin{bmatrix} A & 0 \\ 0 & A_{K} \end{bmatrix}\! +\!
   \begin{bmatrix} B_{2}D_{K} & B_{2} \\ B_{K} & 0 \end{bmatrix}\!
       \begin{bmatrix} C_{2} & 0 \\ 0 & C_{K} \end{bmatrix}, \
       \begin{bmatrix} B_{2}D_{K} & B_{2} \\ B_{K} & 0
       \end{bmatrix}\! =\! \begin{bmatrix} B_{2} & 0 \\ 0 & B_{K}
       \end{bmatrix}\! \begin{bmatrix} D_{K} & I \\ I & 0 \end{bmatrix}
   $$
   and noting that $\sbm{ D_{K} & I \\ I & 0 }$ is invertible.
  \end{proof}

  In both the frequency-domain setting of Section \ref{S:ND-freq} and
  the state-space setting of Section \ref{S:ND-statespace}, the true
  $H^{\infty}$-problem is intractable and we resorted to some
  compromise:  the Schur-Agler-class reformulation in Section
  \ref{S:ND-freq} and the scaled-$H^{\infty}$-problem reformulation
  in Section \ref{S:ND-statespace}.  We would now like to compare
  these compromises for the setting where they both apply, namely,
  where we are given both the transfer function $G$ and the
  state-space representation $\{A,B,C,D\}$ for the plant.

  \begin{theorem} \label{T:compare} Suppose that
      $ G(z) = \sbm{G_{11} & G_{12} \\ G_{21} & 0 }$ is in model-matching form with
      state-space realization $G(z) = D + C (I - Z(z) A)^{-1} Z(z) B$
      as in \eqref{G-real}.  Suppose that the controller $K(z) = D_{K} + C_{K}
      (I - Z_{K}(z) A_{K})^{-1} Z_{K}(z) B_{K}$ solves the scaled
      $H^\infty$-problem.  Then the transfer function
      $\widetilde W(z)$ as in \eqref{Wtilde} is a Schur-Agler-class
      solution of the Model-Matching problem.
      \end{theorem}

      \begin{proof}  Simply note that, under the assumptions of the
      theorem, $\widetilde W(z)$ has a realization
      $ \widetilde W = D_{cl} + C_{cl} (I - Z_{cl}(z)
      A_{cl})^{-1} Z_{cl}(z) B_{cl}$ for which there is a
      state-space change of coordinates $Q \in {\mathcal D}$ transforming the
      realization to a contraction:
      $$\left\| \begin{bmatrix} A' & B' \\ C' & D \end{bmatrix}  \right\|
      < 1 \text{ where } \begin{bmatrix} A' & B' \\ C' & D
      \end{bmatrix} = \begin{bmatrix} Q & 0 \\ 0 & I \end{bmatrix}
      \begin{bmatrix} A_{cl} & B_{cl} \\ C_{cl} & D_{cl} \end{bmatrix} \begin{bmatrix}
      Q^{-1} & 0 \\ 0 & I \end{bmatrix}.
    $$
    Thus we also have $\widetilde W(z) = D + C' (I - Z_{cl}(z)
    A') Z_{cl}(z) B'$ from which it follows that $W$ is in the
    strict Schur-Agler class, i.e., $\|\wtilW(X)\|<1$ for any $d$-tuple
    $X=(X_1,\ldots,X_d)$ of contraction operators $X_j$ on a separable
    Hilbert space $\cX$.
    By construction $\widetilde W$
    necessarily has the model matching form $\widetilde W =
    G_{11} + G_{12} \Lambda G_{21}$ with $\Lambda$ stable.
     \end{proof}

    \begin{remark} \label{R:Schur-Agler/scaledHinf} {\em
        In general a
        Schur-Agler function $S(z)$ can be realized with
        a colligation matrix $\sbm{A & B \\ C & D}$ which is not
        of the form
        \begin{equation}\label{statespacefail}
\begin{bmatrix} A & B \\ C & D \end{bmatrix} =
        \begin{bmatrix} Q^{-1} & 0 \\ 0 & I \end{bmatrix}
        \begin{bmatrix} A' & B' \\ C' & D \end{bmatrix}
            \begin{bmatrix} Q & 0 \\ 0 & I\end{bmatrix}
\end{equation}
    with $\sbm{A' & B' \\ C' & D }$ equal to a strict
    contraction and $Q\in\cD$ invertible. As an example, let $A$ be the
    block $2 \times 2$ matrix given by Anderson-et-al in \cite{Anderson}.
    This matrix has the property that $I - Z(z) A$ is invertible for all
    $z\in\ov{\BD}^2$, but there is no $Q\in\cD$ so that $\|Q^{-1}AQ\|<1$.
    Here $Z(z)$ and $\cD$ are compatible with the block decomposition of
    $A$. Then for $\gamma>0$ sufficiently small the function
    $S(z)=\gamma(I-Z(z)A)^{-1}$ has $\|S(z)\|\le \rho<1$ for some $0<\rho<1$
    and all $z\in\ov{\BD}^2$. Hence $S$ is a strict Schur-class
    function. As mentioned in Section \ref{S:ND-freq}, a consequence
    of the And\^o dilation theorem \cite{Ando} is that the Schur
    class and the Schur-Agler class coincide for $d=2$; it is not
    hard to see that this equality carries over to the strict
    versions and hence $S$ is in the strict Schur-Agler class.
    As a consequence of the strict Bounded-Real-Lemma
    in \cite{BGM3}, $S$ admits a strictly contractive state-space
    realization $\sbm{A'&B'\\C'&D}$.  However, the realization
    $\sbm{A&B\\C&D}=\sbm{A&A\\ \gamma I&\gamma I}$ of $S$, obtained from the fact
    that
    \[
    S(z)=\gamma(I-Z(z)A)^{-1}=\gamma I+\gamma(I-Z(z)A)^{-1}Z(z)A,
    \]
    cannot relate to $\sbm{A'&B'\\C'&D}$ as in \eqref{statespacefail} since
    that would imply the existence of an invertible $Q\in\cD$ so that
    $Q^{-1}AQ=A'$ is a strict contraction.
    } \end{remark}

    \begin{remark} {\em
    Let us assume that the $G(z)$ in Theorem \ref{T:compare} is
    such that $G_{12}$ and $G_{21}$ are square and invertible on
    the distinguished boundary ${\mathbb T}^{d}$ of the polydisk
    ${\mathbb D}^{d}$ so that the Model-Matching problem can be
    converted to a polydisk bitangential Nevanlinna-Pick
    interpolation problem along a subvariety as in \cite{BM}.  As
    we have seen, the solution criterion using the Agler
    interpolation theorem of \cite{Agler-unpublished, BT} then
    involves an LOI (Linear Operator Inequality or infinite LMI).
    On the other hand, if we assume that we are given a stable
    state-space realization $\{A, B, C, D\}$ for $G(z) = \sbm{
    G_{11}(z) & G_{12}(z) \\ G_{21}(z) & 0}$, we may instead
    solve the associated scaled $H^{\infty}$-problem associated
    with this realization data-set.  The associated solution
    criterion in Theorem \ref{T:perform} remarkably involves only
    finite LMIs.  A disadvantage of this state-space approach,
    however, is that in principle one would have to sweep all
    possible (similarity equivalence classes of) realizations of
    $G(z)$; while each non-equivalent realization gives a
    distinct $H^{\infty}$-problem, the associated
    frequency-domain Model-Matching/bitangential
    variety-interpolation problem remains the same.
    } \end{remark}

  \subsection{Notes}  \label{S:NDnotes}
  In \cite{Lin1} Lin conjectured the result stated in Theorem
  \ref{T:Lin} that $G_{22}$-stabilizability is equivalent to the
  existence of a stable coprime factorization for $G_{22}$.  This
  conjecture was settled by Quadrat (see \cite{Q-YKII, Q-LN,
  Q-elementary}) who obtained the equivalence of this property with
  projective-freeness of the underlying ring and noticed the
  applicability of the results from \cite{BST, KKT} concerning the
  projective-freeness of ${\mathbb C}(z)_{ss}$.

  For the general theory of the $N$-D systems, in particular for
  $N$=2, considered in Subsection \ref{S:ND-statespace} we refer
  to \cite{Kaczorek,DuXie}.

  The sufficiency of scaled stability
  for asymptotic/\-Hautus-stability goes back to \cite{ElAgiziFahmy}.
  Theorem \ref{T:stab} was proved in \cite{LZD} for the more general
  LFT models in the context of robust control with structured
  uncertainty. The proof given here is based on the extended
  Finsler's lemma (Lemma \ref{L:FinslerII}), and basically follows
  the proof from \cite{GA} for the solution to the scaled
  $H^\infty$-problem (Theorem \ref{T:perform}). As pointed out in
  \cite{GA}, one of the advantages of the LMI-approach to the
  state-space $H^\infty$ problem, even in the classical setting, is
  that it allows one to seek controllers that solve the scaled
  $H^\infty$-problem with a given maximal order.
  Indeed, it is shown in \cite{GA,AG} (see also \cite{DP}) that
  certain additional rank constraints on the solutions $X$ and $Y$
  of the LMIs \eqref{Y-LMI} and \eqref{X-LMI} enforce the existence of
  a solution with a prescribed maximal order. However, these
  additional constraints destroy the convexity of the solution
  criteria, and are therefore usually not considered as a desirable
  addition.

  An important point in the application of Finsler's lemma in the derivation
  of the LMI solution criteria in Theorems \ref{T:stab} and
  \ref{T:perform} is that the closed-loop system matrix $A_{cl}$ in
  \eqref{Acl'} has an affine expression in terms of the unknown design
  parameters $\{A_{K},B_{K},C_{K},D_{K}\}$.  This is the key point where
  the assumption $D_{22}=0$ is used.  A parallel simplification
  occurs in the frequency-domain setting where the assumption
  $G_{22} = 0$ leads to the Model-Matching form.  The
  distinction however is that the assumption $G_{22}=0$ is considered
  unattractive from a physical point of view while the parallel
  state-space assumption $D_{22}: = G_{22}(0) = 0$ is considered
  innocuous.

  There is a whole array of lemmas of Finsler type; we have only
  mentioned the form most suitable for our application.  It turns out
  that these various Finsler lemmas are closely connected with the
  theory of plus operators and Pesonen operators on an indefinite
  inner product space (see \cite{Bognar}).
  An engaging historical survey on all the Finsler's lemmas
  is the paper of Uhlig \cite{Uhlig}.

  The notions of modally detectable and modally stabilizable introduced
  in Subsection \ref{S:ND-equiv} along with Theorem \ref{T:NDfreq-state}
  seem new, though of somewhat limited use because it is not known if every
  realization can be reduced to a modally detectable and modally stabilizable
  realization. We included the result as an
  illustration of the difficulties with realization theory for $N$-D transfer functions.

  We note that the usual proof of Lemma \ref{L:ND-sing} for the
  classical 1-D case uses the pole-shifting characterization of
  stabilizability/detectability (see \cite[Exercise 2.19]{DP}).
  The proof here using the Hautus
  characterization of stabilizability/detectability provides a
  different proof for the 1-D case.

  \section{Robust control with structured uncertainty: the
  commutative case}  \label{S:com-robust}\setcounter{equation}{0}

  In the analysis of 1-D control systems, an issue is the uncertainty
  in the plant parameters.  As a control goal, one wants the control
  to achieve internal stability (and perhaps also performance) not
  only for the nominal plant $G$ but also for a whole prescribed family of plants
  containing the nominal plant $G$.

A question then is whether the controller can or cannot have (online) access to the
uncertainty parameters. In a state-space context it is possible to find sufficient
conditions for the case that the controller cannot access the uncertainty parameters,
with criteria that are similar to those found in Theorems \ref{T:stab} and
\ref{T:perform} but additional rank constraints need to be imposed as well, which
destroys the convex character of the solution criterion.
The case where the controller can have access to the uncertainty parameters is usually
given the interpretation of gain-scheduling, and fits better with the multidimensional
system problems discussed in Section \ref{S:ND}. In this section we discuss three
formulations of 1-D control systems with uncertainty in the plant parameters, two of
which can be given gain-scheduling interpretation, i.e., the controller has access to
the uncertainty parameters, and one where the controller is not allowed to use the
uncertainty parameters.

    \subsection{Gain-scheduling in state-space coordinates}\label{S:GS-SSC}

  Following \cite{Packard},
  we suppose that we are given a standard linear time-invariant
       input/state/output system
\begin{equation}\label{GS-system}
\Sigma \colon \left\{ \begin{array}{ccc}
    x(t+1) & = & A_{M}(\delta_{U}) x(t)+ B_{M1}(\delta_{U}) w(t) +
    B_{M2}(\delta_{U}) u(t) \\
    z(t) & = & C_{M1}(\delta_{U}) x(t) + D_{M11}(\delta_{U}) w(t) +
    D_{M12}(\delta_{U}) u(t)  \\
    y(t) & = & C_{M2}(\delta_{U}) x(t) + D_{M21}(\delta_{U}) w(t)
    + D_{M22}(\delta_{U}) u(t)
    \end{array}  \right.\ (t\in\BZ_+)
\end{equation}
    but where the system matrix
    $$\begin{bmatrix}  A_{M}(\delta_{U}) & B_{M1}(\delta_{U}) &
    B_{M2}(\delta_{U}) \\ C_{M1}(\delta_{U})  & D_{M11}(\delta_{U}) &
    D_{M12}(\delta_{U}) \\ C_{M2}(\delta_{U}) & D_{M21}(\delta_{U}) &
    D_{M22}(\delta_{U}) \end{bmatrix}:\mat{c}{\cX\\\cW\\\cU}\to\mat{c}{\cX\\\cZ\\\cY}
    $$
    is not known exactly but depends on some
    {\em uncertainty} parameters $\delta_{U}  = (\delta_{1}, \dots,
    \delta_{d})$ in $\BC^d$.  Here the quantities $\delta_{i}$ are viewed as
       uncertain parameters which the controller can measure and
       use in real time.  The goal is to design
       a controller $\Sigma_{K}$ (independent of $\delta_U$)
       off-line
       so that the closed-loop system (with the controller
       accessing the current values of the varying parameters
       $\delta_{1}, \dots, \delta_{d}$ as well as the value of
       the measurement signal $y$ from the plant) has desirable
       properties for all admissible values of $\delta_U$, usually
       normalized to be $|\delta_{k}| \le 1$ for $k = 1, \dots, d$.
       %It may be desirable to allow some full blocks
%      ( $\Delta_{k}$  equal to an arbitrary $n_{k} \times n_{k}$
%      matrix of norm at most 1 rather than the scalar block
%      $\delta_{k} I_{n_{k}}$), but we stick with scalar blocks
%      for consistency and to relate more easily to the analysis
%      in the previous sections.

    The transfer function for the uncertainty
    parameter $\delta_{U}$ can be expressed as
    \begin{align}
     G(\delta) & = \begin{bmatrix} D_{M11}(\delta_{U}) &
    D_{M12}(\delta_{U}) \\ D_{M21}(\delta_{U}) &
    D_{M22}(\delta_{U})\end{bmatrix}  \notag \\
       &   \qquad
       +  \lambda \begin{bmatrix} C_{M1}(\delta_{U}) \\ C_{M2}(\delta_{U})
  \end{bmatrix} (  I_{\cX} - \lambda A_{M}(\delta_{U}))^{-1} \begin{bmatrix}
  B_{M1}(\delta_{U}) & B_{M2}(\delta_{U}) \end{bmatrix} \label{Utransfunc}
       \end{align}
where we have introduced the aggregate variable
       $$  \delta = (\delta_{U}, \lambda) = (\delta_{1}, \dots,
       \delta_{d}, \lambda).
       $$

    It is not too much of a restriction to assume in
    addition that the functional dependence on $\delta_{U}$ is given by a
    linear fractional map (where the subscript $U$ suggests {\em
    uncertainty} and the subscript $S$ suggests {\em shift})
 \[
\begin{array}{l}
\mat{ccc}{A_{M}(\delta_{U}) & B_{M1}(\delta_{U}) &
       B_{M2}(\delta_{U})
  \\ C_{M1}(\delta_{U})  & D_{M11}(\delta_{U}) &
       D_{M12}(\delta_{U}) \\ C_{M2}(\delta_{U}) & D_{M21}(\delta_{U}) &
       D_{M22}(\delta_{U})}=
\mat{ccc}{
       A_{SS} & B_{S1} & B_{S2} \\
       C_{1S} & D_{11} & D_{12} \\
       C_{2S} & D_{21} & D_{22}}+\\[.5cm]
\qquad\qquad\qquad\qquad+\mat{c}{A_{SU}\\C_{1U}\\C_{2U}}(I-Z(\de_U)A_{UU})^{-1}
Z(\de_U)\mat{ccc}{A_{US} & B_{U1} & B_{U2}},
       \end{array}
\]
where $Z(\de_U)$ is defined analogously to $Z(z)$ in \eqref{Z(z)} relative to
a given decomposition of the ``uncertainty'' state-space
$\cX_{U}=\cX_{U,1}\oplus\cdots\oplus\cX_{U,d}$
on which that state operator $A_{UU}$ acts. In that case the transfer function $G(\delta)$
admits a state-space realization
 \begin{equation}\label{GrealGS}
     G(\de) = \begin{bmatrix} G_{11} & G_{12} \\ G_{21} & G_{22}
 \end{bmatrix} =
     \begin{bmatrix} D_{11} & D_{12} \\ D_{21} & D_{22} \end{bmatrix} +
     \begin{bmatrix} C_{1} \\C_{2} \end{bmatrix}  (I - Z(\de)A)^{-1} Z(\de)
         \begin{bmatrix} B_{1} & B_{2} \end{bmatrix}
   \end{equation}
with system matrix given by
      \begin{equation} \label{aggsys}
    \begin{bmatrix} A & B_{1} & B_{2} \\ C_{1} & D_{11} & D_{12} \\
      C_{2} & D_{21} & D_{22} \end{bmatrix} =
       \left[ \begin{array}{cc|cc}
     A_{UU} & A_{US} & B_{U1} & B_{U2} \\
     A_{SU} & A_{SS} & B_{S1} & B_{S2} \\
    \hline  C_{1U} & C_{1S} & D_{11} & D_{12} \\
     C_{2U} & C_{2S} & D_{21} & D_{22} \end{array} \right].
      \end{equation}
Here $Z(\de)$ is again defined analogously to \eqref{Z(z)} but now on the
extended state-space $\cX_{ext}=\cX_{U}\oplus\cX$.

We can then consider this gain-scheduling problem as a problem of the
constructed $N$-D system (with $N=d+1$), and seek for a controller $K$
with a state-space realization
\begin{equation}\label{KrealGS}
K(\de)=D_K+C_K(I-Z_K(\de)A_K)Z_K(\de)B_K
\end{equation}
so that the closed loop system has desirable properties from a gain-scheduling
perspective. Making a similar decomposition of the system matrix for the
controller $K$ as in \eqref{aggsys}, we note that $K(\de)$ can also be written as
\[
K(\de)=D_{M,K}(\de_U)+\la C_{M,K}(\de_U)(I-\la A_{M,K}(\de_U))^{-1}B_{M,K}(\de_U),
\]
where $A_{M,K}(\de_U)$, $B_{M,K}(\de_U)$, $C_{M,K}(\de_U)$ and $D_{M,K}(\de_U)$
appear as the transfer functions of $N$-D systems (with $N=d$), that is, $K(\de)$
can be seen as the transfer function of a linear time-invariant input/state/output
system
$$ \Si_K: \left\{ \begin{array}{ccc}
       x_K(t+1) & = & A_{M,K}(\delta_{U}) x_K(t) + B_{M,K}(\delta_{U})  u(t) \\
       u(t) & = & C_{M,K}(\delta_{U}) x_K(t) + D_{M,K}(\delta_{U}) y(t)
       \end{array}
       \right.\quad (n\in\BZ_+)
       $$
depending on the same uncertainty parameters $\de_U=(\de_1,\ldots,\de_d)$ as the
system $\Si$.

Similarly, the transfer function $G_{cl}(\de)$ of the closed-loop system with
system matrix $\sbm{A_{cl}&B_{cl}\\C_{cl}&D_{cl}}$ as defined in
\eqref{clsysmat} also can be written as a transfer matrix
\[
G_{cl}(\de)=D_{M,cl}(\de_U)+\la C_{M,cl}(\de_U)(I-\la A_{M,cl}(\de_U))^{-1}B_{M,cl}(\de_U)
\]
with $A_{M,cl}(\de_U)$, $B_{M,cl}(\de_U)$, $C_{M,cl}(\de_U)$ and $D_{M,cl}(\de_U)$
transfer functions of $N$-D systems (with $N=d$), and the corresponding linear
time-invariant input/state/output system
$$\Si_{cl}: \left\{ \begin{array}{ccc}
       x(t+1) & = & A_{M,cl}(\delta_{U}) x(t) + B_{M,cl}(\delta_{U})  w(t) \\
       z(t) & = & C_{M,cl}(\delta_{U}) x(t) + D_{M,cl}(\delta_{U}) w(t)
       \end{array}
       \right.\quad(n\in\BZ_+)
       $$
also appears as the closed-loop system of $\Si$ and $\Si_K$.

It then turns out that stability of $A_{cl}$, that is,
$I-Z_{cl}(\de)A_{cl}$ invertible for all $\de$ in $\ov{\BD}^{d+1}$
(with $Z_{cl}$ as defined in Subsection \ref{S:ND-statespace})
corresponds precisely to {\em robust stability of $\Si_{cl}$}, i.e.,
the spectral radius of $A_{M,cl}(\de_U)$ is less than 1 for all
$\de_U=(\de_1,\ldots,\de_d)$ so that $|\de_k|\leq 1$ for
$k=1,\ldots,d$, and $K$ with realization \eqref{KrealGS} solves the
state-space $H^\infty$-problem for $G$ with realization
\eqref{GrealGS} means that the closed loop system $\Si_{cl}$ has
{\em robust performance}, i.e., $\Si_{cl}$ is robustly stable and
the transfer function $G_{cl}$ satisfies
\[
\|G_{cl}(\de)\|\leq 1\text{ for all }\de=(\de_1,\ldots,\de_d,\la)\in\ov{\BD}^{d+1}.
\]
We may thus see the state-space formulation of the gain-scheduling problems considered
in this subsection as a special case of the $N$-$D$ system stabilization and $H^\infty$-problems
of Subsection \ref{S:ND-statespace}. In particular, the sufficiency analysis
given there, and the results of Theorem \ref{T:stab} and \ref{T:perform}, provide
practical methods for obtaining solutions. As the conditions are only sufficient,
solutions obtained in principle may be conservative.

       \subsection{Gain-scheduling: a pure frequency-domain
       formulation}  \label{S:GS-freq}

       In the approach of Helton (see \cite{HeltonACC, HeltonIEEE}),
       one eschews transfer functions and state-space coordinates completely
       and supposes that one is given a plant $G$ whose frequency
       response depends on a load with frequency function $\delta(z)$
       at the discretion of the user; when the load $\delta$ is
       loaded onto $G$, the resulting frequency-response function has
       the form $G(z, \delta(z))$ where $G = G(\cdot, \cdot)$ is a
       function of two variables.  The control problem (for the
       company selling this device $G$ to a user) is to design the
       controller $K = K(\cdot, \cdot)$ so that $K(\cdot,
       \delta(\cdot))$ solves the $H^{\infty}$-problem for the plant
       $G(\cdot, \delta(\cdot))$.  The idea here is that once the
       user loads $\delta$ onto $G$ with known frequency-response
       function, he is also to load  $\delta$ onto the controller $K$
       (designed off-line); in this way the same controller works for
       many customers using many different $\delta$'s.  When the dust
       settles, this problem reduces to the frequency-domain problem
       posed in Section \ref{S:ND-freq} with $d=2$; an application of the
       Youla-Ku\v{c}era parametrization (or simply using the
       function $Q(z) = K(z) (I - G_{22}(z) K(z))^{-1}$ if the plant
       $G$ itself is stable) reduces the problem of designing the
       control $K$ to a Nevanlinna-Pick-type interpolation problem on
       the bidisk.

       \subsection{Robust control with a hybrid frequency-domain/state-space formulation}
       \label{S:hybrid}

We now consider a hybrid frequency-domain/state-space formulation of
the problem considered in Subsection \ref{S:GS-SSC}; the main
difference is that in this case the controller is not granted access
to the uncertainty parameters.

Assume we are given a 1-D-plant $G(\la)$ that depends on uncertainty parameters
$\de_U=(\de_1,\ldots,\de_d)$ via the linear fractional representation
\begin{align}
 & G(\delta_U, \lambda) =\begin{bmatrix} G_{11}(\lambda) & G_{12}(\lambda) \\
      G_{21}(\lambda) & G_{22}(\lambda) \end{bmatrix}+\notag\\
& \qquad\qquad
      + \begin{bmatrix} G_{1U}(\lambda) \\ G_{2U}(\lambda)
  \end{bmatrix} (I - Z(\de_U) G_{UU}(\lambda))^{-1}
  Z(\de_U) \begin{bmatrix} G_{U1}(\lambda) &
  G_{U2}(\lambda) \end{bmatrix}
  \label{Gdeltalam}
  \end{align}
with $Z(\de_U)$ as defined in Subsection \ref{S:GS-SSC}, and
where the coefficients are 1-D-plants independent of $\de_U$:
$$
       G_{aug}(\la) = \begin{bmatrix} G_{UU}(\lambda) & G_{U1}(\lambda)
       & G_{U2}(\lambda) \\
       G_{1U}(\lambda) & G_{11}(\lambda) & G_{12}(\lambda) \\
       G_{2U}(\lambda) & G_{21}(\lambda) & G_{22}(\lambda)
       \end{bmatrix}:\mat{c}{\cX_{U}\\\cW\\\cU}\to\mat{c}{\cX_{U}\\\cZ\\\cY}.
$$
In case $G_{aug}(\la)$ is also given by a state-space realization,
we can write $G(\de_U,\la)$ as in \eqref{GrealGS} with
$\de=(\de_U,\la)$ and $Z(\de)$ acting on the extended state-space
$\cX_{ext}=\cX_U\oplus\cX$.

For this variation of the gain-scheduling problem we seek to design a controller
$K(\la)$ with matrix values representing operators from $\cY$ to $\cU$ so that $K$
solves the $H^{\infty}$-problem for $G(\delta_U,\la)$ for every $\delta_U$ with
$\| Z(\de_U) \| \leq1$, i.e., $|\de_j|\leq 1$ for $j=1,\ldots,d$. For the sequel it is
convenient to assume that $\cZ=\cW$. In that case, using the Main Loop Theorem
\cite[Theorem 11.7 page  284]{ZDG}, it is easy to see that this problem
can be reformulated as: {\em Find  a single-variable transfer matrix $K(\cdot)$ so that
  $\Theta(\widetilde G, K)$ given by \eqref{ThetaGK}, with
  $\widetilde G = \sbm{
  \widetilde G_{11} & \widetilde G_{12} \\
  \widetilde G_{21} & \widetilde G_{22}}$ in \eqref{ThetaGK} taken to be
  $$
  \begin{bmatrix} \wtilG_{11}(\lambda) & \wtilG_{12}(\lambda) \\
      \wtilG_{21}(\lambda) & \wtilG_{22}(\lambda) \end{bmatrix}
      =
      \left[ \begin{array}{cc|c}
      G_{UU}(\lambda) & G_{U1}(\lambda)
           & G_{U2}(\lambda) \\
           G_{1U}(\lambda) & G_{11}(\lambda) & G_{12}(\lambda) \\
           \hline
           G_{2U}(\lambda) & G_{21}(\lambda) & G_{22}(\lambda)
       \end{array} \right],
    $$
    is stable and such that}
   $$\mu_{\DS} \left(
   \widetilde G_{11}(\lambda) + \widetilde G_{12}(\lambda)
   (I - K(\lambda) \widetilde G_{22}(\lambda))^{-1} K(\lambda)
   \widetilde G_{21}(\lambda) \right) < 1.
   $$
Here $\mu_\DS$ is as defined in \eqref{mu} with
$\DS$ the $C^*$-algebra
\[
\DS=\left\{\mat{cc}{Z(\de_U)&0\\0&T}\colon \de_U\in\BC^d,\ T\in\cL(\cZ)\right\}
\subset\cL(\cX_U\oplus\cZ).
\]

Application of the Youla-Ku\v{c}era parametrization of the controllers $K$ that
stabilize $\Th(\wtilG,K)$ as in Subsection \ref{S:HinfToMM} converts the problem
to the following: {\em Given stable 1-variable transfer functions $T_{1}(\lambda)$,
$T_{2}(\lambda)$, and $T_{3}(\lambda)$ with matrix values representing operators
in the respective spaces
   $$\cL(\cX_U \oplus \cW, \cX_U   \oplus  \cZ), \quad
   \cL(\cX_U \oplus \cU, \cX_U \oplus \cZ), \quad
   \cL(\cX_U \oplus \cW, \cX_U \oplus \cY),
   $$
   find a stable 1-variable transfer function $\Lambda(\lambda)$ with
  matrix values representing operators in $\cL(\cX_U\oplus\cY,\cX_U\oplus \cU)$ so that the
  transfer function $S(\lambda)$ given by
  \begin{equation} \label{MMform}
      S(\lambda) = T_{1}(\lambda) + T_{2}(\lambda)  \Lambda(\lambda)
  T_{3}(\lambda)
  \end{equation}
  has $\mu_{\DS}(S(\lambda)) < 1$ for all $\lambda \in {\mathbb D}$.}
If $T_{2}(\zeta)$ and $T_{3}(\zeta)$ are square and invertible for
  $\zeta$ on the boundary ${\mathbb T}$ of the unit disk ${\mathbb
  D}$, the model-matching form  \eqref{MMform} can be converted to
  bitangential interpolation conditions (see e.g.~\cite{BGR}); for
  simplicity, say that these interpolation conditions have the form
  \begin{equation}   \label{interpolation} x_{i} S(\lambda_{i}) =
      y_{i}, \quad  S(\lambda_{j}') u_{j} = v_{j}
  \text{ for } i = 1, \dots, k, \quad j = 1, \dots, k'
  \end{equation}
  for given distinct points $\lambda_{i}, \lambda'_{j}$ in ${\mathbb
  D}$, row vectors $x_{i}, y_{i}$ and column vectors $u_{j}, v_{j}$.
  Then the robust $H^{\infty}$-problem ($H^{\infty}$ rather than
  rational version) can be converted to the
  $\mu$-Nevanlinna-Pick problem:   {\em find holomorphic function $S$
  on the unit disk with matrix values representing operators in
  $\cL(\cX_U \oplus \cW, \cX_U \oplus \cZ)$ satisfying the
  interpolation conditions \eqref{interpolation} such that also}
  $$
   \mu_{\DS}(S(\lambda)) < 1 \text{ for
   all } \lambda \in {\mathbb D}.
  $$

  It is this $\mu$-version of the Nevanlinna-Pick interpolation
  problem which has been studied from various points of view
  (including novel variants of the Commutant Lifting Theorem) by
  Bercovici-Foias-Tannenbaum (see \cite{BFT1,BFT2,BFT3,BFT4}) and Agler-Young
  (see \cite{AY1,AY3,AY5,AY7} and Huang-Marcantognini-Young \cite{HMY}).
  These authors actually
  study only very special cases of the general control problem as
  formulated here; hence the results at this stage are not
  particularly practical for actual control applications.  However
  this work has led to interesting new mathematics in a number of
  directions: we mention in particular the work of Agler-Young on new
  types of dilation theory and operator-model theory (see \cite{AY2,
  AY5}),
  new kinds of realization theorems \cite{AY6},
  the complex geometry of new kinds of domains in
  ${\mathbb C}^{d}$ (see \cite{AY4, AY8, AY9}), and a multivariable extension of
  the Bercovici-Foias-Tannenbaum spectral commutant lifting theorem
  due to Popescu \cite{Po2}.

 \subsection{Notes}  \label{S:com-robust:notes}

 In the usual formulation of $\mu$ (see \cite{PackardDoyle, ZDG}), in
 addition to the scalar blocks $\delta_{i}I_{n_{i}}$ in $Z(\delta)$,
 it is standard to also allow some of the blocks to be full blocks of the
 form $\Delta_{i} = \sbm{\delta^{(i)}_{11} & \cdots &
 \delta^{(i)}_{1n_{i}} \\ \vdots & & \vdots \\ \delta^{(i)}_{n_{i}1}
 & \cdots& \delta^{(i)}_{n_{i}n_{i}} }$.  The resulting transfer
 functions then have domains equal to be (reducible) Cartan domains
 which are more general than the unit polydisk.  The theory of the
 Schur-Agler class has been extended to this setting in
 \cite{AT, BallBoloJFA}.  More generally, it is natural also to
 allow non-square blocks.  A formalism for handling this is given in
 \cite{BGM3}; for this setting one must work with the intertwining
 space of $\Delta$ rather than the commutant of $\Delta$ in the
 definition of $\widehat \mu$ in \eqref{muhat}.  With a formalism for
 such a non-square uncertainty structure available, one can avoid the
 awkward assumption in Subsection \ref{S:hybrid} and elsewhere that
 $\cW = \cZ$.

 \section{Robust control with dynamic time-varying structured uncertainty}
\label{S:NC}\setcounter{equation}{0}

 \subsection{The state-space LFT-model formulation}
  \label{S:NC-statespace}

  Following \cite{LZDProc, LZD, LuThesis, Paganini}, we now introduce
  a variation on the gain-scheduling problem discussed in Section
  \ref{S:GS-SSC} where the uncertainty parameters $\delta_{U} =
  (\delta_{1}, \dots, \delta_{d})$ become operators on $\ell^{2}$, the space of square-summable
sequences of complex numbers indexed by the integers ${\mathbb Z}$,
  and are to be interpreted as dynamic, time-varying uncertainties.  To make
  the ideas precise, we suppose that we are given a system matrix as
  in \eqref{aggsys}.  We then tensor all operators with the identity
  operator $I_{\ell^{2}}$ on $\ell^{2}$ to obtain an enlarged system matrix
  \begin{equation}   \label{aggsys-en}
     \bM=\begin{bmatrix} A & B_{1} & B_{2} \\ C_{1} & D_{11} & D_{12} \\
      C_{2} & D_{21} & D_{22} \end{bmatrix} \otimes I_{\ell^{2}} =
       \left[ \begin{array}{cc|cc}
     A_{UU} & A_{US} & B_{U1} & B_{U2} \\
     A_{SU} & A_{SS} & B_{S1} & B_{S2} \\
    \hline  C_{1U} & C_{1S} & D_{11} & D_{12} \\
     C_{2U} & C_{2S} & D_{21} & D_{22} \end{array} \right]
     \otimes I_{\ell^{2}},
    \end{equation}
which we also write as
\begin{equation}   \label{aggsys-en2}
      \bM=\begin{bmatrix} \bA & \bB_{1} & \bB_{2} \\ \bC_{1}
    & \bD_{11} & \bD_{12} \\ \bC_{2} & \bD_{21} & \bD_{22}
    \end{bmatrix}
    \colon \begin{bmatrix} (\cX_{U}\oplus \cX_{S}) \otimes \ell^{2} \\
    \cW \otimes \ell^{2} \\ \cU \otimes
    \ell^{2} \end{bmatrix} \to \begin{bmatrix} (\cX_{U}\oplus \cX_{S} )\otimes
    \ell^{2} \\ \cZ \otimes \ell^{2} \\
    \cY \otimes \ell^{2} \end{bmatrix}.
    \end{equation}
    Given a decomposition $\cX_U=\cX_{U1}\oplus\cdots\oplus\cX_{Ud}$ of the uncertainty
    state space $\cX_U$, we define the matrix pencil $\bZ_U(\bdelta_U)$ with argument
    equal to a $d$-tuple $\bdelta_U = (\bdelta_{1}, \dots, \bdelta_{d})$ of (not
    necessarily commuting) operators on $\ell^{2}$ by
    $$
    \bZ_U(\bdelta_U) = \begin{bmatrix} I_{\cX_{U1}} \otimes \bdelta_{1} &
    & \\ & \ddots & \\ & & I_{\cX_{Ud}} \otimes \bdelta_{d}
    \end{bmatrix},
    $$
    In addition we let $\bS$ denote the bilateral shift operator on
    $\ell^{2}$; we sometimes will also view $\bS$ as an operator on the space
    $\ell$ of all sequences of complex numbers or on the subspace $\ell^2_{\tu{fin}}$
    of $\ell^2$ that consists of all sequences in $\ell^2$ with finite support.
    We obtain an uncertain linear system of the form
\begin{equation}\label{NC-system}
\bSigma \colon \left\{ \begin{array}{ccc}
     \bS^*\vec{x} & = & A_{M}(\bdelta_{U}) \vec{x}+ B_{M1}(\bdelta_{U}) \vec{w} +
     B_{M2}(\bdelta_{U}) \vec{u} \\
     \vec{z} & = & C_{M1}(\bdelta_{U}) \vec{x} + D_{M11}(\bdelta_{U}) \vec{w} +
     D_{M12}(\bdelta_{U}) \vec{u} \\
     \vec{y} & = & C_{M2}(\bdelta_{U}) \vec{x} + D_{M21}(\bdelta_{U}) \vec{w}
     + D_{M22}(\bdelta_{U}) \vec{u}
     \end{array}  \right.
\end{equation}
     where the system matrix
     $$\begin{bmatrix}  A_{M}(\bdelta_{U}) & B_{M1}(\bdelta_{U}) &
    B_{M2}(\bdelta_{U}) \\ C_{M1}(\bdelta_{U})  & D_{M11}(\bdelta_{U}) &
    D_{M12}(\bdelta_{U}) \\ C_{M2}(\bdelta_{U}) & D_{M21}(\bdelta_{U}) &
    D_{M22}(\bdelta_{U}) \end{bmatrix}
    :\mat{c}{\cX_S\otimes \ell^{2}_{\tu{fin}}\\\cW \otimes \ell^{2}_{\tu{fin}}
    \\\cU \otimes \ell^{2}_{\tu{fin}}}
    \to\mat{c}{\cX_S\otimes \ell\\\cZ \otimes \ell \\\cY \otimes \ell}
    $$
is obtained from the feedback connection
\[
\mat{c}{\vec{\wtil{x}}_U\\\vec{\wtil{x}}_S\\\vec{z}\\\vec{y}}
=\bM\mat{c}{\vec{x}_U\\\vec{x}_S\\\vec{z}\\\vec{y}},\quad
\text{subject to}\quad \vec{x}_U=\bZ_U(\bdelta_{U})\vec{\wtil{x}}_U,
\]
that is,
\begin{equation}   \label{LFTform}
 \begin{array}{l}
 \mat{ccc}{\bA_{M}(\bdelta_{U}) & \bB_{M1}(\bdelta_{U}) &
    \bB_{M2}(\bdelta_{U})
   \\ \bC_{M1}(\bdelta_{U})  & \bD_{M11}(\bdelta_{U}) &
    \bD_{M12}(\bdelta_{U}) \\ \bC_{M2}(\bdelta_{U}) & \bD_{M21}(\bdelta_{U}) &
    \bD_{M22}(\bdelta_{U})}=
 \mat{ccc}{
    \bA_{SS} & \bB_{S1} & \bB_{S2} \\
    \bC_{1S} & \bD_{11} & \bD_{12} \\
    \bC_{2S} & \bD_{21} & \bD_{22}}+\\[.5cm]
 \qquad\qquad+\mat{c}{\bA_{SU}\\\bC_{1U}\\\bC_{2U}}(I-\bZ_U(\bdelta_U) \bA_{UU})^{-1}
 \bZ_U(\bdelta_U)\mat{ccc}{\bA_{US} & \bB_{U1} & \bB_{U2}}.
    \end{array}
 \end{equation}

As this system is time-varying, due to the presence of the time-varying
  uncertainty parameters $\bdelta_{U}$, it is not convenient to work
  with a transfer-function acting on the frequency-domain; instead we
  stay in the time-domain and work with the input-output operator
  which has the form
  \begin{align}
      \bG(\bdelta) & = \begin{bmatrix} \bD_{M11}(\bdelta_{U}) &
     \bD_{M 12}(\bdelta_{U}) \\ \bD_{M 21}(\bdelta_{U}) &
     \bD_{M22}(\bdelta_{U})\end{bmatrix} +   \begin{bmatrix} \bC_{M1}(\bdelta_{U}) \\ \bC_{M2}(\bdelta_{U})
   \end{bmatrix}\times \label{UIO} \\
    &   \qquad\times
     (  I_{\cX_S\otimes\ell^2} -  (I_{\cX_S}\otimes\bS) \bA_{M}(\bdelta_{U}))^{-1}
   (I_{\cX_S}\otimes\bS)
   \begin{bmatrix}\bB_{M1}(\bdelta_{U}) & \bB_{M2}(\bdelta_{U}) \end{bmatrix},
\notag
    \end{align}

Now write $\bdelta$  for  the collection $(\bdelta_{U}, \bS)$ of $d+1$ operators on $\ell^{2}$.
Then the input-output operator $\bG(\bdelta)$ given by \eqref{UIO}
 has the noncommutative transfer-function realization
 \begin{equation}\label{Greal-noncom}
     \bG(\bdelta) = \begin{bmatrix} \bG_{11}(\bdelta) &
     \bG_{12}(\bdelta) \\ \bG_{21}(\bdelta) &
     \bG_{22}(\bdelta) \end{bmatrix} =
     \begin{bmatrix} \bD_{11} & \bD_{12} \\ \bD_{21} & \bD_{22} \end{bmatrix} +
     \begin{bmatrix} \bC_{1} \\ \bC_{2} \end{bmatrix}  (I -
     \bZ(\bdelta) \bA)^{-1} \bZ(\bdelta)
     \begin{bmatrix} \bB_{1} & \bB_{2} \end{bmatrix}
   \end{equation}
   with system matrix as in \eqref{aggsys-en} and
   $\bZ(\bdelta)=\sbm{\bZ_U(\bdelta_U)&0\\0&I_{\cX_S}\otimes\bS}$.
   In the formulas \eqref{LFTform}-\eqref{Greal-noncom} the inverses may have to be
   interpreted as the algebraic inverses of the corresponding infinite block matrices;
   in that way, the formulas make sense at least for the nominal plant, i.e.,
   with $\bdelta_U=(0,\ldots,0)$.

   More generally, the transfer-function $\bG$ can be extended to a function of
   $d+1$ variables in $\cL(\ell^2)$ by replacing $\bS$ with another variable
   $\de_{d+1}\in \cL(\ell^2)$. In that case, the transfer-function can be viewed as
   an LFT-model with structured uncertainty, as studied in \cite{LZD,DP}. However,
   as a consequence of the Sz.-Nagy dilation theory, without loss of generality
   it is possible in this setting
   of LFT-models to  fix one of the variables to be the shift
   operator $\bS$; in this way the LFT-model results developed for
   $d+1$ free variable contractions apply equally well to the case of
   interest where one of the variables is fixed to be the shift
   operator.

   Such an input/state/output system $\bSigma$ with structured
   dynamic time-varying uncertainty $\bdelta_{U}$ is said to be {\em
   robustly stable} (with respect to the  dynamic time-varying
   uncertainty structure $\bZ_U(\bdelta_{U})$) if the state-matrix
   $\bA_M(\bdelta_{U})$ is stable for all choices of $\bdelta_{U}$
   subject to $\| \bZ_U(\bdelta_{U}) \| \le 1$, that is, if
   $I_{\cX_S\otimes\ell^2} -  (I_{\cX_S}\otimes\bS) \bA_{M}(\bdelta_{U})$ is
   invertible as an operator on $\cX_S\oplus\ell^2$ for all $\bdelta_{U}$
   with $\| \bZ_U(\bdelta_{U}) \| \le 1$. Since
   $$
   \bA_{M}(\bdelta_U) = \bA_{SS} + \bA_{SU} (I - \bZ_U(\bdelta_{U})
   \bA_{UU})^{-1} \bZ(\bdelta_{U}) \bA_{US},
   $$
   it follows from  the Main Loop Theorem
   \cite[Theorem 11.7 page  284]{ZDG}, that this condition in turn reduces
   to:
   \begin{equation}  \label{NCHautusstab}
    I_{\boldcX} - \bZ(\bdelta) \bA \text{ is invertible for all }
    \bdelta = (\bdelta_{U}, \bS) \text{
    with } \| \bZ(\bdelta) \| \le 1.
   \end{equation}
   Note that this condition amounts to a noncommutative version of
   the Hautus-stability criterion for the matrix $A$ (where $\bA = A
   \otimes I_{\ell^{2}}$). We shall therefore call the state matrix $\bA$
   {\em nc-Hautus-stable} if \eqref{NCHautusstab} is satisfied (with nc indicating
   that we are in the noncommutative setting).
   The input/state/output system $\bSigma$ is said to have {\em nc-performance}
   (with respect to the  dynamic time-varying
      uncertainty structure $\bZ_U(\bdelta_U)$) if it is robustly
      stable (with respect to this dynamic time-varying uncertainty
      structure) and in addition the input-output operator
      $G(\bdelta)$ has norm strictly less than 1 for all
      choices of $\bdelta=(\bdelta_U,\bS)$ with $\| \bZ(\bdelta) \| \le 1$.

One of the key results from the thesis of Paganini \cite{Paganini} which makes the
noncommutative setting
of this section more in line with the 1-D case is that, contrary to what is the case
in Subsection \ref{S:ND-statespace}, for operators $\bA=A\oplus I_{\ell^2}$ on
$\cX\oplus \ell^2$ we do have
$\mu_\DS(\bA)=\widehat\mu_\DS(\bA)$ when we take $\DS$ to be the $C^*$-algebra
\begin{equation}\label{NC-DS}
\DS=\left\{\mat{cc}{\!\!\bZ(\bdelta_U)\!\!\!&0\!\!\!\\0&I_{\cX_S}\otimes\bdelta_{d+1}\!\!}\colon
\bdelta_U=(\bdelta_1,\ldots,\bdelta_d),\ \bdelta_j\in\cL(\ell^2),\ j=1,\ldots,d+1\right\}.
\end{equation}
Write $\bcD$ for the commutant of $\DS$ in $\cL((\cX_U\oplus\cX_S)\otimes\ell^2)$. Then
the main implication of the fact that $\mu_\DS(\bA)=\widehat\mu_\DS(\bA)$
is that nc-Hautus-stability of $\bA$ is now the same as
the existence of an invertible operator $\bQ\in\bcD$ so that $\|\bQ^{-1}\bA\bQ\|<1$ or,
equivalently, the existence of a solution $\bX\in\bcD$ to the LMIs $\bA^*\bX\bA-\bA<0$ and
$\bX>0$.
However, it is not hard to see that $\bX$ is an element of $\bcD$ if and only if
$\bX=X\otimes I_{\ell^2}$ with $X$ being an element of the $C^*$-algebra $\cD$
in \eqref{Zcom}. Thus, in fact, we find that $\bA=A\oplus I_{\ell^2}$ is nc-Hautus-stable
precisely when $A$ is scaled stable, i.e., when there exists a solution
$X\in\cD$ to the LMIs $A^*XA-A<0$ and $X>0$.

These observations can also be seen as a special case (when $C_2=0$ and $B_2=0$)
of the following complete analogue of Theorem \ref{T:1.3} for this
noncommutative setting due to Paganini \cite{Paganini}.

%%%%%%%%%%%%%%%%%%%%%%%%%%%%%%%%%%%%%%%%%%%%%%%%%%%%%%%%%%%%%%%%%%%%%%%%
\begin{proposition}\label{P:NCdetect}
Given a system matrix as in \eqref{aggsys-en}-\eqref{aggsys-en2}, then:
\begin{itemize}
\item[(i)] The output pair $\{\bC_2,\bA\}$ is nc-Hautus-detectable, that is,
for every $\bdelta=(\bdelta_1,\ldots,\bdelta_{d+1})$, with $\bdelta_j\in\cL(\ell^2)$
for $j=1,\ldots,d+1$, so that $\|\bZ(\bdelta)\|\leq1$ the operator
\[
\mat{c}{I-\bZ(\bdelta)\bA\\\bC_2}:
(\cX_U\oplus\cX_S)\otimes\ell^2\to\mat{c}{(\cX_U\oplus\cX_S)\otimes\ell^2\\\cY\oplus\ell^2}
\]
has a left inverse, if and only if $\{\bC_2,\bA\}$ is nc-operator-detectable, i.e.,
there exists an operator $\bL=L\otimes I_{\ell^2}$, with $L:\cY\to\cX$, so that
$\bA+ \bL \bC_2$ is nc-Hautus-stable,  if and only if there exists a
solution $X\in\cD$ to the LMIs
\begin{equation}\label{X-stableLMI2}
A^*XA-X-C_2^*C_2<0,\quad X>0.
\end{equation}

\item[(ii)] The input pair $\{\bA,\bB_2\}$ is nc-Hautus-stabilizable, that is,
for every $\bdelta=(\bdelta_1,\ldots,\bdelta_{d+1})$, with $\bdelta_j\in\cL(\ell^2)$
for $j=1,\ldots,d+1$, so that $\|\bZ(\bdelta)\|\leq1$ the operator
\[
\mat{cc}{I-\bZ(\bdelta)\bA&\bB_2}:\mat{c}{(\cX_U\oplus\cX_S)\otimes\ell^2\\\cU\oplus\ell^2}\to
(\cX_U\oplus\cX_S)\otimes\ell^2
\]
has a left inverse, if and only if $\{\bA, \bB_{2}\}$ is nc-operator-stabilizable, i.e.,
there exists an operator $\bF=F\otimes I_{\ell^2}$, with $F:\cX\to\cU$, so that
$\bA+\bB_2\bF$ is nc-Hautus-stable, which happens if and only if there exists a solution
$Y\in\cD$ to the LMIs
\begin{equation}\label{Y-stableLMI2}
AYA^*-Y-B_2B_2^*<0,\quad Y>0.
\end{equation}
\end{itemize}
\end{proposition}

In case the input/state/output system $\bSigma$ is not stable and/or does
not have performance, we want to remedy this by means of a feedback with
a controller $\bK$, which we assume has on-line access to the structured
dynamic time-varying uncertainty operators $\bdelta_{U}$ in addition to being dynamic,
i.e., $\bK=\bK(\bdelta)=\bK(\bdelta_U,\bS)$. More specifically, we shall restrict to
controllers of the form
\begin{equation}\label{ncstatecontroller}
 \bK(\bdelta) = \bD_{K} + \bC_{K}(I - \bZ_{K}(\bdelta)
   \bA_{K})^{-1} \bZ_{K}(\bdelta) \bB_{K}
\end{equation}
 where
   $$\bZ_{K}(\bdelta)=\mat{cc}{\bZ_{KU}(\bdelta_U)&0\\0&I_{\cX_{KS}}\otimes \bS},\
   \bZ_{KU}(\bdelta_U) = \begin{bmatrix} I_{\cX_{K1}} \otimes \bdelta_{1}
   & & \\ & \ddots & \\ & & I_{\cX_{Kd}} \otimes \bdelta_{d}
\end{bmatrix},
   $$
   with system matrix $\bM_{\bK}$ of the form
\begin{equation}\label{ncsysmat1}
   \bM_{K} = \begin{bmatrix} \bA_{K} & \bB_{K} \\ \bC_{K} &
   \bD_{K} \end{bmatrix} \colon \begin{bmatrix} (\cX_{KU}\oplus\cX_{KS}) \otimes
   \ell^{2} \\ \cY \otimes \ell^{2} \end{bmatrix} \to
   \begin{bmatrix} (\cX_{KU}\oplus\cX_{KS}) \otimes \ell^{2} \\ \cU \otimes \ell^{2}
       \end{bmatrix}
\end{equation}
where $\cX_{KU}=\cX_{KU1}\oplus\cdots\oplus\cX_{KUd}$, and where the matrix
entries in turn have a tensor-factorization
\begin{equation}\label{ncsysmat2}
 \begin{bmatrix} \bA_{K} & \bB_{K} \\ \bC_{K} & \bD_{K}
    \end{bmatrix} = \begin{bmatrix}  A_{K} \otimes I_{\ell^{2}} &
    B_{K} \otimes I_{\ell^{2}} \\ C_{K} \otimes I_{\ell^{2}} & D_{K}
    \otimes I_{\ell^{2}} \end{bmatrix}.
\end{equation}

    If such a controller $\bK(\bdelta)$ is put in feedback connection with
    $\bG(\bdelta)$, where we impose the usual assumption $D_{22}=0$ to guarantee
    well-posedness, the resulting closed-loop system input-output
    operator $\bG_{cl}(\bdelta)$, as a function
    of the operator uncertainty parameters $\bdelta_{U} =
    (\bdelta_{1}, \dots, \bdelta_{d})$ and the shift $\bS$, has a
    realization which is formally exactly as in \eqref{cltransfer}, that is
    $$
    \bG_{cl}(\bdelta) = \bD_{cl} + \bC_{cl}(I - \bZ_{cl}(\bdelta)
    \bA_{cl})^{-1} \bZ_{cl}(\bdelta) \bC_{cl}
    $$
    with system matrix
    \begin{equation}  \label{clsysmat-NC}
       \begin{bmatrix}  \bA_{cl} & \bB_{cl} \\ \bC_{cl} & \bD_{cl} \end{bmatrix} =
       \left[ \begin{array}{cc|c}
       \bA+\bB_{2} \bD_{K} \bC_{2} & \bB_{2} \bC_{K} & \bB_{1} +
       \bB_{2} \bD_{K} \bD_{21} \\
       \bB_{K} \bC_{2} & \bA_{K} & \bB_{K} \bD_{21} \\
       \hline \bC_{1} + \bD_{12} \bD_{K} \bC_{2} & \bD_{12}
       \bC_{K} & \bD_{11} +
       \bD_{12} \bD_{K} \bD_{21} \end{array} \right],
    \end{equation}
    which is the same as the system matrix \eqref{clsysmat} tensored with $I_{\ell^2}$, and
\begin{equation}\label{Zcl}
  \bZ_{cl}(\bdelta)=\mat{cc}{\bZ(\bdelta)&0\\0&\bZ_{K}(\bdelta)}\quad\text{where}\quad
    \bdelta=(\bdelta_U,\bS).
\end{equation}

    The {\em state-space nc-stabilization problem} (with respect to the
    given dynamic time-varying uncertainty structure $\bdelta_{U}$)
    then is to design a controller $\bK$ with
    state-space realization $\{ \bA_{K}, \bB_{K}, \bC_{K}, \bD_{K}\}$
    as above so that the closed-loop system $\bSigma_{cl}$ defined by the
    system matrix \eqref{clsysmat-NC} is robustly stable.
    The {\em state-space nc-$H^{\infty}$-problem} is to design a controller $\bK$ with
    state-space realization $\{ \bA_{K}, \bB_{K}, \bC_{K}, \bD_{K}\}$
    as above so that the closed-loop system $\bSigma_{cl}$ also has robust performance.

Since the closed-loop state-operator $\bA_{cl}$ is equal to $A_{cl}\otimes I_{\ell^2}$
with $A_{cl}$ defined by \eqref{clsysmat}, it follows as another
implication of the fact that $\mu_\DS$ is equal to $\widehat\mu_\DS$ for operators that
are tensored with $I_{\ell^2}$ (with respect to the appropriate $C^*$-algebra $\DS$)
that $\bA_{cl}$ is nc-Hautus-stable precisely when $A_{cl}$ is scaled stable, i.e.,
we have the following result.

%%%%%%%%%%%%%%%%%%%%%%%%%%%%%%%%%%%%%%%%%%%%%%%%%%%%%%%%%%%%%%%%%%%%%%%%
\begin{proposition}\label{P:NCstable1}
Let $\bSigma$ and $\Si$ be the systems given by \eqref{NC-system} and \eqref{GS-system},
respectively, corresponding to a given system matrix \eqref{aggsys}. Then $\bSigma$ is
nc-Hautus-stabilizable if and only if $\Sigma$ is scaled-stabilizable.
\end{proposition}

Thus, remarkably, the solution criterion given in Section \ref{S:ND-statespace} for the
scaled state-space stabilization problem turns out to be necessary and sufficient for
the solution of the dynamic time-varying structured-uncertainty version of the problem.

%%%%%%%%%%%%%%%%%%%%%%%%%%%%%%%%%%%%%%%%%%%%%%%%%%%%%%%%%%%%%%%%%%%%%%%%
\begin{theorem}\label{T:NCstable2}
Let $\bSigma$ be the system given by \eqref{NC-system} corresponding to a given
system matrix \eqref{aggsys-en}. Then $\bSigma$ is nc-Hautus-stabilizable if and
only if the output pair $\{\bC_2,\bA\}$ is nc-Hautus-detectable and the input pair
$\{\bA,\bB_2\}$ is nc-Hautus-stabiliz\-able, i.e., if there exist solutions $X,Y\in\cD$,
with $\cD$ the $C^*$-algebra given in \eqref{Zcom}, to the LMIs \eqref{X-stableLMI2}
and \eqref{Y-stableLMI2}. In this case $\bK \sim \sbm{A_{K} & B_{K}
\\ C_{K} & D_{K}} \otimes I_{\ell^{2}}$ with $\sbm{A_{K} & B_{K} \\
C_{K} & D_{K}}$ as in \eqref{stabcon} is a controller solving the
nc-Hautus stabilization problem for $\bSigma$.
\end{theorem}

In a similar way, the state-space nc-$H^{\infty}$-problem corresponds to the
scaled $H^\infty$-problem of Subsection \eqref{S:ND-statespace}.

%%%%%%%%%%%%%%%%%%%%%%%%%%%%%%%%%%%%%%%%%%%%%%%%%%%%%%%%%%%%%%%%%%%%%%%%
\begin{theorem}\label{T:NCHinfty}
Let $\bSigma$ be the system given by \eqref{NC-system} for a given system matrix
\eqref{aggsys-en}. Then there exists a solution $\bK$, with realization
\eqref{ncstatecontroller}, to the state-space nc-$H^{\infty}$-problem for the
non-commutative system $\bSigma$ if and only if there exist $X,Y\in\cD$ that
satisfy the LMIs \eqref{X-LMI} and \eqref{Y-LMI} and the coupling condition
\eqref{coupling}.
\end{theorem}

%%%%%%%%%%%%%%%%%%%%%%%%%%%%%%%%%%%%
\begin{proof}
Let $\bSigma$ and $\Si$ be the systems given by \eqref{NC-system} and \eqref{GS-system},
respectively, corresponding to a given system matrix \eqref{aggsys}. Using the strict
bounded real lemma from \cite{BGM3} in combination with similar arguments as used above
for the nc-stabilizability problem, it follows that a transfer-function $\bK$ with
realization \eqref{ncstatecontroller}-\eqref{ncsysmat2} is a solution to the state-space
nc-$H^{\infty}$-problem for $\bSigma$ if and only if the transfer function $K$ with
realization \eqref{statecontroller} is a solution to the scaled $H^\infty$-problem for
the system $\Sigma$. The statement then follows from Theorem \ref{T:perform}.
\end{proof}

 \subsection{A noncommutative frequency-domain formulation}
 \label{S:NC-FD}

    In this subsection we present a frequency-domain version of
    the noncommutative state-space setup of the previous subsection
    used to model linear  input/state/output systems with LFT-model
    for dynamic time-varying structured uncertainty.
    The frequency-domain setup here is analogous to that of Section
    \ref{S:ND-freq} but the unit polydisk $\overline{\mathbb D}^{d}$
    is replaced by the noncommutative polydisk $\overline{\mathbb
    D}^{d}_{nc}$ consisting of all $d$-tuples $\bdelta =
    (\bdelta_{1}, \dots, \bdelta_{d})$ of contraction operators on a
    fixed separable infinite-dimensional Hilbert space $\cK$.

    We need a few preliminary definitions.
    We define $\cF_{d}$ to be the free semigroup consisting of all words
    $\alpha = i_{N} \cdots i_{1}$
    in the letters $\{1, \dots, d\}$. When $\alpha = i_{N} \cdots
    i_{1}$ we write $N = |\alpha|$ for the number of letters in the
    word $\alpha$. The multiplication of two words is given by
    concatenation:
    $$ \alpha \cdot \beta = i_{N} \cdots i_{1} j_{M} \cdots j_{1}
    \text{ if } \alpha = i_{N} \cdots i_{1}  \text{ and } \beta =
    j_{M} \cdots j_{1}.
    $$
    The unit element of $\cF_d$ is the {\em empty word} denoted by $\emptyset$
    with $|\emptyset| = 0$.
    In addition, we let  $z = (z_{1}, \dots, z_{d})$ stands for a $d$-tuple of {\em
    noncommuting} indeterminates, and  for any $\alpha = i_{N} \cdots i_{1}
    \in \cF_{d}-\{\emptyset\}$, we let $z^{\alpha}$ denote the noncommutative
    monomial $z^{\alpha} = z_{i_{N}} \cdots z_{i_{1}}$, while $z^\emptyset=1$.
    If $\alpha$ and $\beta$ are two words in $\cF_{d}$, we multiply the associated
    monomials $z^{\alpha}$ and $z^{\beta}$ in the natural way:
    $$
      z^{\alpha} \cdot z^{\beta} = z^{\alpha\cdot\beta}.
    $$
    Given two Hilbert spaces
    $\cU$ and $\cY$, we let $\cL(\cU, \cY)\langle \langle z \rangle
    \rangle$ denote the
    collection of  all noncommutative formal power series $S(z)$ of the form $S(z) =
    \sum_{\alpha \in \cF_{d}} S_{\alpha} z^{\alpha}$ where the
    coefficients $S_{\alpha}$ are operators in $\cL(\cU, \cY)$
    for each $\alpha \in \cF_{d}$.
    Given a formal power series
    $S(z) = \sum_{\alpha \in  \cF_{d}} S_{\alpha} z^{\alpha}$
    together with a $d$-tuple of linear operators $\bdelta = (\bdelta_{1}, \dots,
    \bdelta_{d})$ acting on $\ell^{2}$, we define
    $S(\bdelta)$ by
    $$
      S(\bdelta) = \lim_{N \to \infty} \sum_{\alpha \in \cF_{d}
      \colon |\alpha| = N} S_{\alpha} \otimes
      \bdelta^{\alpha}
      \in \cL(\cU \otimes \cK, \cY \otimes \cK)
    $$
    whenever the limit exists in the operator-norm topology;
    here we use the notation $\bdelta^\al$ for the operator
    $$
      \bdelta^{\alpha} = \bdelta_{i_{N}} \cdots \bdelta_{i_{1}} \text{ if } \alpha =
      i_{N} \cdots i_{1} \in \cF_{d}-\{\emptyset\}\text{ and }
      \bdelta^\emptyset=I_\cK.
    $$
    We define the {\em  noncommutative Schur-Agler class}
    $\mathcal{ S A}_{nc,d}(\cU, \cY)$ ({\em strict noncommutative
    Schur-Agler class $\mathcal{S A}^{o}_{nc,d}(\cU, \cY)$}) to consist of all formal power
    series in $\cL(\cU, \cY)\langle \langle z \rangle \rangle$ such that
   $\| S(\bdelta))\| \le 1$ ($\|S(\bdelta)\| < 1$)  whenever $\bdelta
   = (\bdelta_{1}, \dots, \bdelta_{d})$ is a
   $d$-tuple of operators on $\cK$ with  $\|\bdelta_{j}\| < 1$
   ($\|\bdelta_{j} \| \le 1$) for $j = 1, \dots, d$. Let
\begin{align}
 \BD_{nc,d}:=\{\bdelta=(\bdelta_1,\ldots, \bdelta_d) \colon
 \bdelta_j\in\cL(\cK),\ \|\bdelta_j\|< 1,\ j=1,\ldots,d\},\notag\\
 \ov{\BD}_{nc,d}:=\{\bdelta=(\bdelta_1,\ldots, \bdelta_d) \colon
 \bdelta_j\in\cL(\cK),\ \|\bdelta_j\|\leq 1,\ j=1,\ldots,d\}.\notag
\end{align}
   We then define the {\em strict noncommutative $H^{\infty}$-space}
   $H^{\infty,o}(\cL(\cU,\cY))$ to consist of all
   functions $F$ from $\overline{\mathbb D}_{nc,d}$ to $\cL(\cU \otimes
   \cK, \cY \otimes \cK)$ which can be expressed in the form
  $$
  F(\bdelta) = S(\bdelta)
  $$
  for all $\bdelta \in \overline{\mathbb D}_{nc,d}$
  where $\rho^{-1} S$ is in the strict noncommutative Schur-Agler class
  $\mathcal{S A}_{nc,d}^o(\cU, \cY)$ for some real number $\rho > 0$.
  We write $H^{\infty}_{nc, d}(\cL(\cU, \cY))$ for the set of functions $G$
  from $\BD_{nc,d}$ to $\cL(\cU \otimes \cK, \cY \otimes \cK)$ that are also
  of the form $G(\bdelta) = S(\bdelta)$, but now for $\bdelta\in\BD_{nc,d}$ and $\rho^{-1} S$
  in $\mathcal{S A}_{nc,d}(\cU, \cY)$ for some $\rho>0$. Note that
  $\mathcal{ S A}_{nc,d}(\cU, \cY)$ amounts to
  $\mathcal{S A}_{nc,d}({\mathbb C}, {\mathbb C}) \otimes \cL(\cU, \cY)$.
In the sequel we abbreviate the notation $\mathcal{S A}_{nc, d}({\mathbb
  C}, {\mathbb C})$ for the scalar Schur-Agler class to simply
  $\mathcal{S A}_{nc, d}$. Similarly, we simply write $\mathcal{SA}_{nc,d}^o$,
  $H^{\infty,o}_{nc,d}$ and $H^\infty_{nc,d}$ instead of $\mathcal{SA}_{nc,d}^o(\BC,\BC)$,
  $H^{\infty,o}_{nc,d}(\BC,\BC)$ and $H^\infty_{nc,d}(\BC,\BC)$, respectively.
    Thus we also have $H^{\infty,o}_{nc,
  d}({\mathcal L}(\cU, \cY)) = H^{\infty,o}_{nc,d} \otimes \cL(\cU,
  \cY)$, etc.  We shall be primarily interested in the strict
  versions $\mathcal{ S A}^{o}_{nc, d}$ and $H^{\infty,o}_{nc, d}$ of
  the noncommutative Schur-Agler class and $H^{\infty}$-space.

  We have the following characterization of the space
  $H^{\infty,o}_{nc, d}(\cL(\cU, \cY))$.
  For the definition of {\em completely positive kernel} and more
  complete details, we refer to \cite{BtH}.  The formulation given
  here does not have the same form as in Theorem 3.6(2) of
  \cite{BtH}, but one can use the techniques given there to convert
  to the form given in the following theorem.

  \begin{theorem}   \label{T:cp-diskalg}
    The function $F \colon \overline{\mathbb D}_{nc, d} \to \cL(\cU
    \otimes \cK, \cY \otimes \cK)$ is in the strict noncommutative
    $H^{\infty}$-space $H^{\infty,o}_{nc, d}(\cL(\cU,\cY))$ if and only
    if there are $d$ strictly completely positive kernels
    $$
      K_{k} \colon (\overline{{\mathbb D}}_{nc,d} \times
      \overline{\mathbb D}_{nc, d}) \times \cL(\cK)  \to \cL(Y \otimes
      \cK) \text{ for } k = 1, \dots, d
      $$ and a positive real number $\rho$ so that
      the following Agler decomposition holds:
      \begin{align*}
      &  \rho^{2} \cdot (I \otimes B) - S(\bdelta) \left( I \otimes B
      \right) S(\btau)^{*} =
      \sum_{k=1}^{d}  K_{k}(\bdelta,\btau)[B - \bdelta_{k}B\btau_{k}^{*}]
     \end{align*}
     for all $B \in \cL(\cK)$ and  $\bdelta= (\bdelta_{1}, \dots, \bdelta_{d})$,
     $\btau =(\btau_{1}, \dots, \btau_{d})$ in $\overline{\mathbb D}_{nc, d}$.
     \end{theorem}

     One of the main results of \cite{BGM2} is that the noncommutative
     Schur-Agler class has a contractive Givone-Roesser realization.

     \begin{theorem}  \label{T:BGM2} (See \cite{BGM2, BGM3}.)  A given
     function $F \colon \overline{\mathbb D}_{nc, d} \to \cL(\cU
     \otimes \cK, \cY \otimes \cK)$
     is in the strict noncommutative Schur-Agler class $\mathcal{S
     A}_{nc, d}^o(\cU, \cY)$ if and only if
     there exists a strictly contractive colligation matrix
     $$
     M = \begin{bmatrix} A & B \\ C & D \end{bmatrix}
     \colon \begin{bmatrix} \oplus_{j=1}^{d} \cX_{j} \\ \cU \end{bmatrix} \to
     \begin{bmatrix} \oplus_{j=1}^{d} \cX_{j} \\ \cY \end{bmatrix}
    $$
  for some Hilbert state space $\cX=\cX_1\oplus\cdots\oplus\cX_d$ so that the
  evaluation of $F$ at
  $\bdelta  = (\bdelta_{1}, \dots, \bdelta_{d}) \in \overline{\mathbb D}_{nc,d}$
  is given by
  \begin{equation}   \label{ncpolydisk-real}
   F(\bdelta) = D \otimes I_{\cK} + ( C \otimes I_{\cK}(
  ( I -  \bZ(\bdelta) (A \otimes I_{\cK}))^{-1} \bZ(\bdelta) (B \otimes I_{\cK})
  \end{equation}
  where
  $$
    \bZ(\bdelta) =  \begin{bmatrix}I_{\cX_1} \otimes \bdelta_{1} & & \\ &
    \ddots & \\ & & I_{\cX_d} \otimes\bdelta_{d}
  \end{bmatrix}.
  $$
  Hence a function $F \colon \overline{\mathbb D}_{nc, d} \to \cL(\cU
  \otimes \cK, \cY \otimes \cK)$ is in the strict noncommutative
  $H^{\infty}$-space
  $H^{\infty,o}_{nc, d}(\cL(\cU, \cY))$ if and only if there is a bounded
  linear operator
  $$
   \begin{bmatrix} A & B \\ C & D \end{bmatrix} \colon \begin{bmatrix}
       \oplus_{k=1}^{d} \cX_{k} \\ \cU \end{bmatrix} \to
       \begin{bmatrix} \oplus_{k=1}^{d} \cX_{k} \\ \cY \end{bmatrix}
  $$
  such that
  $$ \left\| \begin{bmatrix} A & B \\ \rho^{-1} C & \rho^{-1} D
  \end{bmatrix} \right\| < 1 \text{ for some } \rho > 0
  $$
  so that $F$ is given as in \eqref{ncpolydisk-real}.
  \end{theorem}

If $\cU$ and $\cY$ are finite-dimensional Hilbert spaces, we may
view $\mathcal{S A}^{o}_{nc, d}(\cU, \cY)$ and $H^{\infty,o}_{nc,
d}(\cL(\cU, \cY))$ as matrices
over the respective scalar-valued classes $\mathcal{ S A}^{o}_{nc,d}$ and
$H^{\infty,o}_{nc, d}$.  When this is the case,
it is natural to define {\em rational} versions of $\mathcal{S A}_{nc,d}^o$ and
$H^{\infty,o}_{nc, d}$ to consist of those functions in $\mathcal{S
A}_{nc, d}^{o}$ (respectively, $H^{\infty,o}_{nc, d}$) for which the
realization \eqref{ncpolydisk-real} can be
taken with the state spaces $\cX_{1}, \dots, \cX_{d}$ also
finite-dimensional; we denote the rational versions of $\mathcal{S
A}_{nc, d}^o$ and $H^{\infty,o}_{nc, d}$ by $\cR \mathcal{S A}_{nc, d}^o$
and $\cR H^{\infty,o}_{nc, d}$, respectively.  We remark that as a
consequence of Theorem 11.1 in \cite{BGM1}, this rationality
assumption on a given function $F$ in $H^{\infty,o}_{nc, d}$
can be expressed intrinsically in terms of the finiteness of rank for
a finite collection of Hankel matrices  formed from the  power-series
coefficients $F_{\alpha}$ of $F$, i.e., the operators $F_{\alpha} \in
\cL(\cU, \cY)$ such that
$$
  F(\bdelta) = \sum_{\alpha \in \cF_{d}} F_{\alpha} \otimes
  \bdelta^{\alpha}.
$$

In general, the embedding of a noncommutative integral domain into a
skew-field is difficult (see e.g.~\cite{HMcCV, KVV09}).  For the case of
$\cR H^{\infty,o}$, the embedding issue becomes tractable if
we restrict to denominator functions $D(\bdelta) \in
H^{\infty,o} \in \cL(\cU)$ for which $D(0)$ is invertible.
If $D$ is given in terms of a strictly contractive
realization $D(\bdelta) =  \bD + \bC(I
- \bZ(\bdelta) \bA)^{-1} \bZ(\bdelta) \bB$ (where $\bA = A \otimes
I_{\cK}$ and similarly for $\bB$, $\bC$ and $\bD$),
  then $D(\bdelta)^{-1}$ can be calculated, at least for $\|\bZ(\bdelta)\|$
  small enough, via
the familiar cross-realization formula for the inverse:
$$
 D(\bdelta)^{-1} = \bD^{-1} - \bD^{-1} \bC (I - \bZ(\bdelta)
 \bA^{\times})^{-1} \bZ(\bdelta) \bB \bD^{-1}
 $$
where $\bA^{\times} = A^{\times} \otimes I_{\cK}$ with $A^{\times} =
 A - B D^{-1} C$.
We define $Q(\cR H^{\infty,o}_{nc, d})(\cL(\cU, \cY))_{0}$ to be the
smallest linear space of functions from some neighborhood of $0$ in
$\overline{\mathbb D}_{nc,d}$ (with respect to the Cartesian product
operator-norm topology on $\overline{{\mathbb D}}_{nc,d} \subset
{\mathcal L}({\mathcal K})^{d}$)
to $\cL(\cU, \cY)$ which is invariant under multiplication on the left by
elements of $\cR H^{\infty,o}_{nc, d}(\cL(\cY))$ and by inverses of elements
of $\cR H^{\infty,o}_{nc, d}(\cL(\cY))$ having invertible value at $0$,
and invariant under multiplication on the right by the corresponding
set of functions with $\cU$ in place of $\cY$. Note that the final
subscript $0$ in the notation
$Q(\cR H^{\infty,o}_{nc, d})(\cL(\cU, \cY))_{0}$
is suggestive of the requirement that functions of this class
are required to be
analytic in a neighborhood of the origin $0 \in {\mathbb D}_{nc, d}$.

Let us denote by $\cR \cO^{0}_{nc, d}(\cL(\cU, \cY))$ the space
of functions defined as follows: we say that the
function $G$ defined on a neighborhood of the origin in ${\mathbb
D}_{nc, d}$ with values in $\cL(\cU, \cY)$ is in the space $\cR
\cO^{0}_{nc, d}(\cL(\cU, \cY))$ if $G$ has a realization of the form
$$
  G(\bdelta) = \bD + \bC (I - \bZ(\bdelta) \bA)^{-1} \bZ(\bdelta) \bB
$$
for a colligation matrix $\bM: = \sbm{ \bA & \bB \\ \bC & \bD }$
of the form  $\bM = M \otimes I_{\cK}$ where
$$
 M = \begin{bmatrix} A & B \\ C  & D \end{bmatrix} \colon
 \begin{bmatrix} \oplus_{k=1}^{d} \cX_{k} \\ \cU \end{bmatrix} \to
     \begin{bmatrix} \oplus_{k=1}^{d} \cX_{k} \\ \cY \end{bmatrix}
$$
for some finite-dimensional state-spaces $\cX_{1}, \dots, \cX_{d}$.
Unlike the assumptions in the case of a realization for a
Schur-Agler-class function in Theorem \ref{T:BGM2}, there is no
assumption that $M$ be contractive or that $A$ be stable.
It is easily seen that $Q(\cR H^{\infty,o}_{nc, d}(\cL(\cU, \cY)))_0$ is a
subset of $\cR \cO^{0}_{nc, d}(\cL(\cU, \cY))$; whether these two
spaces are the same or not we leave as an open question.
We also note that the class $\cR \cO^{0}_{nc, d}(\cL(\cU, \cY))$ has
an intrinsic characterization:  $F$ is in $\cR \cO^{0}_{nc, d}(\cL(\cU, \cY))$
if and only if some rescaled version $ \widetilde F(\bdelta) = F(r
\bdelta)$ (where $r \bdelta = (r \bdelta_{1}, \dots, r \bdelta_{d})$
if $\bdelta = (\bdelta_{1}, \dots, \bdelta_{d})$) is in the
rational noncommutative $H^{\infty}$-class $ \cR H^{\infty,o}_{nc, d}(\cL(\cU,
\cY))$ for some $r > 0$ and hence  has the intrinsic characterization
in terms of a completely positive Agler decomposition and
finite-rankness of a finite collection of Hankel matrices as
described above for the class $\cR H^{\infty,o}_{nc, c}(\cL(\cU,
\cY))$.

We may then pose the following
control problems:

\smallskip
   \noindent
  \textbf{Noncommutative polydisk
       internal-stabilization/$H^{\infty}$-control problem:}
       We suppose that we are given finite-dimensional
       spaces $\cW$, $\cU$, $\cZ$, $\cY$ and a block-matrix $G =
       \sbm{G_{11} & G_{12}  \\ G_{21} & G_{22}}$  in
       $ \cR \cO^{0}_{nc, d}(\cL( \cW \oplus \cU, \cZ \oplus \cY))$.
     We seek to find a controller $K$ in $\cR \cO^{0}_{nc, d}(\cL(\cY,
     \cU))$
     which solves the (1) {\em internal stabilization problem}, i.e.
     so that the closed-loop system is {\em internally
     stable} in the sense that
       all matrix entries of the block matrix $\Theta(G,K)$ given by
       \eqref{ThetaGK} are in $\cR H^{\infty,o}_{nc, d}$,
       and which
       possibly also solves the (2) {\em $H^{\infty}$-problem}, i.e., in
       addition to internal stability, the closed-loop system has
       {\em performance} in the sense that
       $ T_{zw} = G_{11} + G_{12} K (I - G_{22} K)^{-1} G_{21}$ is in
       the rational strict noncommutative Schur-Agler class $\cR \mathcal{S
       A}^{o}_{nc, d}(\cW, \cZ)$.

    \smallskip

    Even though our algebra of scalar plants $\cR \cO^{0}_{nc, d}$ is
    noncommutative, the parameterization result Theorem \ref{T:stab3}
    still goes through in the following form; we leave it to the reader to
    check that the same algebra as used for the commutative case
    leads to the following noncommutative analogue.

    \begin{theorem} \label{T:NCstab3}  Assume that $G \in \cR
    \cO^{0}_{nc, d}(\cL(\cW \oplus \cU, \cZ \oplus \cY))$ is given
    and that $G$ has at least one stabilizing controller $K_*$.
     Define $U_{*} =
       (I - G_{22} K_{*})^{-1}$, $V_{*} = K_{*}(I - G_{22} K_{*})^{-1}$,
       $\widetilde U_{*} = (I - K_{*} G_{22})^{-1}$ and $\widetilde V_{*} = (I -
       K_{*} G_{22})^{-1} K_{*}$.
       Then the set of all stabilizing controllers $K$ for $G$ is given by
       either of the two formulas
       \begin{align*}
       & K = (V_{*} + Q) (U_{*} + G_{22} Q)^{-1} \text{ subject to }
       (U_{*} + G_{22} Q)(0) \text{ is invertible,} \\
       &K=(\widetilde U_{*} + Q G_{22})^{-1} (\widetilde V_{*} + Q)
       \text{ subject to } (\widetilde U_{*} + Q G_{22})(0) \text{ is
       invertible},
     \end{align*}
     where in addition $Q$ has the form $Q = \widetilde L \Lambda L$
     where $\widetilde L$ and $L$ are
     given by \eqref{tildeLL} and $\Lambda$ is a free stable parameter
     in $H^{\infty,o}_{nc, d}(\cL(\cY \oplus \cU, \cU \oplus
     \cY))$. Moreover, if $Q = \widetilde L \Lambda L$ with $\La$ stable,
     then $(U_{*} + G_{22} Q)(0)$ is invertible if and only if
     $(\widetilde U_{*} + Q G_{22})(0)$ is invertible, and
     both formulas give rise to the same controller $K$.
     \end{theorem}

     Given a transfer matrix $G_{22} \in \cR \cO^{0}_{nc, d}(\cL(\cU, \cY))$,
     we say that $G_{22}$ has a {\em stable double coprime factorization}
     if there exist transfer matrices $D(\bdelta)$, $N(\bdelta)$, $X(\bdelta)$,
     $Y(\bdelta)$, $\widetilde D(\bdelta)$, $\widetilde N(\bdelta)$,
     $\widetilde X(\bdelta)$, and $\widetilde Y(\bdelta)$ of
     compatible sizes with stable matrix entries (i.e., with matrix
     entries in $\cR H^{\infty,o}_{nc, d}$) subject also to
     $$ D(0), \, \widetilde D(0), \, X(0), \, \widetilde X(0) \text{
     all invertible }
     $$
     so that the noncommutative version of condition
     \eqref{coprime} holds:
    \begin{equation}\label{NCcoprime}
    \begin{array}{c}
    G_{22}(\bdelta) = D(\bdelta)^{-1} N(\bdelta) = \widetilde N(\bdelta) \widetilde D^{-1}(\bdelta),
    \\[.2cm]
     \begin{bmatrix} D(\bdelta) & -N(\bdelta) \\ -\widetilde Y(\bdelta) & \widetilde X(\bdelta)
     \end{bmatrix} \begin{bmatrix} X(\bdelta) & \widetilde N(\bdelta) \\ Y(\bdelta) &
     \widetilde D(\bdelta) \end{bmatrix} =\mat{cc}{I_{n_{\cY}}&0\\0&I_{n_{\cU}}}.
     \end{array}
    \end{equation}
    Then we leave it to the reader to check that the same algebra as
    used for the commutative case leads to the following
    noncommutative version of Theorem \ref{T:stab4}.

    \begin{theorem}   \label{T:NCstab4}
    Assume that $G\in \cR \cO^{0}_{nc, d}$ is stabilizable and that $G_{22}$
    admits a double coprime factorization \eqref{NCcoprime}. Then the set of all
    stabilizing controllers is given by
\begin{eqnarray*}
    K(\bdelta)&=&(Y(\bdelta)+\wtilD(\bdelta)\La(\bdelta))(X(\bdelta)+\wtilN(\bdelta)\La(\bdelta))^{-1}\\
    &=&(\wtilX(\bdelta)+\La(\bdelta) N(\bdelta))^{-1}(\wtilY(\bdelta)+\La(\bdelta) D(\bdelta)),
\end{eqnarray*}
    where $\La$ is a free stable parameter from $H^{\infty,0}_{nc,
    d}(\cL(\cU, \cY)$ such that
    $X(0)-\wtilN(0)\La(0)$ is invertible and $\wtilX(0)+\La(0) N(0)$ is
    invertible.
   \end{theorem}

   Just as in the commutative case, consideration of the
   $H^\infty$-control problem for a given transfer matrix $G \in \cR
   \cO^{0}_{nc, d}(\cL(\cW \oplus \cU, \cZ \oplus \cY))$ after the
   change of the design parameter from the controller
   $K$ to the free-stable parameter $\Lambda$ in either of the two
   parameterizations of Theorems \ref{T:NCstab3} and \ref{T:NCstab4}
   leads to the following noncommutative version of the
   Model-Matching problem; we view this problem as a noncommutative
   version of a Sarason interpolation problem.

     \smallskip
     \noindent
     \textbf{Noncommutative-polydisk Sarason interpolation problem:}
     Given matrices $T_{1}$, $T_{2}$, $T_{3}$ of compatible sizes over
     $\cR H^{\infty,o}_{nc, d}$,
      find a matrix $\Lambda$ (of
      appropriate size) over $\cR H^{\infty,o}_{nc, d}$
      so that the matrix $S = T_{1} + T_{2} \Lambda T_{3}$
      is in the strict rational noncommutative Schur-Agler class $\cR \mathcal { S A}_{nc,
      d}^o(\cW, \cZ)$.

      While there has been some work on left-tangential
      Nevanlinna-Pick-type interpolation for the noncommutative
      Schur-Agler class (see \cite{BallBoloJOT}), there does not seem to
      have been any work on  a Commutant Lifting theorem for this setup
      or on how to convert a Sarason problem as above to an
      interpolation problem as formulated in \cite{BallBoloJOT}.  We
      leave this area to future work.

   \subsection{Equivalence of state-space noncommutative LFT-model
   and noncommutative frequency-domain formulation}
   \label{S:NCstate-freq}

In order to make the connections between the results in the previous two subsections,
we consider functions as in Subsection \ref{S:NC-FD}, but we normalize the infinite
dimensional Hilbert space $\cK$ to be $\ell^2$ and work with $d+1$ variables
$\bdelta=(\bdelta_1,\ldots,\bdelta_{d+1})$ in $\cL(\ell^2)$ instead of $d$. As pointed
out in Subsection \ref{S:NC-statespace}, we may without loss of generality assume that
the last variable $\bdelta_{d+1}$ is fixed to be the shift operator $\bS$ on $\ell^2$.

   The following is an improved analogue of Lemma \ref{L:ND-sing} for
   the noncommutative setting.

   \begin{theorem}   \label{T:noncom-sing}  Suppose that the matrix
       function $W \in \cR \cO^{0}_{nc, d+1}(\cL(\cU, \cY))$ has a
       finite-dimensional realization
       $$
  W(\bdelta) = \bD + \bC (I - \bZ(\bdelta) \bA)^{-1} \bZ(\bdelta)\bB,
  $$
 where
 $$ \bA = A \otimes I_{\ell^{2}}, \quad \bB = B \otimes I_{\ell^{2}},
         \quad \bC = C \otimes I_{\ell^{2}}, \quad \bD = D
         \otimes I_{\ell^{2}},
$$
 which is both nc-Hautus-detectable and nc-Hautus-observable.
 Then $W$ is stable in the noncommutative frequency-domain sense
 (i.e., all matrix entries of $W$ are in $H^{\infty,o}_{nc, d+1}$) if
 and only if $W$ is stable in the state-space sense, i.e., the matrix
 $\bA$ is nc-Hautus-stable.
    \end{theorem}

       \begin{proof}
  If the matrix $\bA$ is nc-Hautus-stable, it is trivial that then all
  matrix entries of $W$ are in $H^{\infty,o}_{nc, d+1}$.  We therefore
  assume that all matrix entries of $W$ are in $H^{\infty,o}_{nc, d+1}$.
  It remains to show that, under the assumption that $\{C, A\}$ is
  nc-Hautus detectable and that $\{A, B\}$ is nc-Hautus stabilizable,
  it follows that $A$ is nc-Hautus stable.

  The first step is to observe the identity
  \begin{equation} \label{NCidentity1}
    S_{1}(\bdelta): =   \begin{bmatrix} I -  \bZ(\bdelta) \bA \\ \bC \end{bmatrix} (I -
      \bZ(\bdelta) \bA)^{-1} \bZ(\bdelta) \bB = \begin{bmatrix}
      \bZ(\bdelta) \bB \\ W(\bdelta) - \bD \end{bmatrix}.
  \end{equation}
  Since $W(\bdelta) - \bD$ is in $H^{\infty,o}_{nc,d+1}(\cL(\cU, \cY))$
  by assumption and trivially $\bZ(\bdelta) \bB$ is in
  $H^{\infty,o}_{nc, d+1}(\cL(\cU, \cX))$, it follows that
  $S_{1}(\bdelta)$ is in $H^{\infty,o}_{nc, d+1}(\cL(\cU, \cX \oplus
  \cY))$.  By the detectability assumption and Proposition \ref{P:NCdetect}
  it follow that there exists an operator $\bL=L\otimes I_{\ell^2}$
  with $L:\cY\to\cX$ so that $\bA+\bL\bC$ is nc-Hautus-stable. Thus
  \[
F_1(\bdelta)=(I-\bZ(\bdelta)(\bA+\bL\bC))^{-1}\mat{cc}{I&-\bZ(\bdelta)\bL}
  \]
is in $H^{\infty,o}_{nc, d+1}(\cL(\cX \oplus\cY,\cX))$. Note that
$F_1(\bdelta)S_1(\bdelta)=(I- \bZ(\bdelta) \bA)^{-1} \bZ(\bdelta) \bB$.
The fact that both $F_1$ and $S_1$ are transfer-functions over $H^{\infty,o}_{nc, d+1}$
implies that $S_2(\bdelta)=(I- \bZ(\bdelta) \bA)^{-1} \bZ(\bdelta) \bB$ is in
$H^{\infty,o}_{nc, d+1}(\cL(\cU,\cX))$.

We next use the identity
   \begin{align}
        \begin{bmatrix} \bZ(\bdelta) & S_{2}(\bdelta) \end{bmatrix}:
      =  &
      \begin{bmatrix} \bZ(\bdelta) & (I - \bZ(\bdelta) \bA)^{-1}
       \bZ(\bdelta) \bB \end{bmatrix} \notag \\
       =  & \bZ(\bdelta) (I - \bA
       \bZ(\bdelta))^{-1} \begin{bmatrix} I - \bA \bZ(\bdelta) & \bB
       \end{bmatrix}.
       \label{NCidentity2}
  \end{align}
Now the nc-Hautus-stabilizability assumption and the second part
of Proposition \ref{P:NCdetect} imply in a similar way that
$S_3(\bdelta)=\bZ(\bdelta)(I-\bZ(\bdelta)\bA)^{-1}$ is in
$H^{\infty,o}_{nc, d+1}(\cL(\cX,\cX))$.
%
 %  by the detectability assumption, $\sbm{ I - \bZ(\bdelta) \bA
%  \\ \bC }$ has a left inverse $L(\bdelta)$ for each $\bdelta \in
%  \overline{\mathbb D}_{nc, d}$.  Moreover, by Theorem B.3 in
%  \cite{Paganini}, it follows that we may assume that
%  $$
%    \sup \{ \| L(\bdelta) \| \colon \bdelta \in \overline{\mathbb
%    D}_{nc, d} \} < \infty.
%   $$
%   From the fact that $S_{1}$ in \eqref{NCidentity1} is in
%   $H^{\infty,0}_{nc, d}(\cL(\cU, \cX \oplus \cY))$, we may now
%   conclude that $S_{2}(\bdelta): = (I - \bZ(\bdelta) \bA)^{-1}
%   \bZ(\bdelta) \bB$ is in $H^{\infty,o}_{nc, d}(\cL(\cU, \cX))$.
%
%   We
%   next use the identity
%   \begin{align}
%        \begin{bmatrix} \bZ(\bdelta) & S_{2}(\bdelta) \end{bmatrix}:
%      =  &
%      \begin{bmatrix} \bZ(\bdelta) & (I - \bZ(\bdelta) A)^{-1}
%       \bZ(\bdelta) B \end{bmatrix} \notag \\
%       =  & \bZ(\bdelta) (I - \bA
%       \bZ(\bdelta))^{-1} \begin{bmatrix} I - \bA \bZ(\bdelta) & B
%       \end{bmatrix}.
%       \label{NCidentity2}
%  \end{align}
%  By the stabilizability assumption combined with the dual version of
%  Theorem B.3 from \cite{Paganini}, the fact that $S_{2}(\bdelta)$ is in
%  $H^{\infty,o}_{nc, d}(\cL(\cU,
%  \cX))$ implies that $S_{3}(\bdelta): = \bZ(\bdelta) (I - \bA
%  \bZ(\bdelta))^{-1}$ is in $H^{\infty,o}_{nc, d}(\cL(\cX))$.
  Note
  that $S_{3}$ in turn has the trivial realization
  $$
    S_{3}(\bdelta) = \bD' + \bC' (I - \bZ(\bdelta) \bA')^{-1}
    \bZ(\bdelta) \bB'
  $$
  where $\sbm{\bA' & \bB' \\ \bC' & \bD' } = \sbm{ A' & B' \\ C' & D'
  } \otimes I_{\ell^{2}}$ and $\sbm{ A' & B' \\ C' & D' } = \sbm{ A &
  I \\ I & 0 }$.  Thus $(A', B', C', D') = (A, I, I, 0)$ is trivially
  GR-controllable and GR-observable in the sense of \cite{BGM1}.  On
  the other hand, by Theorem \ref{T:BGM2} there exists a strictly
  contractive matrix $\sbm{A'' & B'' \\ C'' & 0 }$ so that
  $$
    S_{3}(\bdelta) = r'' \bC''(I - \bZ(\bdelta) \bA'')^{-1} \bZ(\bdelta)
    \bB''
  $$
  for some $r< \infty$.  Moreover, by the Kalman decomposition for
  noncommutative GR-systems given in \cite{BGM1}, we may assume
  without loss of generality that $(A'', B'', C'', 0)$ is
  GR-controllable and GR-observable.  Then, by the main result of
  Alpay--Kaliuzhnyi-Verbovetskyi in \cite{A-KV}, it is known that the
  function $S(\bdelta) = \sum_{\alpha \in \cF_{d}} S_{\alpha} \otimes
  \bdelta^{\alpha}$ uniquely determines the formal power series $S(z)
  = \sum_{\alpha \in \cF_{d}} S_{\alpha} z^{\alpha}$.
  It now follows from the
  State-Space Similarity Theorem for noncommutative GR-systems in \cite{BGM1} that
  there is an invertible block diagonal similarity transform $Q  \in
  \cL(\cX', \cX'')$ so that
  $$
  \begin{bmatrix} A & I \\ I & 0 \end{bmatrix} : = \begin{bmatrix} A'
      & B' \\ C' & 0 \end{bmatrix} = \begin{bmatrix} Q^{-1} & 0 \\ 0
      & I \end{bmatrix} \begin{bmatrix} A'' & B'' \\ C'' & 0
  \end{bmatrix} \begin{bmatrix} Q & 0 \\ 0 & I \end{bmatrix}.
   $$
   In particular, $A = Q^{-1} A'' Q$ where $A''$ is a strict contraction
   and $Q$ is a structured similarity from which it
   follows that $A$ is also nc-Hautus-stable as wanted.
\end{proof}

We can now obtain the equivalence of the frequency-domain and
state-space formulations of the internal stabilization problems for
the case where the state-space internal stabilization problem is
solvable.

\begin{theorem} \label{T:NCfreq-state}
   Suppose that we are given a realization
  $$
  G(\bdelta) = \begin{bmatrix} G_{11}(\bdelta) & G_{12}(\bdelta) \\
  G_{21}(\bdelta) & G_{22}(\bdelta) \end{bmatrix} =
  \begin{bmatrix} \bD_{11} & \bD_{12} \\ \bD_{21} & 0 \end{bmatrix}
      + \begin{bmatrix} \bC_{1} \\ \bC_{2} \end{bmatrix}
      ( I - \bZ(\bdelta) \bA)^{-1} \bZ(\bdelta) \begin{bmatrix} \bB_{1} &
      \bB_{2} \end{bmatrix}
  $$
  for an element $G \in \cR \cO^{0}_{nc, d+1}(\cL(\cW \oplus \cU, \cZ
  \oplus \cY))$ such that the state-space internal stabilization
  problem has a solution.  Suppose also that we are given a controller
  $K \in \cR \cO^{0}_{nc, d+1}(\cL(\cY, \cU))$ with state-space
  realization
  $$
   K(\bdelta) = \bD_{K} + \bC_{K}(I - \bZ_{K}(\bdelta) \bA_{K})^{-1}
   \bZ_{K}(\bdelta) \bB_{K}.
  $$
  which is both nc-Hautus-stabilizable and nc-Hautus-detectable.
  Then the controller $K\sim\{\bA_{K}, \bB_{K}, \bC_{K}, \bD_{K}\}$
  solves the state-space internal
  stabilization problem associated with $\{\bA, \begin{bmatrix} \bB_{1} & \bB_{2} \end{bmatrix},
  \sbm{ \bC_{1} \\ \bC_{2} }, \sbm{ \bD_{11} & \bD_{12} \\ \bD_{21} & 0}\}$
  if and only if $K(\bdelta)$ solves the noncommutative
  frequency-domain internal stabilization problem associated with
  \[
  G(\bdelta) = \sbm{ G_{11}(\bdelta) & G_{12}(\bdelta) \\
  G_{21}(\bdelta) & G_{22}(\bdelta) }.
  \]
  \end{theorem}

  \begin{proof}  By Theorem \ref{T:NCstable2}, the assumption that that
  the state-space internal stabilization problem is solvable means that
  $\{\bC_{2},\bA\}$ is nc-Hautus-detectable and $\{\bA,\bB_{2}\}$ is
  nc-Hautus-stabilizable. We shall use this form of the standing assumption.
  Moreover, in this case, a given controller $K \sim \{\bA_{K}, \bB_{K}, \bC_{K}, \bD_{K}\}$
  solves the state-space internal stabilization problem if and only if
  $K$ stabilizes $G_{22}$.

      Suppose now that $K \sim \{\bA_{K}, \bB_{K}, \bC_{K},\bD_{K}\}$ solves
      the state-space internal stabilization problem, i.e., the state operator
      $\bA_{cl}$ in \eqref{clsysmat-NC} is nc-Hautus-stable. Note that the
      $3 \times 3 $ noncommutative transfer matrix $\boldsymbol \Theta(G, K)$ has realization
      $\boldsymbol \Theta(G, K) = \bD_{\Th} + \bC_{\Th}(I - \bZ_{\Th}(\bdelta)
      \bA_{\Th})^{-1} \bZ_{\Th}(\bdelta) \bB_{\Th}$ with $\bZ_{\Th}(\bdelta)=\bZ_{cl}(\bdelta)$
      as in \eqref{Zcl} where
      \[
      \mat{cc}{ \bA_{\Th} & \bB_{\Th} \\ \bC_{\Th} & \bD_{\Th}} =
      \mat{cc}{A_{\Th} & B_{\Th} \\ C_{\Th} & D_{\Th}} \otimes I_{\ell^{2}}
      \]
      with
      \begin{align}
     &  A_{\Th} = \begin{bmatrix} A + B_{2} D_{K}
          C_{2} & B_{w}C_{K} \\ B_{K} C_{2} & A_{K}
      \end{bmatrix}, \quad
      B_{\Th} = \begin{bmatrix} B_{1} + B_{2} D_{K} D_{2} & B_{2}
      & B_{2} D_{K} \\ B_{K} D_{21} & 0 & B_{K} \end{bmatrix},
      \notag \\
      & C_{\Th} = \begin{bmatrix} C_{1} + D_{12} D_{K} C_{2} & D_{12}
      C_{K} \\ D_{K} C_{2} & C_{K} \\ C_{2} & 0 \end{bmatrix}, \quad
     D_{\Th} =  \begin{bmatrix} D_{1} + D_{12} D_{K} D_{21} & D_{12} & D_{12}
      D_{K} \\ D_{K} D_{21} & I & K_{K} \\ D_{21} & 0 & I
      \end{bmatrix}.
      \label{realThetaGK}
      \end{align}
      Now observe that $\bA_{\Th}$ is equal to $\bA_{cl}$, so that all nine
      transfer matrices in $\boldsymbol \Theta(G, K)$ have a realization with
      state operator $\bA_{\Th}=\bA_{cl}$ nc-Hautus-stable. Hence all matrix
      entries of $\boldsymbol \Theta(G, K)$ are in $H^{\infty,o}_{nc, d+1}$.

       Suppose that $K(\bdelta)$ with realization $K \sim
       \{\bA_{K}, \bB_{K}, \bC_{K}, \bD_{K}\}$ internally stabilizes $G$ in
       the frequency-domain sense.  This means that all nine transfer
       matrices in $\Theta(G, K)$ are stable.  In particular, the $2
       \times 2$ transfer matrix
       $\widetilde W : = \boldsymbol \Theta(G_{22}, K) - \boldsymbol \Theta(G_{22}, K)(0)$ is
       stable.  From  \eqref{realThetaGK} we read off that
      $ \widetilde W$ has realization
      $$
       \widetilde W(\bdelta) = \begin{bmatrix} \bD_{K}\bC_{2} &
       \bC_{K} \\ \bC_{2} & 0 \end{bmatrix}
       (I - \bZ_{\Th}(\bdelta) \bA_{\Th})^{-1} \begin{bmatrix} \bB_{2} &
       \bB_{2} \bD_{K} \\ 0 & \bB_{K} \end{bmatrix}.
     $$
     By Theorem \ref{T:noncom-sing}, to show that $\bA_{cl}=\bA_{\Th}$ is
     nc-Hautus-stable, it suffices to show that
     $\left\{ \sbm{\bD_{K}\bC_{2} & \bC_{K} \\ \bC_{2} & 0 }, \bA_{cl}\right\}$
     is nc-Hautus-detectable and that $\left\{ \bA_{cl},
     \sbm{ \bB_{2} & \bB_{2}\bD_{K} \\ 0 & \bB_{K}} \right\}$ is
     nc-Hautus-stabilizable.  By using our assumption
     that $\{\bA_{K}, \bB_{K}, \bC_{K}, \bD_{K}\}$ is both nc-Hautus-detectable and
     nc-Hautus-stabilizable, one can now follow the
     argument in the proof of Theorem \ref{T:NDfreq-state} to deduce
     that  $\left\{ \sbm{\bD_{K}\bC_{2} & \bC_{K} \\ \bC_{2} & 0 },
     \bA_{cl}\right\}$ is noncommutative detectable and that
     $\left\{ \bA_{cl},
      \sbm{ \bB_{2} & \bB_{2}\bD_{K} \\ 0 & \bB_{K}} \right\}$ is
      noncommutative Hautus-stabilizable as needed.
   \end{proof}

   We do not know as of this writing whether any given controller $K$ in the space
   $\cR \cO^{0}_{nc, d+1}(\cL(\cY, \cU))$ has a nc-Hautus-detectable/stabilizable
   realization (see the discussion in
   the Notes below).  However, for the Model-Matching
   problem, internal stabilizability in the frequency-domain sense
   means that all transfer matrices $T_{1}, T_{2}, T_{3}$ are stable
   (i.e., have all matrix entries in $H^{\infty,o}_{nc, d+1}$) and hence
   the standard plant matrix $G = \sbm{ T_{11} & T_{12} \\ T_{22} & 0
   }$ has a stable realization.  A given controller $K$ solves the
   internal stabilization problem exactly when it is stable; thus we
   may work with realizations $K \sim \{\bA_{K}, \bB_{K}, \bC_{K}, \bD_{K}\}$
   with $\bA_{K}$ nc-Hautus-stable, and hence a fortiori
   with both $\{ \bC_{K}, \bA_{K} \}$ nc-Hautus-detectable
   and $\{\bA_{K},\bB_{K}\}$ nc-Hautus-stabilizable.  In this
   scenario Theorem \ref{T:NCfreq-state} tells us that a controller
   $K(\bdelta)$ solves the frequency-domain internal stabilization
   problem exactly when any stable realization $K \sim \{\bA_{K}, \bB_{K},
   \bC_{K}, \bD_{K}\}$ solves the state-space internal stabilization
   problem.  Moreover, the frequency-domain performance measure
   matches with the state-space performance measure, namely: that the
   closed-loop transfer matrix $T_{zw} = G_{11} + G_{12} (I - K
   G_{22})^{-1} K G_{21}$ be in the strict noncommutative Schur-Agler
   class $\mathcal{ S A}^{o}_{nc, d+1}(\cW, \cZ)$.  In this way we
   arrive at a solution of the noncommutative Sarason interpolation
   problem posed in Section \ref{S:NC-FD}.

   \begin{theorem} \label{T:noncomSarason}
       Suppose that we are given a transfer matrix of the
       form  $G = \sbm{ T_{1} & T_{2} \\ T_{3} & 0} \in
       H^{\infty,o}_{nc, d+1}(\cL(\cW \oplus \cU, \cZ \oplus \cY))$ with
       a realization
       $$
       \begin{bmatrix} T_{1}(\bdelta) & T_{2}(\bdelta) \\
       T_{3}(\bdelta) & 0 \end{bmatrix} =
       \begin{bmatrix} \bD_{11} & \bD_{12} \\ \bD_{21} & 0
       \end{bmatrix} + \begin{bmatrix} \bC_{1} \\ \bC_{2}
       \end{bmatrix} (I - \bZ(\bdelta) \bA)^{-1} \bZ(\bdelta)
       \begin{bmatrix} \bB_{1} & \bB_{2} \end{bmatrix}
    $$
    (so $\bC_{2} (I - \bZ(\bdelta) \bA)^{-1} \bZ(\bdelta) \bB = 0$ for
    all $\bdelta$) where
    $$
    \begin{bmatrix} \bA & \bB_{1} & \bB_{2} \\ \bC_{1} & \bD_{11} &
    \bD_{12} \\ \bC_{2} & \bD_{21} & 0 \end{bmatrix} =
    \begin{bmatrix} A & B_{1} & B_{2} \\ C_{1} & D_{11} & D_{12}
        \\ C_{2} & D_{21} & 0 \end{bmatrix} \otimes I_{\ell^{2}}
     $$
     as usual.  Then there exists a $K \in H^{\infty,o}_{nc, d+1}$ so
     that $T_{1} + T_{2} K T_{3}$ is in the strict noncommutative
     Schur-Agler class $\mathcal{S A}^{o}_{nc, d+1}$ if and only if
there exist $X,Y\in\cD$, with $\cD$ as in \eqref{Zcom}, satisfying LMIs:
\begin{align*}
       & \begin{bmatrix} N_{c} & 0 \\ 0 & I \end{bmatrix} ^{*}
            \begin{bmatrix} AYA^{*} - Y & AYC_{1}^{*} & B_{1} \\
     C_{1}Y A^{*} & C_{1} Y C_{1}^{*} - I & D_{11} \\
       B_{1}^{*} & D_{11}^{*} & -I \end{bmatrix}
     \begin{bmatrix} N_{c} & 0 \\ 0 & I \end{bmatrix} <  0,\quad Y>0, \\
         & \begin{bmatrix} N_{o} & 0 \\ 0 & I \end{bmatrix} ^{*}
        \begin{bmatrix} A^{*} X A - X & A^{*} X B_{1} & C_{1} ^{*} \\
       B_{1}^{*} X A & B_{1}^{*} X B_{1} - I & D_{11}^{*} \\
       C_{1} & D_{11} & -I \end{bmatrix}
       \begin{bmatrix} N_{o} & 0 \\ 0 & I \end{bmatrix} < 0,\quad X>0,
      \end{align*}
and the coupling condition
\begin{equation*}
\mat{cc}{X&I\\I&Y}\geq 0.
\end{equation*}
Here $N_{c}$ and $N_{o}$ are matrices chosen so that
    \begin{align*}
        & N_{c} \text{ is injective and }
     \operatorname{Im } N_{c} = \operatorname{Ker } \begin{bmatrix}
     B_{2}^{*} & D_{12}^{*} \end{bmatrix} \text{ and } \\
     & N_{o} \text{ is injective and } \operatorname{Im } N_{o} =
     \operatorname{Ker } \begin{bmatrix} C_{2} & D_{21} \end{bmatrix}.
     \end{align*}
     \end{theorem}

    \subsection{Notes}  \label{S:NC-notes}

    \textbf{1.} The equality of $\mu_{\boldsymbol \Delta}(\bA)$ with
    $\widehat \mu_{\Delta}(A)$ where $\boldsymbol \Delta$ is as in
    \eqref{NC-DS} appears in Paganini's thesis \cite{Paganini}; as
    mentioned in the Introduction, results of the same flavor have
    been given in \cite{BFKT, BFT5, FM00, MT, Shamma}.
    Ball-Groenewald-Malakorn \cite{BGM3} show how this result is
    closely related to the realization theory for the noncommutative
    Schur-Agler class obtained in \cite{BGM2}. There it is shown that
    $\mu_{\boldsymbol \Delta}(\bA) \le \overline{\mu}_{\boldsymbol
    \Delta}(\bA) = \widehat \mu_{\Delta}(A)$, where
    $\overline{\mu}_{\boldsymbol \Delta}(\bA)$ is a uniform version of
    $\mu_{\boldsymbol \Delta}(\bA)$.  The fact that
    $\mu_{\boldsymbol \Delta}(\bA) = \overline{\mu}_{\boldsymbol
    \Delta}(\bA)$ is the content of Theorem B.3 in \cite{Paganini}.
    Paganini's analysis is carried out in the more general
    form required to obtain the result of Proposition \ref{P:NCdetect}.

    The thesis of Paganini also includes some alternate
    versions of Proposition \ref{P:NCdetect}.  Specifically, rather
    than letting each $\bdelta_{j}$ be an arbitrary operator on
    $\ell^{2}$, one may restrict to such operators which are {\em
    causal} (i.e., lower-triangular) and/or {\em slowly
    time-varying} in a precise quantitative sense.  With any
    combination of these refined
    uncertainty structures in force, all the results developed in
    Section \ref{S:NC} continue to hold.  With one or more of these
    modifications in force, it is more plausible to argue that the
    assumption made in Section \ref{S:NC-statespace} that the
    controller $K$ has on-line access to the uncertainties
    $\bdelta_{i}$ is physically realistic.

    The replacement of the condition $\mu(\Delta) < 1$ by $\widehat
    \mu(\Delta) < 1$ can be considered as a {\em relaxation} of the
    problem:  while one really wants $\mu(\Delta) < 1$, one is
    content to analyze $\widehat \mu(\Delta) < 1$ since $\widehat
    \mu(\Delta) $ is easier to compute.  Necessary and sufficient
    conditions for $\widehat \mu(\Delta) < 1$ then provide sufficient
    conditions for $\mu(\Delta) < 1$ (due to the general inequality
    $\mu(\Delta) \le \widehat \mu(\Delta)$).  In the setting of the
    enhanced uncertainty structure discussed in this section, by the
    discussion immediately preceding  Proposition \ref{P:NCdetect} we
    see in this case that the relaxation is {\em exact} in the sense
    that $\widehat \mu(\Delta) < 1$ is necessary as well as
    sufficient for $\mu(\Delta) < 1$.  In Remark 1.2 of the paper of
    Megretsky-Treil \cite{MT}, it is shown how the
    $\mu$-singular-value approach can be put in the following general
    framework involving quadratic constraints (called the {\em
    S-procedure} for obscure reasons).  One is given
    quadratic functionals $\sigma_{0}, \sigma_{1}, \dots,
    \sigma_{\ell}$ defined on some set $L$ and one wants to know
    when it is the case that
    \begin{equation}  \label{want}
   \sigma_{j}(x) \ge 0 \text{ for } j = 1, \dots,
    \ell  \Longrightarrow \sigma_{0}(x) \le 0 \text{ for } x \in L.
    \end{equation}
    A computable sufficient condition (the relaxation) is the existence
    of nonnegative real numbers $\tau_{1}, \dots, \tau_{\ell}$
    ($\tau_{j} \ge 0$ for $j = 1, \dots, \ell$) so that
    \begin{equation}  \label{relax}
	\sigma_{0}(x) + \sum_{j=1}^{\ell} \tau_{j} \sigma_{j}(x) \le
	0 \text{ for all } x \in L.
  \end{equation}
  The main result of \cite{MT} is that there is a particular case of
  this setting (where $L$ is a linear shift-invariant subspace of
  vector-valued $L^{2}(0, \infty)$ (or more generally
  $L_{\text{loc}}^{2}(0, \infty)$) and the quadratic constraints are
  shift-invariant) where the relaxation is again exact (i.e., where
  \eqref{want} and \eqref{relax} are equivalent); this result
  is closely related to Proposition \ref{P:NCdetect} and the work of
  \cite{Paganini}.  A nice survey of the S-procedure and its
  applications to a variety of other problems is the paper of
  P\'olik-Terlaky \cite{PT}.

    \textbf{2.} It is of interest to note that the type of
    noncommutative system theory developed in this section (in
    particular, nc-detectability/stabilizability and nc-coprime
    representation as in \eqref{NCcoprime}) has been used in the work
    of Beck \cite{BeckSCL06} and Li-Paganini \cite{LiPaganini} in
    connection with model reduction for linear systems with
    LFT-modelled structured uncertainty.

    \textbf{3.}  We note that Theorem \ref{T:NCstab4} gives a
    Youla-Ku\v cera-type parametrization for the set of stabilizing
    controllers for a given plant $G \in \cR \cO^{0}_{nc,d}(\cL(\cW
    \oplus \cU, \cZ \oplus \cY))$ under
    the assumption that $G_{22}$ has a double coprime factorization.  In
    connection with this result, we formulate a noncommutative
    analogue of the conjecture of Lin:  {\em If $G \in \cR
    \cO^{0}_{nc, d}(\cL(\cW \oplus \cU, \cZ \oplus \cY))$ is
    stabilizable, does it follow that $G_{22}$ has a double-coprime
    factorization?}  If $G_{22}$ has a realization
    $$
     G_{22}(\bdelta) = \bC_{2}(I - Z(\bdelta) \bA)^{-1} Z(\bdelta)
     \bB_{2}
    $$
    with $\sbm {\bA & \bB \\ \bC & 0} = \sbm{ A & B \\ C & 0}
    \otimes I_{\ell^{2}}$ nc-Hautus stabilizable and nc-Hautus
    detectable, then one can adapt the state-space formulas for the
    classical case (see \cite{NJB, KharSon}) to arrive at state-space
    realization formulas for a double-coprime factorization of
    $G_{22}$.  If it is the case that one can always find a nc-Hautus
    stabilizable/detectable realization for $G_{22}$, it follows that
    $G_{22}$ in fact always has a double-coprime factorization and
    hence the noncommutative Lin conjecture is answered in the affirmative.
    However, we do not know at this time whether nc-Hautus
    stabilizable/detectable realizations always exist for a given
    $G_{22} \in \cR \cO^{0}_{nc, d}(\cL(\cU, \cY))$.  From the
    results of \cite{BGM1}, it is known that {\em minimal} i.e., {\em
    controllable} and {\em observable} realizations exist for a given
    $G_{22}$.  However, here controllable is in the sense that a
    certain finite collection of control operators be
    surjective and observable is in the sense that a certain finite
    collection of observation operators be injective.  It is not
    known if this type of controllability is equivalent to {\em
    nc-Hautus controllability}, i.e., to the operator pencil
    $\begin{bmatrix} I - Z(\bdelta) \bA  & \bB \end{bmatrix}$ being
    surjective for all $\bdelta \in \cL(\ell^{2})^{d+1}$ (not just
    $\bdelta$ in the noncommutative polydisk $\overline{\mathbb
    D}_{nc, d}$).  Thus it is unknown if controllable implies
    nc-Hautus stabilizable in this context.  Dually, we
    do not know if observable implies nc-Hautus detectable.

    \textbf{4.} Theorem \ref{T:noncom-sing} can be viewed as saying
    that, under a stabilizability/detect\-ability hypothesis, any stable
    singularity of the noncommutative function $W$ must show up
    internally as a singularity in the resolvent $(I - \bZ(\bdelta)
    \bA)^{-1}$ of the state matrix $\bA$.  A variant on this theme is
    the well known fact for the classical case that, under a
    controllability/observability assumption, any singularity
    (stable or not) of the rational matrix function $W(\lambda) = D +
    \lambda C (I - \lambda A)^{-1} B$ necessarily must show up internally as a
    singularity in the resolvent $(I - \lambda A)^{-1}$ of the state
    matrix $A$.  A  version of this result for the noncommutative
    case has now appeared in the paper of
    Kaliuzhnyi-Verbovetskyi--Vinnikov \cite{KVV09};  however the
    notion of controllable and observable there is not quite the same
    as the notion of controllable and observable for non-commutative
    Givone-Roesser systems as given in \cite{BGM1}.

    \textbf{5.}
      Given a function $S(z) = \sum_{n \in {\mathbb Z}^{d}_{+}} S_{n}
      z^{n}$ (where $z = (z_{1}, \dots, z_{d})$ is the variable in
      the commutative polydisk $\overline{\mathbb D}^{d}$ and we use
      the standard multivariable notation $z^{n} = z_{1}^{n_{1}}
      \cdots z_{d}^{n_{d}}$ if $n = (n_{1}, \dots, n_{d}) \in
      {\mathbb Z}^{d}_{+}$), we know from the results of
      \cite{Agler-Hellinger, AgMcC99, BT}  that $S$ has a contractive
      realization $S(z) = D + C (I - Z(z) A) Z(z) B$.  In light of
      the work of \cite{BGM2}, we see that any such contractive
      system matrix $\sbm{A & B \\ C & D} \colon (\oplus_{k=1}^{d}\cX_k
      \oplus \cU) \to (\oplus_{k=1}^{d} \cX_{k} \oplus \cY)$ can also
      be used to define an element ${\mathbf S}$ of the
      noncommutative Schur-Agler class $\mathcal{S A}_{nc, d}(\cU,
      \cY)$:
      $${\mathbf S}(\bdelta) = \bD + \bC(I - \bZ(\bdelta) \bA)^{-1}
      \bZ(\bdelta) \bB
      $$
      where  $\sbm{ \bA & \bB \\ \bC & \bD } = \sbm{A & B \\ C & D}
      \otimes I_{\ell^{2}}$.  Thus a choice of contractive
      realization $\{A, B, C, D\}$ for the commutative
      Schur-Agler-class function $S$ can be viewed as a choice of
      noncommutative lifting to a noncommutative Schur-Agler-class
      function ${\mathbf S}(\bdelta)$; the lifting property is that
      $$
       {\mathbf S}(z {\mathbf I}) = S(z) \otimes I_{\ell^{2}} \text{ where }
       z{\mathbf I} = (z_{1}I_{\ell^{2}}, \dots, z_{d}I_{\ell^{2}})
       \in \overline{\mathbb D}_{nc, d}
       \text{ if } z = (z_{1}, \dots, z_{d}) \in \overline{\mathbb
       D}^{d}.
     $$
     While the realization for the commutative function is highly
     non-unique, the realization for the noncommutative function is
     unique up to state-space similarity if arranged to be minimal
     (i.e., controllable and observable as in \cite{BGM1}).
     Philosophically one can say that evaluation of the function on
     the commutative polydisk ${\mathbb D}^{d}$ does not give enough
     frequencies to detect the realization; enlarging the frequency
     domain (or points of evaluation) to the noncommutative polydisk
     ${\mathbb D}^{d}_{nc, d}$ does give enough frequencies to detect
     the realization in an essentially unique way.


\begin{thebibliography}{99}



   \bibitem{Agler-unpublished} J.~Agler, Interpolation, unpublished
      manuscript, 1988.

           \bibitem{Agler-Hellinger} J.~Agler, On the representation of
      certain holomorphic functions defined on a polydisk, in: {\em
      Topics in Operator Theory: Ernst D. Hellinger Memorial Volume}
      (Ed. L.~de Branges, I.~Gohberg, and J.~Rovnyak) pp. 47--66,
      OT \textbf{48} Birkh\"auser, Basel-Berlin-Boston, 1990.

      \bibitem{AgMcC99} J.~Agler and J.E.~McCarthy, Nevanlinna-Pick
      interpolation on the bidisk, {\em J.~reine angew.~Math.}
      \textbf{506} (1999), 191--124.

      \bibitem{AgMcCbook}  J.~Agler and J.E.~McCarthy, {\em Pick
      Interpolation and Hilbert Function Spaces}, Graduate Studies in
      Mathematics Vol. \textbf{44}, American Mathematical Society,
      Providence, 2002.

\bibitem{AY1} J.~Agler and N.J.~Young, A commutant lifting
      theorem for a domain in ${\mathbb C}^{2}$ and spectral
      interpolation, {\em J.~Funct.~Anal.} \textbf{161} (1999)
      No.~2,  452--477.

      \bibitem{AY2} J.~Agler and N.J.~Young, Operators having the
      symmetrized bidisc as  spectral set, {\em Proc. Edinburgh
      Math.~Soc.} (2) \textbf{43} (2000) No.~1, 195--210.

      \bibitem{AY3} J.~Agler and N.J.~Young, The two-point spectral
      Nevanlinna-Pick problem, {\em Integral Equations Operator
      Theory} \textbf{37} (2000) No.~4, 375--385.

      \bibitem{AY4} J.~Agler and N.J.~Young, A Schwarz lemma for the
      symmetrized bidisc, {\em Bull.~London Math.~Soc.} \textbf{33}
      (2001) No~2, 175--186.

      \bibitem{AY5} J.~Agler and N.J.~Young, A model theory for
      $\Gamma$-contractions, {\em J.~Operator Theory} \textbf{49}
      (2003) No.~1, 45--60.

      \bibitem{AY6} J.~Agler and N.J.~Young, Realization of functions
      into the symmetrised bidisc,  in: {\em Reproducing Kernel Spaces and
      Applications}, pp. 1--37, OT \textbf{143}, Birkh\"auser,
      Basel-Berlin-Boston, 2003.

      \bibitem{AY7} J.~Agler and N.J.~Young, The two-by-two spectral
      Nevanlinna-Pick problem, {\em Trans.~Amer. Math.~Soc.}
      \textbf{356} (2004) No.~2, 573--585.

      \bibitem{AY8} J.~Agler and N.J.~Young, The hyperbolic geometry
      of the symmetrized bidisc, {\em J.~Geomet.~Anal.}
      \textbf{14} (2004) No.~3, 375--403.

      \bibitem{AY9} J.~Agler and N.J.~Young, The complex geodesics of
      the symmetrized bidisc, {\em Internat.~J.~Math.}
      \textbf{17} (2006) No.~4, 375--391.

  \bibitem{A-KV} D.~Alpay and D.S.~Kalyuzhny\u{\i}-Verbovetzki\u{\i},
  On the intersection of null spaces for matrix substitutions in a
  non-commutative rational formal power series, {\em
  C.R.~Acad. Sci.~Paris} Ser. I \textbf{339} (2004), 533--538.

  \bibitem{AT} C.-G.~Ambrozie and D.~Timotin, A von Neumann type
  inequality for certain domains in ${\mathbf C}^{n}$, {\em
  Proc~Amer.~Math.~Soc.} \textbf{131} (2003) No.~3, 859--869.

  \bibitem{Anderson}  B.D.O.~Anderson, P.~Agathoklis, E.I.~Jury and
  M.~Mansour, Stability and the matrix Lyapunov equation for discrete
  2-dimensional systems, {\em IEEE Trans.~Circuits \& Systems}
  \textbf{33} (1986) No.~3, 261--267.

\bibitem{Ando} T. And\^o, On a pair of commutative contractions,
{\em Acta Sci.~Math.} {\bf 24} (1963), 88--90.

  \bibitem{AG}  P.~Apkarian and P.~Gahinet, A convex characterization of
  gain-scheduled $H^{\infty}$ controllers, {\em  IEEE Trans. Automat.~Control},
  \textbf{40} (1995) No.~5, 853--864.

  \bibitem{ArPo} A.~Arias and G.~Popescu, Noncommutative
  interpolation and Poisson transforms, {\em Israel J.~Math.}
  \textbf{115} (2000), 205--234.

  \bibitem{BallBoloJFA} J.A.~Ball and V.~Bolotnikov, Realization and
  interpolation for Schur-Agler-class functions on domains with
  matrix polynomial defining function in ${\mathbb C}^{n}$, {\em
  J.~Funct.~Anal.} \textbf{213} (2004), 45--87.

  \bibitem{BallBoloNYJ} J.A.~Ball and V.~Bolotnikov, Nevanlinna-Pick
  interpolation for Schur-Agler class functions on domains with
  matrix polynomial defining function, {\em New York J.~Math.}
  \textbf{11} (2005), 245--209.

  \bibitem{BallBoloJOT}
  J.A.~Ball and V.~Bolotnikov, Interpolation in the noncommutative
  Schur-Agler class, {\em J. Operator Theory} \textbf{58} (2007)
  No.~1, 83--126.

  \bibitem{BallChu} J.A.~Ball, J.~Chudoung, and M.V.~Day,  Robust optimal
  switching control for nonlinear systems, {\em SIAM J.~Control
  Optim.} \textbf{41} (2002) No.~3, 900--931.

  \bibitem{BallCohen} J.A.~Ball and N.~Cohen, Sensitivity minimization
  in an $H_{\infty}$ norm: Parametrization of all solutions, {\em
  Internat.~J.~Control} \textbf{46} (1987), 785--816.

\bibitem{BFGtH09}
 J.A.~Ball, Q.~Fang, G.~Groenewald, and S.~ter Horst,
Equivalence of robust stabilization and robust performance via
feedback, {\em Math.~Control Signals Systems} {\bf 21} (2009), 51--68.

  \bibitem{BGR} J.A.~Ball, I.~Gohberg, and L.~Rodman, {\em
  Interpolation of Rational Matrix Functions}, OT \textbf{44},
  Birkh\"auser, Basel-Berlin-Boston, 1990.

 \bibitem{BGM1} J.A.~Ball, G.~Groenewald and T.~Malakorn,
  Structured noncommutative multidimensional linear systems,
  {\em SIAM J.~Control Optim.} \textbf{44} (2005) No.~4,
  1474--1528.

  \bibitem{BGM2} J.A.~Ball, G.~Groenewald and T.~Malakorn, Conservative
  structured noncommutative multidimensional linear systems, in: {\em
  The State Space Method Generalizations and Applications} (D.~Alpay
  and I.~Gohberg, ed.), pp. 179--223, OT \textbf{161}, Birkh\"auser,
  Basel-Berlin-Boston, 2005.

  \bibitem{BGM3} J.A.~Ball, G.~Groenewald and T.~Malakorn, Bounded real
  lemma for structured noncommutative multidimensional linear systems
  and robust control, {\em Multidimens.~Sys.~Signal
  Process.} \textbf{17} (2006), 119--150.

    \bibitem{BtH} J.A.~Ball and S.~ter Horst, Multivariable
  operator-valued Nevanlinna-Pick interpolation: a survey, {\em Proceedings
  of IWOTA (International Workshop on Operator Theory and
  Applications)} 2007, Potchefstroom, South Africa, Birkh\"auser,
  volume to appear.

  \bibitem{BLTT} J.A.~Ball, W.S.~Li, D.~Timotin and T.T.~Trent, A commutant
  lifting theorem on the polydisc: interpolation problems for the
  bidisc, {\em Indiana Univ. Math. J.} \textbf{48} (1999), 653-675.

  \bibitem{BM} J.A.~Ball and T.~Malakorn, Multidimensional linear
  feedback control systems and interpolation problems for multivariable
  holomorphic functions, {\em Multidimens.~Sys.~Signal
  Process.} \textbf{15} (2004), 7--36.

  \bibitem{BR} J.A.~Ball and A.C.M.~Ran, Optimal Hankel norm model
  reductions and Wiener-Hopf factorization I: The canonical case,
  {\em SIAM J.~Control Optim.} \textbf{25} (1987) No.~2,
  362--382.

  \bibitem{BSV} J.A.~Ball, C.~Sadosky, and V.~Vinnikov, Scattering
  systems with several evolutions and multidimensional
  input/state/output linear systems, {\em Integral Equations
  Operator Theory} \textbf{52} (2005), 323--393.

  \bibitem{BT} J.A.~Ball and T.T.~Trent, Unitary colligations,
  reproducing kernel Hilbert spaces, and Nevanlinna-Pick interpolation
  in several variables, {\em J.~Funct.~Anal.} \textbf{157}
  (1998), 1--61.

  \bibitem{BeckSCL06} C.L.~Beck, Coprime factors reduction methods for
  linear parameter varying and uncertain systems, {\em Systems
  Control Lett.} \textbf{55} (2006), 199--213.


  \bibitem{BFKT} H.~Bercovici, C.~Foias, P.P.~Khargonekar, and A.~Tannenbaum,
  On a lifting theorem for the structured singular value, {\em
  J.~Math.~ Anal.~Appl.} \textbf{187} (1994),
  617--627.

 \bibitem{BFT1}  H.~Bercovici, C.~Foias, and A.~Tannenbaum,
  Structured interpolation theory, in: {\em Extensions and
  Interpolation of Linear Operators and Matrix Functions} pp.~195--220,
  OT \textbf{47}, Birkh\"auser, Basel-Berlin-Boston, 1990.

  \bibitem{BFT2} H.~Bercovici, C.~Foias, and A.~Tannenbaum, A spectral
  commutant lifting theorem, {\em Trans.~Amer.~Math.~Soc.}
  \textbf{325} (1991) No.~2, 741--763.

  \bibitem{BFT3} H.~Bercovici, C.~Foias, and A.~Tannenbaum, On
  spectral tangential Nevanlinna-Pick interpolation, {\em
  J.~Math.~Anal.~Appl.} \textbf{155} (1991) No.~1, 156--176.

  \bibitem{BFT4} H.~Bercovici, C.~Foias, and A.~Tannenbaum,
  On the optimal solutions in spectral commutant lifting theory, {\em
  J.~Funct.~Anal.} \textbf{101} (1991) No.~1, 38--49.

  \bibitem{BFT5} H.~Bercovici, C.~Foias, and A.~Tannenbaum, The
  structured singular value for linear input/output operators, {\em
  SIAM J.~Control Optim.} \textbf{34} (1996) No.~4, 1392--1404.

  \bibitem{BD06} V. Bolotnikov and H. Dym, On Boundary Interpolation for Matrix
  Valued Schur Functions, {\em Mem.~Amer.~Math.~Soc.} {\bf 181} (2006), no. 856.

  \bibitem{Bognar} J.~Bogn\'ar, {\em Indefinite Inner Product Spaces},
  Springer-Verlag, New York-Heidelberg-Berlin, 1974.

  \bibitem{Bose} N.K.~Bose, Problems and progress in multidimensional
  systems theory, {\em Proc.~IEEE} \textbf{65} (1977) No.~6, 824--840.

  \bibitem{BST} C.I.~Byrnes, M.W.~Spong, and T.-J.~Tarn, A several
  complex variables approach to feedback stabilization of linear
  neutral delay-differential systems, {\em Math.~Systems Theory}
  \textbf{17} (1984), 97--133.

  \bibitem{ChenFran} T.~Chen and B.A.~Francis, {\em Optimal
  Sampled-Data Control Systems}, Springer-Verlag, London, 1996.

  \bibitem{CZ} R.F.~Curtain and H.J.~Zwart, {\em An Introduction to
  Infinite-Dimensional Linear Systems Theory}, Texts in Applied
  Mathematics \textbf{21}, Springer-Verlag, Berlin, 1995.

  \bibitem{DavP} K.R.~Davidson and D.R.~Pitts, Nevanlinna-Pick
  interpolation for noncommutative analytic Toeplitz algebras,
  {\em Integral Equations and Operator Theory} \textbf{31} (1998)
  No.~3, 321--337.

  \bibitem{DLMS} C.A.~Desoer, R.-W.~Liu, and R.~Saeks, Feedback system
  design:  The fractional approach to analysis and synthesis, {\em IEEE
  Trans.~Automat.~Control} \textbf{25} (1980) No.~3, 399--412.


\bibitem{D66}
R.G. Douglas, On majorization, factorization, and range inclusion of operators on
Hilbert space,
{\em Proc.~Amer.~Math.~Soc.} {\bf 17} (1966), 413--415.

\bibitem{DoyleMu} J.C.~Doyle, Analysis of feedback systems with
structured uncertainties, {\em IEE Proceedings} \textbf{129} (1982),
242--250.

\bibitem{Doyle-notes} J.C.~Doyle, {\em Lecture notes in advanced
multivariable control},  ONR/Honeywell Workshop, Minneapolis, 1984.

\bibitem{DGKF}  J.C.~Doyle, K.~Glover, P.P.~Khargonekar, and
B.A.~Francis, State-space solutions to standard $H_{2}$ and
$H_{\infty}$ control problems, {\em IEEE Trans.~Automat.~Control}
\textbf{34} (1989), 831--847.

  \bibitem{DuXie} C.~Du and L.~Xie,  {\em $H_{\infty}$ Control and
  Filtering of Two-dimensional Systems}, Lecture Notes in Control and
  Information Sciences \textbf{278}, Springer, Berlin, 2002.

  \bibitem{DXZ} C.~Du, L.~Xie and C.~Zhang, $H_{\infty}$ control and
  robust stabilization of two-dimensional systems in Roesser models,
  {\em Automatica} \textbf{37} (2001), 205--211.

  \bibitem{DP}
  G.E.~Dullerud and F.~Paganini, {\em A Course in Robust Control
  Theory: A Convex Approach}, Texts in Applied Mathematics Vol.~{\bf 36},
  Springer-Verlag, New York, 2000.

  \bibitem{Dym}  H.~Dym, {\em $J$ Contractive Matrix Functions,
  Reproducing Kernel Hilbert Spaces and Interpolation}, CBMS No.
  \textbf{71}, American Mathematical Society, Providence, 1989.

\bibitem{ElAgiziFahmy}
N.G. El-Agizi, M.M. Fahm\.y, Two-dimensional digital filters with no
overflow oscillations, {\em IEEE Trans. Acoustical. Speech Signal
Process.} {\bf 27} (1979), 465--469.

 \bibitem{FM00} A. Feintuch and A. Markus, The structured norm of
 a Hilbert space operator with respect to a given algebra
 of operators, in: {\em  Operator Theory and Interpolation}, pp.~163--183,
\textbf{OT 115}, Birkh\"auser-Verlag, Basel-Berlin-Boston, 2000.

 \bibitem{Finsler} P. Finsler, \"Uber das volkommen definiter und semidefiniter
Formen in Scharen quadratischer Formen, Comment.~Math.~Helv.~{\bf 9} (1937),
188--192.

  \bibitem{FF}  C.~Foias and A.E.~Frazho, {\em The Commutant Lifting
  Approach to Interpolation Problems}, OT \textbf{44},
  Birkh\"auser-Verlag, Basel-Berlin-Boston, 1990.

  \bibitem{FFGK} C.~Foias, A.E.~Frazho, I.~Gohberg, and
  M.A.~Kaashoek, {\em Metric Constrained Interpolation, Commutant
  Lifting and Systems}, OT \textbf{100}, Birkh\"auser-Verlag,
  Basel-Berlin-Boston, 1998.

  \bibitem{Francis}
  B.A.~Francis, {\em A Course in $H_{\infty}$ Control
  Theory}, Lecture Notes in Control and
  Information Sciences \textbf{88}, Springer, Berlin, 1987.

  \bibitem{FHZ}
  B.A.~Francis, J.W.~Helton, and G.~Zames,
  $H^{\infty}$-optimal feedback controllers for linear multivariable
  systems,  {\em IEEE Trans.~Automat.~Control} \textbf{29} (1984)
  No.~10, 888--900.


  \bibitem{GA}
  P.~Gahinet and P.~Apkarian, A linear matrix inequality
  approach to $H^{\infty}$ control, {\em Internat.~J. of
  Robust Nonlinear Control} \textbf{4} (1994), 421-448.

  \bibitem{GR}
  D.D.~Givone and R.P.~Roesser, Multidimensional linear
  iterative circuits---General properties, {\em IEEE Trans.~Compt.},
  \textbf{21} (1972) , 1067--1073.

  \bibitem{GN}
  L.~El Ghaoui and S.-I.~Niculescu (editors), {\em
  Advances in Linear Matrix Inequality Methods in Control}, SIAM,
  Philadelphia, 2000.

  \bibitem{Glover} K.~Glover, All optimal Hankel-norm approximations
  of linear multivariable systems and their $L_{\infty}$-error
  bounds, {\em Int.~J.~Control} \textbf{39} (1984) No.~6, 1115--1193.

  \bibitem{Green} M.~Green, $H_{\infty}$ controller synthesis by
  $J$-lossless coprime factorization, {\em SIAM J.~Control
  Optim.} \textbf{28} (1992), 522--547.

  \bibitem{GGLD} M.~Green, K.~Glover, D.J.N.~Limebeer, and
  J.C.~Doyle, A $J$-spectral factorization approach to
  $H_{\infty}$-control, {\em SIAM J.~Control  Optim.}
  \textbf{28} (1990), 1350--1371.

  \bibitem{GL} M.~Green and D.J.N.~Limebeer, {\em Linear Robust
  Control}, Prentice Hall, London, 1995.


  \bibitem{HeltonACC} J.W.~Helton, A type of gain scheduling which converts to
  a ``classical'' problem in several complex variables, {\em
  Proc.~Amer.~Control Conf.} 1999, San Diego, CA.

  \bibitem{HeltonIEEE}
  J.W.~Helton, Some adaptive control problems
  which convert to a ``classical'' problem in several complex
  variables, {\em IEEE Trans.~Automat.~Control} \textbf{46} (2001)
  No.~12, 2038--2043.

  \bibitem{HMcCV}  J.W.~Helton, S.A.~McCullough and V.~Vinnikov,
  Noncommutative convexity arises from Linear Matrix Inequalities,
  {\em J.~Funct.~Anal.} \textbf{240} (2006), 105--191.

  \bibitem{HinPrit} D.~Hinrichsen and A.J.~Pritchard, Stochastic
  $H^{\infty}$, {\em SIAM J.~Control Optim.} \textbf{36} (1998)
  No.~5, 1504--1538.

  \bibitem{HMY} H.-N.~Huang, S.A.M.~Marcantognini and N.J.~Young, The
  spectral Carath\'eodory-Fej\'er problem, {\em Integral Equations
  Operator Theory} \textbf{56} (2006) No.~2, 229--256.

  \bibitem{IS}
  T.~Iwasaki and R.E.~Skelton, All controllers for the
  general $H_{\infty}$ control problem:  LMI existence conditions and
  state space formulas, {\em Automatica} \textbf{30} (1994) No.~8,
  1307--1317.

  \bibitem{James}  M.R.~James, H.I.~Nurdin, and I.R.~Petersen,
  $H^{\infty}$ control of linear quantum stochastic systems, {\em IEEE
  Trans.~Automat.~Control} \textbf{53} (2008) No.~8, 1787--1803.

\bibitem{Jury}
E.I. Jury, Stability of multidimensional scalar and matrix polynomials,
Proc.~IEEE, vol.~{\bf 66} (1978), 1018--1047.

  \bibitem{Kaczorek}
  T.~Kaczorek, {\em Two-Dimensional Linear Systems},
  Lecture Notes in Control and Information Sciences \textbf{68},
  Springer-Verlag, Berlin, 1985.

  \bibitem{KVV09}
  D.S. Kaliuzhnyi-Verbovetskyi and V. Vinnikov,
  Singularities of rational functions and minimal factorizations: The noncommutative
  and commutative setting, {\em Linear Algebra Appl.} {\bf 430} (2009),
  869--889.

  \bibitem{KKT} E.W.~Kamen, P.P.~Khargonekar and A.~Tannenbaum,
  Pointwise stability and feedback control of linear systems with
  noncommensurate time delays, {\em Acta Appl.~Math.}
  \textbf{2} (1984), 159--184.

  \bibitem{KharMun}
  V.L.~Kharitonov and J.A.~Torres-Mu\~noz, Robust
  stability of multivariate polynomials. Part 1:  small coefficient
  perturbations, {\em Multidimens.~Sys.~Signal Process.}
  \textbf{10} (1999), 7--20.

  \bibitem{KharSon}
  P.P.~Khargonekar and E.D.~Sontag, On the relation
  between stable matrix fraction factorizations and regulable
  realizations of linear systems over rings, {\em IEEE
  Trans.~Automat.~Control} \textbf{27} (1982) No.~3, 627--638.

  \bibitem{Kimura87} H.~Kimura, Directional interpolation approach to
  $H_{\infty}$-optimization and robust stabilization, {\em IEEE
  Trans.~Automat.~Control} \textbf{32} (1987), 1085--1093.

  \bibitem{Kimura89} H.~Kimura, Conjugation, interpolation and
  model-matching in $H^{\infty}$, {\em Int.~J.~Control} \textbf{49}
  (1989), 269--307.

  \bibitem{2D-II}
  S.Y.~Kung, B.C.~L\'evy, M.~Morf and T.~Kailath, New
  results in 2-D systems theory, Part II: 2-D state-space
  models---realization and the notions of controllability,
  observability, and minimality, {\em Proceedings of the IEEE}
  \textbf{65} (1977) No.~6, 945--961.


  \bibitem{LiPaganini} L.~Li and F.~Paganini, Structured coprime factor
  model reduction based on LMIs, {\em Automatica} \textbf{41} (2005)
  No.~1, 145--151.

  \bibitem{LA} D.J.N.~Limebeer and B.D.O.~Anderson, An
  interpolation theory approach to $H_{\infty}$ controller degree
  bounds, {\em Linear Algebra Appl.} \textbf{98}
  (1988), 347--386.

  \bibitem{LH} D.J.N.~Limebeer and G.~Halikias, An analysis of pole
  zero cancellations in $H_{\infty}$ control problems of the second
  kind, {\em SIAM J.~Control Optim.} \textbf{25} (1987),
  1457--1493.

  \bibitem{Lin1} Z.~Lin, Feedback stabilization of MIMO $n$-D linear systems,
  {\em Multidimens.~Sys.~Signal Process.} \textbf{9}
  (1998), 149--172.

  \bibitem{Lin2} Z.~Lin, Feedback stabilization of MIMO 3-D linear
  systems, {\em IEEE Trans.~Automat.~Control} \textbf{44} (1999),
  1950--1955.

\bibitem{Lin3} Z. Lin, Output Feedback Stabilizability and Stabilization
of Linear $n$D Systems, In: <{\em Multidimensional Signals, Circuits and Systems},
(J. Wood and K. Galkowski eds.), pp.~59-76, Chapter 4, Taylor \& Francis, London, 2001.

 \bibitem{LF81} J.H. Lodge and M.M. Fahmy, Stability and overflow
oscillations in 2-D state-space digital filters, IEEE Trans.~Acoustical.~Speech
Signal Processing, vol. ASSP-{\bf 29} (1981), 1161--1171.

  \bibitem{LuThesis}
  W.-M.~Lu, {\em Control of Uncertain Systems:
  State-Space Characterizations}, Thesis submitted to California
  Institute of Technology, Pasadena, 1995.


  \bibitem{LZDProc} W.-M.~Lu, K.~Zhou and J.C.~Doyle, Stabilization of
  $LFT$ systems, {\em Proc.~30th Conference on Decision and Control},
  Brighton, England, December 1991, 1239--1244.

  \bibitem{LZD}  W.-M.~Lu, K.~Zhou and J.C.~Doyle, Stabilization
  of uncertain linear systems:  An LFT approach, {\em IEEE
  Trans.~Auto.~Contr.} \textbf{41} (1996) No.~1 , 50-65.

  \bibitem{MT} A.~Megretsky and S.~Treil, Power distribution
  inequalities in optimization and robustness of uncertain systems,
  {\em J.~Mathematical Systems, Estimation, and Control}
  \textbf{3} (1993) No.~3, 301--319.

  \bibitem{McFGl} D.C.~McFarlane and K.~Glover, {\em Robust
  Controller Design Using Normalized Coprime Factor Plant
  Descriptions}, Lecture Notes in Control and Information Sciences
  \textbf{138}, Springer-Verlag, Berlin-New York, 1990.

  \bibitem{2D-I} M.~Morf, B.C.~L\'evy, and S.-Y.Kung, New results in
  2-D systems theory, Part I: 2-D polynomial matrices, factorization,
  and coprimeness, {\em Proceedings of the IEEE} \textbf{65} (1977)
  No.~6, 861--872.

  \bibitem{Mori02} K.~Mori, Parameterization of stabilizing controllers
  over commutative rings with application to multidimensional systems,
  {\em IEEE Trans.~Circuits and Systems---I} \textbf{49} (2002) No.~6,
  743--752.

  \bibitem{Mori04} K.~Mori, Relationship between standard control
  problem and model-matching problem without coprime factorizability,
  {\em IEEE Trans.~Automat.~Control} \textbf{49} (2004) No.~2, 230--233.

  \bibitem{NJB}  C.N.~Nett, C.A.~Jacobson, and M.J.~Balas, A
  connection between state-space and doubly coprime fractional
  representations, {\em IEEE Trans.~Automat.~Control} \textbf{29}
  (1984) No.~9, 831--832.

  \bibitem{NevPick} R.~Nevanlinna, \"Uber beschr\"ankte Funktionen, die in gegebene
  Punkten vor\-ge\-schrie\-bene Werte annehmen, {\em Ann.
  ~Acad.~Sci.~Fenn.~Ser.~A} \textbf{13} (1919) No.~1.

  \bibitem{Packard} A.~Packard, Gain scheduling via linear fractional
  transformations, {\em Systems \& Control Letters} \textbf{22}
  (1994), 79--92.

  \bibitem{PackardDoyle}  A.~Packard and J.C.~Doyle, The complex
  structured singular value, {\em Automatica} \textbf{29} (1993) No.~1,
  71--109.

  \bibitem{Paganini} F.~Paganini, {\em Sets and Constraints in the
  Analysis of Uncertain Systems}, Thesis submitted to California
  Institute of Technology, Pasadena, 1996.

  \bibitem{Paulsen}  V.~Paulsen, {\em Completely Bounded Maps and
  Operator Algebras}, Cambridge Studies in Advanced Mathematics
  \textbf{78}, 2002.

  \bibitem{Pick}  G.~Pick, \"Uber die beschr\"ankungen analytischer
  Funktionen, welche durch vor\-ge\-ge\-bene Funktionswerte bewirkt werden,
  {\em Math.~Ann.} \textbf{7} (1916), 7--23.

\bibitem{Piekarski} M.S. Piekarski, Algebraic characterization of matrices
whose multivariable characteristic polynomial is Hurwitzian, in:
{\em Proc.~Int.~Symp.~Operator Theory} Lubbock, TX, Aug. 1977, 121--126.

\bibitem{PT}  I.~P\'olik and T.~Terlaky, A survey of the S-lemma,
{\em SIAM Review} \textbf{49} (2007) No.~3, 371--418.

\bibitem{Po1} G.~Popescu, Interpolation problems in several
variables, {\em J.~Math.~Anal.~Appl.} \textbf{227} (1998) NO.~1,
227--250.

\bibitem{Po2} G.~Popescu, Spectral lifting in Banach algebras and
interpolation in several variables, {\em Trans.~Amer.~Math.~So.}
\textbf{353} (2001) No.~7, 2843--2857.

\bibitem{Po3} G.~Popescu, Free holomorphic functions on the unit ball
of $B({\mathcal H})^{n}$, {\em J.~Funct.~Anal.} \textbf{241}
(2006) No.~1, 268--333.

\bibitem{Po4} G.~Popescu, Noncommutative transforms and free
pluriharmonic functions, {\em Advances in Mathematics} \textbf{220}
(2009), 831--893.


  \bibitem{Q-LN}
  A.~Quadrat, An introduction to internal stabilization
  of infinite-dimensional linear systems, {\em Lecture notes of the
   International School in Automatic Control of Lille:  Control of
  Distributed Parameter Systems:  Theory \& Applications} (organized by
  M.~Fliess \& W. Perruquetti), Lille (France) September 2--6, 2002.

  \bibitem{Q-YKI}
  A.~Quadrat, On a generalization of the
  Youla-Ku\v{c}era parametrization. Part I: The fractional ideal
  approach to SISO systems, {\em Systems Control Lett.}
  \textbf{50} (2003) No~2, 135--148.

  \bibitem{Q-Leuven}  A.~Quadrat, Every internally stabilizable
  multidimensional system admits a doubly coprime factorization, {\em
  Proceedings of the International Symposium on the Mathematical Theory
  of Networks and Systems}, Leuven, Belgium, July, 2004.

  \bibitem{Q-elementary}
  A.~Quadrat, An elementary proof of the general
  $Q$-parametrization of all stabilizing controllers, {\em Proc. 16th
  IFAC World Congress}, Prague (Czech Republic), July 2005.

  \bibitem{Q-Lattice}
  A.~Quadrat, A lattice approach to analysis and
  synthesis problems, {\em Math.~Control Signals Systems}
  \textbf{18} (2006) No.~2, 147--186.

  \bibitem{Q-YKII}
  A.~Quadrat, On a generalization of the
  Youla-Ku\v{c}era parametrization. Part II: The lattice approach to
  MIMO systems, {\em Math.~Control Signals Systems}
  \textbf{18} (2006) No.~3, 199--235.

  \bibitem{RogersGalOw} E.~Rogers, K.~Galkowski, and D.H.~Owens, {\em
  Control Systems Theory and Applications for Linear Repetitive
  Processes}, Lecture Notes in Control and Information Sciences
  \textbf{349}, Springer, Berlin-Heidelberg, 2007.

  \bibitem{SafonovMu} M.G.~Safonov, {\em Stability Robustness of
  Multivariable Feedback Systems}, MIT Press, Cambridge, MA, 1980.


  \bibitem{Sarason} D.~Sarason, Generalized interpolation in
  $H^{\infty}$, {\em Trans.~Amer.~Math.~Soc.} \textbf{127} (1967)
  No.~2, 179--203.

  \bibitem{vdS} A.J.~van der Schaft, {\em $L_{2}$-Gain and Passivity
  Techniques in Nonlinear Control},
  Second Edition, Springer-Verlag,
  London, 2000.


  \bibitem{Scherer} C.W.~Scherer, $H^{\infty}$-optimization without
  assumptions on finite or infinite zeros, {\em SIAM J.~Control
  Optim.} \textbf{30} (1992) No.~1, 143--166.

  \bibitem{Shabat} B.V.~Shabat, {\em Introduction to Complex Analysis
  Part II: Functions of Several Variables}, Translations of Mathematical
  Monographs vol.~\textbf{110}, American Mathematical Society, 1992.

  \bibitem{Shamma} J.S.~Shamma, Robust stability with time-varying
  structured uncertainty, {\em IEEE Trans.~Automat.~Control}
  \textbf{39} (1994) No.~4, 714--724.

  \bibitem{Smith}  M.C.~Smith, On stabilization and existence of
  coprime factorizations, {\em IEEE Trans.~Automat.~Control}
  \textbf{34} (1989), 1005--1007.

  \bibitem{SRP}
  M.N.S.~Swamy, L.M.~Roytman, and E.I.~Plotkin, On
  stability properties of three- and higher dimensional linear
  shift-invariant digital filters, {\em IEEE Trans.~Circuits and Systems}
  \textbf{32} (1985) No.~9, 888--892.

  \bibitem{Sule}
  V.R.~Sule, Feedback stabilization over commutative
  rings:  the matrix case, {\em SIAM J.~Control Optim.}
  \textbf{32} (1994) No.~6, 1675--1695.

  \bibitem{Treil}
  S. Treil, The gap between the complex structures singular value $\mu$ and its
  upper bound is infinite, preprint.

  \bibitem{TrenWil} H.L.~Trentelman and J.C.~Willems, $H_{\infty}$
  control in a behavioral context: the full information case, {\em
  IEEE Trans.~Automat.~Control} \textbf{44} (1999) No.~3, 521--536.

  \bibitem{Uhlig} F.~Uhlig, A recurring theorem about pairs of
  quadratic forms and extensions:  a survey, {\em Linear Algebra and
  its Applications} \textbf{25} (1979), 219--237.

  \bibitem{Vid} M. Vidyasagar, {\em Control System Synthesis:  A
  Factorization Approach}, MIT Press, Cambridge, 1985.

  \bibitem{VSF}
  M.~Vidyasagar, H.~Schneider and B.A.~Francis, Algebraic
  and topological aspects of feedback stabilization, {\em IEEE
  Trans.~Automat.~Control} \textbf{27} (1982) No.~4, 880--894.

  \bibitem{Youla}
  D.C.~Youla and G.~Gnavi, Notes on $n$-dimensional
  system theory, {\em IEEE Trans.~Circuits and Systems} \textbf{26}
  (1979) No.~2, 105--111.

  \bibitem{Zames}
  G.~Zames, Feedback and optimal sensitivity:  Model
  reference transformations, multiplicative seminorms, and approximate
  inverses, {\em IEEE Trans.~Automat.~Control} \textbf{26} (1981)
  No.~2, 301--320.


  \bibitem{ZF}
  G.~Zames and B.A.~Francis, Feedback, minimax
  sensitivity, and optimal robustness, {\em IEEE
  Trans.~Automat.~Control} \textbf{28} (1983) No.~5, 585--601.

  \bibitem{ZDG}
  K.~Zhou, J.C.~Doyle and K.~Glover, {\em Robust and
  Optimal Control}, Prentice-Hall, Upper Saddle River, NJ, 1996.

  \end{thebibliography}
\end{document}